\documentclass[2p]{article}
\usepackage{amssymb,amsmath,amsthm}
\usepackage{indentfirst}
\usepackage{exscale}
\usepackage{relsize}
\usepackage{xcolor}
\usepackage[numbers,sort&compress]{natbib}

\usepackage{geometry}

\theoremstyle{definition}

\allowdisplaybreaks
\theoremstyle{remark}

\numberwithin{equation}{section}

\begin{document}

\title{\Large\bf{Multiplicity of solutions for a nonhomogeneous quasilinear elliptic equation with concave-convex nonlinearities }
 }
\date{}
\author {
\ Wanting Qi$^{1}$,
\ Xingyong Zhang$^{1,2}$\footnote{Corresponding author, E-mail address: zhangxingyong1@163.com}\\
{\footnotesize $^{1}$Faculty of Science, Kunming University of Science and Technology, Kunming, Yunnan, 650500, P.R. China.}\\
{\footnotesize $^{2}$Research Center for Mathematics and Interdisciplinary Sciences, Kunming University of Science and Technology,}\\
 {\footnotesize Kunming, Yunnan, 650500, P.R. China.}\\
 }
 \date{}
 \maketitle

 \begin{center}
 \begin{minipage}{15cm}

\small  {\bf Abstract:}
We investigate the multiplicity of solutions for a quasilinear scalar field equation with a nonhomogeneous differential operator defined by
\begin{eqnarray}\label{0.1.1}
 Su:=-\mbox{div}\left\{\phi \left(\frac{u^{2}+|\nabla u|^{2}}{2}\right)\nabla u\right\}+\phi\left (\frac{u^{2}+|\nabla u|^{2}}{2}\right)u,
 \end{eqnarray}
where $\phi:[0,+\infty)\rightarrow\mathbb{R}$ is a positive continuous function.
This operator is introduced by C. A. Stuart [Milan J. Math. 79 (2011), 327-341]
and depends on not only $\nabla u$  but also $u$.
This particular quasilinear term generally appears in the study of nonlinear optics model
which describes the propagation of self-trapped beam in a cylindrical optical fiber made from a self-focusing dielectric material.
When the reaction term is concave-convex nonlinearities,
by using the Nehari manifold and doing a fine analysis associated on the fibering map,
we obtain that the equation admits at least one positive energy solution  and negative energy solution  which is also the ground state solution of the equation. We overcome two main difficulties which are caused by the nonhomogeneity of the differential operator $S$: (i) the almost everywhere convergence of the gradient for the minimizing sequence $\{u_{n}\}$; (ii) seeking the reasonable restrictions about $S$.
 \par
 {\bf Keywords:}
Quasilinear elliptic equation; concave-convex nonlinearities; Nehari manifold method; fibering maps; ground state solution; nonlinear optics.
\par
 {\bf 2010 Mathematics Subject Classification.} 35J20; 35J50; 35J55.
\end{minipage}
 \end{center}
  \allowdisplaybreaks
 \vskip2mm
\section{Introduction and main results}
 In this paper, we investigate the multiplicity of solutions for the following quasilinear elliptic problem:
\begin{equation}\label{eq1}
 \left\{
  \begin{array}{ll}
 -\mbox{div}\left\{\phi \left(\frac{u^{2}+|\nabla u|^{2}}{2}\right)\nabla u\right\}+\phi\left (\frac{u^{2}+|\nabla u|^{2}}{2}\right)u
 =
 \lambda a(x)|u|^{q-1}u+b(x)|u|^{p-1}u, &\mbox{in } \Omega,\\
 u=0, &\mbox{on } \partial\Omega,
  \end{array}
 \right.
 \end{equation}
where $\lambda>0$, $\Omega\subset \mathbb{R}^{N}$ is a smooth bounded domain with $N>2$, $1<q+1<2<p+1<2^{\ast}=\frac{2N}{N-2}$ and the map $\phi:[0,+\infty)\rightarrow \mathbb{R}$ is of class $ C^{2}$ and satisfies the following conditions:
\begin{itemize}
 \item[$(\phi_1)$] there exist two constants $0< \rho_0\leq \rho_1$ such that
 $ 0<\max\{\frac{q+1}{2}, \frac{2}{p+1}\}\rho_{1}<\rho_0\leq \phi(s)\leq \rho_1$ for all $s\in [0,+\infty)$;
 \item[$(\phi_2)$] there exists a constant $ \rho_2>0$ such that
$|\phi'(s)|s+|\phi''(s)|s^{2} \leq \rho_2 $ for all $s\in [0,+\infty)$;
 \item[$(\phi_3)$] there exist two constants $0< \rho_3\leq \rho_4$ such that
 $0<\rho_{3}\leq (1-q)\phi(s)+2\phi'(s)s\leq\rho_4$ for all $s\in [0,+\infty)$;
 \item[$(\phi_4)$] there exists a constant $ \rho_5>0$ such that
 $\rho_{5}+2\phi'(s)s\leq (p-1)\phi(s)$ for all $s\in [0,+\infty)$;
 \item[$(\phi_5)$] there exists a constant $ \rho_6>0$ such that
 $(1-q)(1-p)\phi\left( s\right)+2(4-p-q)\phi'\left( s\right)s+4\phi''\left( s\right)s^{2}\leq -\rho_6<0$ for all $s\in [0,+\infty)$;
 \item[$(\phi_6)$] the map $t\rightarrow \Phi(t^{2})$ is strictly convex on $[0,+\infty)$, where $\Phi(s):=\int_{0}^{s}\phi(\tau)d\tau$;
 \item[$(\phi_7)$] there exists $\phi(\infty)>0$ such that $\phi(s)\rightarrow\phi(\infty)$ as $s\rightarrow +\infty$.
 \end{itemize}
Furthermore, for $a(x)$ and $b(x)$ we assume that
\begin{itemize}
 \item[$(H)$] $a,b:\Omega\rightarrow\mathbb{R}$ are continuous and sign changing on $\Omega$, and $a,b\in L^{\infty}(\Omega)$.
\end{itemize}
\par
Clearly, the quasilinear term of the operator (\ref{0.1.1}) depends on not only $\nabla u$  but also $u$.
This particular quasilinear term generally appears in the study of nonlinear optics model
which describes the propagation of self-trapped beam in a cylindrical optical fiber made from a self-focusing dielectric material
(See \cite{Stuart1996, Stuart1997,Stuart2001,Stuart2003,Stuart2005,Stuart2010} for the derivation).
\par
At present, there are limit achievements in the study of the existence and multiplicity of solutions to problem involving the operator (\ref{0.1.1}) (for example, see \cite{Stuart2011,Jeanjean2022,Pomponio2021} and references therein).
In \cite{Stuart2011}, Stuart considered the following second-order quasilinear elliptic problem:
\begin{equation}\label{1.3.2}
  -div\left\{\phi \left(\frac{u^{2}+|\nabla u|^{2}}{2}\right)\nabla u\right\}+\phi\left (\frac{u^{2}+|\nabla u|^{2}}{2}\right)u
  =
 f(x,u), \;\mbox{in }\; \Omega
 \end{equation}
with $f(x,u)=\lambda u+h$ under the Dirichlet boundary condition $ u(x)=0$,
$ x\in \partial\Omega$,
where $\Omega$ is a bounded open subset of $\mathbb{R}^{N}$, $\lambda\in\mathbb{R}$, $h\in L^{2}(\Omega)$ and $h\geq0$ a.e. in $\Omega$.
By using an improved mountain pass theorem and the definetion of localizing the Palais-Smale sequence,
they obtained the existence of mountain pass type non-trivial solution for problem (\ref{1.3.2}).
Furthermore, the other solution is obtained as a local minimum of problem (\ref{1.3.2}) in a neighbourhood of the origin.
In \cite{Jeanjean2022}, Jeanjean and Radulescu further considered the cases when $f(x,u)$ has either a sublinear or a linear growth.
In the sublinear case, they assumed that $\phi$ satisfies some reasonable conditions,
by using a minimization procedure,
they obtained the existence of non-negative solutions.
Furthermore, they also proved a nonexistence result.
In the linear growth case, under stronger assumptions on $\phi$, by using the mountain pass theorem together with the Ekeland's variational theorem,
they obtained the existence of at least one or two non-negative solutions.
Inspired by \cite{ Stuart2011,Jeanjean2022}, in \cite{Pomponio2021}, Pomponio and Watanabe considered the cases when $\Omega=\mathbb{R}^{N}$, $N>2$ and $f(x,u)$ is of Berestycki-Lions' type.
Under appropriate assumptions on $\phi$, by the mountain pass theorem together with a technique of adding one dimension of space,
they proved the existence of a non-trivial weak solution.
Then by establishing the Pohozaev identity, they also obtained the existence of a radial ground state solution and a ground state solution under stronger assumptions on $\phi$.
\par
When $\phi(t)\equiv1$ (or more generally when $\phi(t)$ is a positive constant), $a(x)\equiv1$, $b(x)\equiv1$ and $u>0$ in $\Omega$, then problem (\ref{eq1}) reduces to the classical semilinear elliptic problem involving concave-convex term:
\begin{equation}\label{eq2}
 \left\{
  \begin{array}{ll}
 -\Delta u+u =\lambda|u|^{q-1}u+|u|^{p-1}u, &\mbox{in } \Omega,\\
  u>0, &\mbox{in } \Omega,\\
 u=0, &\mbox{on } \partial\Omega,
  \end{array}
 \right.
 \end{equation}
which has been treated by Ambrosetti, Brezis and Cerami\cite{Ambrosetti1994}, where $\Omega$ is a bounded regular domain of $\mathbb{R}^{N}$($N>2$), with smooth boundary and $1<q<2<p<2^{\ast}$.
In \cite{Ambrosetti1994}, combining the method of sub and super-solutions with the variational method, the authors proved the existence of $\lambda_{0}>0$ such that there are two solutions when $\lambda\in(0,\lambda_{0})$, one solution when $\lambda=\lambda_{0}$, and no solution when $\lambda>\lambda_{0}$.
The idea in \cite{Ambrosetti1994} has been extened to a lot of different problems involving concave-convex term in bounded domains, for examples,
the semilinear elliptic equation \cite{Bartsch1995},
$p$-Laplacian problem \cite{Azorero2000,Papageorgiou2019},
elliptic problem involving the fractional Laplacian operator \cite{C2013},
uniformly elliptic fully nonlinear equations \cite{Charro2009} and
Kirchhoff-type problems \cite{Chen2021}.
Moreover, the seminal work \cite{Ambrosetti1994} motivated many authors to study the relationship of the sign of the weight functions and the number of solutions for elliptic problems with concave-convex nonlinearities in bounded domains.
For example, \cite{Huang2013,Brown2007}.
In \cite{Huang2013}, Huang-Chen-Xiu studied the $p_{th}$ order Kirchhoff equation involving with concave-convex nonlinearities.
When the weight functions are both strictly positive, they proved that the problem has at least one positive solution via Mountain Pass lemma.
Moreover, by the Fountain theorem, they also obtained the problem has infinitely many solutions.
Especially, in \cite{Brown2007}, Brown and Wu studied the second order Laplace equation involving concave-convex nonlinearities and sign-changing weight functions.
By using the method of Nehari manifold and the fibering map which is introduced by \cite{Pohozaev1979,Pohozaev1988}, they obtained that the problem has at least two nontrivial solutions.
Subsequently, the idea in \cite{Brown2007} has been applied to a lot of different problems, for examples,
$p$-Laplacian equation \cite{Silva2018,Afrouzi2009},
Fractional $p$-Laplacian problems \cite{Silva2023,Fareh},
poly-Laplacian system on local finite graphs \cite{Yang},
quasilinear elliptic problems involving $\Phi$-Laplacian operator \cite{Carvalho2017,Silva2019},
double phase problem \cite{Mishra2023} and so on. Here we also refer to
\cite{Pohozaev1990,Tarantello1992,Drabek1997,Benkirane1997,Brown2003,Wu2006} where the authors established a precise description on the fibering map.
\par
Inspired by \cite{Stuart2011,Jeanjean2022,Pomponio2021,Brown2007,Carvalho2017}, in this paper, we shall investigate the multiplicity of weak solutions for equation (\ref{eq1}) involving concave-convex nonlinearities by using the Nehari manifold and doing a fine analysis associated on the fibering map.
We assume that the weight functions are both continuous and sign changing.
We split the Nehari manifold into two parts when the parameter $\lambda$ is small enough, and then by standard minimization procedure in each part, we obtain that the problem has at least two nontrivial weak solutions and one of them is the ground state solution.
Our results develop those in some known references in the following sense:
\begin{itemize}
\item[(\uppercase\expandafter{\romannumeral1})]
 Obviously, the reaction term in this paper is formally different from the reaction term in \cite{ Stuart2011} and does not satisfy the assumptions in \cite{Jeanjean2022,Pomponio2021}.
 For example, it is easy to verify that $f(x,u)=\lambda a(x)|u|^{q-1}u+b(x)|u|^{p-1}u$ does not satisfy the assumptions $(g_2)$ and $(g_4)$ in \cite{Pomponio2021}, where $\lambda\in \mathbb{R}$, $0<q<1<p<2^{\ast}$, $a,b$ are two continuous and sign-changing functions;
\item[(\uppercase\expandafter{\romannumeral2})]
When $\phi(t)\equiv1$(or more generally when $\phi(t)$ is a positive constant), the operator (\ref{0.1.1}) reduces to the usual Laplace operator, which has been studied in \cite{Brown2007}.
Therefore, our results contain the results in \cite{Brown2007} and cover more situations (see Section 5).
Moreover, our hypothesis $(\phi_1)$ means that the function $\phi(t)$ lies between two positive constants, which prohibits the consideration of certain types of operators, such as $p$-Laplacian, $(p,q)$-Laplacian and $\Phi$-Laplacian operators as defined in \cite{Carvalho2017,Silva2019}.
Thus, our problem is different from those problems in \cite{Silva2023,Fareh,Carvalho2017,Silva2019};
\item[(\uppercase\expandafter{\romannumeral3})]
 Because of the inhomogeneity of the operator (\ref{0.1.1}), our proofs become more difficult and complex than those in \cite{Brown2007}.
 Moreover, different from the $p$-Laplacian as defined in \cite{Silva2023,Fareh} and $\Phi$-Laplacian operators as defined in \cite{Carvalho2017,Silva2019}, the operator (\ref{0.1.1}) depends not only on $\nabla u$  but also on $u$.
 Therefore, our proofs are different from those in \cite{Silva2023,Fareh,Carvalho2017,Silva2019}, especially in terms of the scaling of various inequalities related to the operator (\ref{0.1.1}).
\item[(\uppercase\expandafter{\romannumeral4})]
 Similar to the approach in \cite{Brown2007}, we split the Nehari manifold into two parts when the parameter $\lambda$ is small enough, and then obtain a bounded minimizing sequence $\{u_{n}\}$ in each part.
 It is not hard to see that the boundedness of $\{u_{n}\}$ guarantees $u_{n}\rightharpoonup u$ in $H^{1}_{0}(\Omega)$.
 In order to prove that $u_{n}\rightarrow u$ in $H^{1}_{0}(\Omega)$,
 we need to use the Brezis-Lieb Lemma in \cite{Brezis1983}.
 However, since the operator (\ref{0.1.1}) is quasi-linear and it depends on not only $u$ but also $\nabla u$,  the weak convergence of $\{u_{n}\}$ in $H^{1}_{0}(\Omega)$ is not enough to guarantee that the Brezis-Lieb Lemma holds, and to this end, it is necessary to obtain the almost everywhere convergence of  $\nabla u_{n}$, which is actually not an easy work.
 We apply Lemma 3.5 in \cite{Jeanjean2022} to overcome this difficulty, and we observe that
$J'_{\lambda}(u_{n})\rightarrow 0$ in $(H^{1}_{0}(\Omega))'$ is the critical property to prove that $\{u_{n}\}$ satisfies those assumptions in Lemma 3.5 in \cite{Jeanjean2022}. Hence, solving this difficulty is transformed into proving that $J'_{\lambda}(u_{n})\rightarrow 0$ in $(H^{1}_{0}(\Omega))'$. We combine the Ekeland's variational principle \cite{Ekeland1974}  with the ideas in \cite{Tarantello1992} to prove that the minimizing sequence $\{u_{n}\}$ is actually a Palais-Smale sequence which implies the fact that $J'_{\lambda}(u_{n})\rightarrow 0$ in $(H^{1}_{0}(\Omega))'$.
\end{itemize}
\par
Next, we state our main results.
\vskip2mm
 \noindent
{\bf Theorem 1.1. } {\it Assume that $(\phi_1)$-$(\phi_7)$ and $(H)$ hold.
Then for each $\lambda$ satisfying $0<\lambda<\lambda_{0}=\min\{\lambda_{1},\lambda_{2}\}$, where
\begin{eqnarray}
\label{aa1}& &
\lambda_{1}=\frac{ \rho_{5}}{(p-q)\|a\|_{\infty}S_{q+1}^{q+1}}\left(  \frac{\rho_{3}}{(p-q)\|b\|_{\infty}S_{p+1}^{p+1}}\right)^{\frac{1-q}{p-1}},
\;\;
\lambda_{2}=\delta^{1-\frac{q+1}{2}}/c_{1},\\
\label{aa2}& &
c_{1}=\frac{\|a\|_{\infty}S_{q+1}^{q+1}}{(q+1)(\frac{\rho_{0}}{2}-\frac{\rho_{1}}{p+1})^{\frac{q+1}{2}}},
\;\;
\delta=\left(\frac{\rho_{0}}{2}-\frac{\rho_{1}}{p+1}\right)\left(  \frac{\rho_{0}}{\|b\|_{\infty}S_{p+1}^{p+1}} \right)^{\frac{2}{p-1}},
\end{eqnarray}
problem (\ref{eq1}) admits at least one positive energy solution $u_{0}$ and negative energy solution $\tilde{u}_{0}$ in $H^{1}_{0}(\Omega)$.
Furthermore,  $\tilde{u}_{0}$ is the ground state solution of (\ref{eq1}).
}

\section{Notations and preliminaries}
Let $\Omega\subset \mathbb{R}^{N}$ be  a smooth bounded domain with $N>2$.
\par
Let $1\leq p<\infty$. $ L^{p}(\Omega)=\{u :\Omega\rightarrow\mathbb{R}:\int_{\Omega}|u|^{p}dx<\infty\}$ with norm
\begin{eqnarray*}
\|u\|_{p}=\left(\int_{\Omega}|u|^{p}dx\right)^{\frac{1}{p}}.
\end{eqnarray*}
\par
$L^{\infty}(\Omega)=\{u :\Omega\rightarrow\mathbb{R}: u \;\mbox{ is measurable on}\;\Omega \;\mbox{and}\; \inf_{E\subset\Omega, |E|=0}\sup_{\Omega\backslash E}|f(x)|<\infty\}$ with norm
\begin{eqnarray*}
\|u\|_{\infty}=\inf_{|\Omega_{0}|=0,\Omega_{0}\subset \Omega}\left( \sup_{\Omega\setminus \Omega_{0}}|u| \right).
\end{eqnarray*}
\par
$W^{1,2}(\Omega)=\{u\in L^{2}(\Omega):\frac{\partial u}{\partial x_{i}}\in L^{2}(\Omega),i=1,2,...,N\}$.
\par
$C^{\infty}(\Omega)=\{ u :\Omega\rightarrow\mathbb{R} :u \;\mbox{is infinitely differentiable on}\;\Omega\}$.
\par
$C_{0}^{\infty}(\Omega)=\{ u\in C^{\infty}(\Omega):u\;\mbox{has compact support on }\;\Omega\}$.
\par
We denote by $W_{0}^{1,2}(\Omega)$ the closure of $C_{0}^{\infty}(\Omega)$ in $W^{1,2}(\Omega)$, usually recorded as $H_{0}^{1}(\Omega)$, with norm
\begin{eqnarray}\label{2.1}
\|u\|=\|\nabla u\|_{2}=\left(\int_{\Omega}|\nabla u|^{2}dx\right)^{\frac{1}{2}}.
\end{eqnarray}
\par
We denote the dual space of $H^{1}_{0}(\Omega)$ by $(H^{1}_{0}(\Omega))^{'}$.
\par
We denote the weak convergence by $\rightharpoonup$ and denote the strong convergence by $\rightarrow$.
\par
 $S_{i}$ denotes the best Sobolev constants for the imbeddings $H_{0}^{1}(\Omega)\hookrightarrow L^{i}(\Omega)$, where $1<i<2^{\ast}$.
\vskip2mm
 \noindent
{\bf Proposition 2.1. }(Sobolev imbedding theorem \cite{Adams2003})
Let $\Omega$ be an open subset of $\mathbb{R}^{N}$.
If $|\Omega|<\infty$, the following imbeddings are continuous:
\begin{eqnarray*}
H_{0}^{1}(\Omega)\hookrightarrow L^{\gamma}(\Omega),\;\;1\leq\gamma\leq2^{\ast}.
\end{eqnarray*}
Furthermore, the following imbeddings are compact:
\begin{eqnarray*}
H_{0}^{1}(\Omega)\hookrightarrow L^{\gamma}(\Omega),\;\;1\leq\gamma<2^{\ast},
\end{eqnarray*}
where
\begin{eqnarray*}
2^{\ast}= \begin{cases}
               \frac{2N}{N-2}, & \mbox{  if } N>2,\\
               +\infty, & \mbox{ if } N\leq2.
\end{cases}
\end{eqnarray*}
\par
We define the variational functional $J_{\lambda}:H_{0}^{1}(\Omega)\rightarrow \mathbb{R}$ corresponding to problem (\ref{eq1}) by
\begin{eqnarray*}
J_{\lambda}(u)
=
\int_{\Omega}\Phi\left(\frac{u^{2}+|\nabla u|^{2}}{2}\right)dx
-\frac{\lambda}{q+1}\int_{\Omega}a(x)|u|^{q+1}dx-\frac{1}{p+1}\int_{\Omega}b(x)|u|^{p+1}dx,
\;\;
u\in H_{0}^{1}(\Omega).
\end{eqnarray*}
As can be seen from the results in Appendix A.1,
assumptions $(\phi_1)$, $(H)$ and $1<q+1<2<p+1<2^{\ast}$ ensure that the functional $J_{\lambda}$ is well defined in $H_{0}^{1}(\Omega)$ for any $\lambda>0$.
Further, if $\phi\in C([0,+\infty),\mathbb{R})$, $(\phi_1)$, $(H)$ and $1<q+1<2<p+1<2^{\ast}$ hold,
the functional $J_{\lambda}\in C^{1}(H_{0}^{1}(\Omega),\mathbb R)$ for any $\lambda>0$,
and the derivation of functional $J_{\lambda}$ satisfies
\begin{eqnarray}\label{2.1.1}
\langle J'_{\lambda}(u),v\rangle
=
\int_{\Omega}\phi\left(\frac{u^{2}+|\nabla u|^{2}}{2}\right)(uv+\nabla u\cdot \nabla v )dx
-\lambda\int_{\Omega}a(x)|u|^{q-1}uvdx-\int_{\Omega}b(x)|u|^{p-1}uvdx, \ \ \forall u,v\in H_{0}^{1}(\Omega).
\end{eqnarray}
It is easy to see that (\ref{eq1}) is the Euler-Lagrange equation of the functional $J_{\lambda}(u)$. So finding a weak solution for the problem (\ref{eq1}) is equivalent to finding a critical point for the functional $J_{\lambda}$.
\par
The super-quadratic nonlinear term leads to the functional $J_{\lambda}$ is not bounded below on $H_{0}^{1}(\Omega)$ and so the minimization method on $H_{0}^{1}(\Omega)$ failed.
Moreover, since $a,b$ are sign-changing on $\Omega$, the well known Ambrosetti-Rabinowitz condition (\cite{Ambrosetti1973}) and even the non-quadratic condition introduced by Costa-Magalh\~aes \cite{Costa1994} do not work any more.
It is well known that if $u$ is a non-trivial weak solution of the problem (\ref{eq1}),
then it must belong to the Nehari manifold $N_{\lambda}$ given by
\begin{eqnarray}\label{2.0.1}
N_{\lambda}:
&= &
 \{u\in H_{0}^{1}(\Omega)\backslash\{0\}:\langle J'_{\lambda}(u),u\rangle=0\} \nonumber\\
&= &
 \left \{  u\in H_{0}^{1}(\Omega)\backslash\{0\}:
 \int_{\Omega}\phi\left(\frac{u^{2}+|\nabla u|^{2}}{2}\right)(u^{2}+|\nabla u|^{2} )dx
  =\lambda\int_{\Omega}a(x)|u|^{q+1}dx+\int_{\Omega}b(x)|u|^{p+1}dx\right\}.
 \end{eqnarray}
It follows from $\langle J'_{\lambda}(u),u\rangle=0$ that the functional $J_{\lambda}$ has two equivalent forms in $N_{\lambda}\subset H_{0}^{1}(\Omega)$:
\begin{eqnarray}\label{2.0.2}
      J_{\lambda}(u)=\int_{\Omega}\Phi\left(\frac{u^{2}+|\nabla u|^{2}}{2}\right)dx
      -\frac{1}{q+1}\int_{\Omega}\phi\left(\frac{u^{2}+|\nabla u|^{2}}{2}\right)(u^{2}+|\nabla u|^{2} )dx
      +\left(\frac{1}{q+1}-\frac{1}{p+1}\right)\int_{\Omega}b(x)|u|^{p+1}dx
\end{eqnarray}
and
\begin{eqnarray}\label{2.0.3}
      J_{\lambda}(u)=\int_{\Omega}\Phi\left(\frac{u^{2}+|\nabla u|^{2}}{2}\right)dx
      -\frac{1}{p+1}\int_{\Omega}\phi\left(\frac{u^{2}+|\nabla u|^{2}}{2}\right)(u^{2}+|\nabla u|^{2} )dx
      -\lambda\left(\frac{1}{q+1}-\frac{1}{p+1}\right)\int_{\Omega}a(x)|u|^{q+1}dx.
\end{eqnarray}
Note that (\ref{2.0.3}) ensures that the functional $J_{\lambda}$ is bounded below on $N_{\lambda}$ (see Proposition 2.2 below).
It permits us to look for the least energy critical point for the functional $J_{\lambda}$.
\vskip2mm
 \noindent
{\bf Proposition 2.2. } {\it Assume that $(\phi_{1})$ and $(H)$ hold. Then the functional $J_{\lambda}$ is coercive and bounded below on $N_{\lambda}$.\\}
{\bf Proof.} It follows from (\ref{2.0.3}), $(\phi_1)$, $(H)$, $1<q+1<2<p+1<2^{\ast}$ and Proposition 2.1 that
\begin{eqnarray*}
J_{\lambda}(u)
& \geq &
\left(\frac{\rho_{0}}{2}-\frac{\rho_{1}}{p+1}\right)(\|u\|_{2}^{2}+\|\nabla u\|_{2}^{2})
-\lambda\left(\frac{1}{q+1}-\frac{1}{p+1}\right)\|a\|_{\infty}\|u\|_{q+1}^{q+1}\\
& \geq &
\left(\frac{\rho_{0}}{2}-\frac{\rho_{1}}{p+1}\right)\|\nabla u\|_{2}^{2}
-\lambda\left(\frac{1}{q+1}-\frac{1}{p+1}\right)\|a\|_{\infty}S^{q+1}_{q+1}\|\nabla u\|_{2}^{q+1}.
\end{eqnarray*}
This also shows that $J_{\lambda}(u)\rightarrow+\infty$ as $ \|u\|\rightarrow\infty$ when $q+1<2$.
Thus, the functional $J_{\lambda}$ is coercive and bounded below on $N_{\lambda}$.
\qed
\par
In order to get the multiplicity of solutions,
we will split $N_{\lambda}$ with the help of the fibering map.
For any $u\in H_{0}^{1}(\Omega)$, we define the fibering map $\gamma_{u}:(0,+\infty)\rightarrow\mathbb{R}$ given by
\begin{eqnarray}\label{2.0.4}
 \gamma_{u}(t)
 :=
 J_{\lambda}(tu)
 =
 \int_{\Omega}\Phi\left(\frac{u^{2}+|\nabla u|^{2}}{2}t^{2}\right)dx
 -\frac{\lambda}{q+1}t^{q+1}\int_{\Omega}a(x)|u|^{q+1}dx
 -\frac{1}{p+1}t^{p+1}\int_{\Omega}b(x)|u|^{p+1}dx.
\end{eqnarray}
As can be seen from Appendix A.2,
assumptions $\phi\in C([0,+\infty), \mathbb{R})$, $(\phi_1)$, $(H)$ and $1<q+1<2<p+1<2^{\ast}$ ensure that $\gamma_{u}\in C^{1}((0,+\infty), \mathbb{R})$ with the derivative given by
\begin{eqnarray}\label{2.0.25}
\gamma_{u}'(t)
=
\int_{\Omega} t\phi\left( \frac{u^{2}+|\nabla u|^{2}}{2}t^{2}\right)(u^{2}+|\nabla u|^{2})dx
-\lambda t^{q}\int_\Omega a(x)|u|^{q+1}dx-t^{p}\int_\Omega b(x)|u|^{p+1}dx.
 \end{eqnarray}
 \noindent
{\bf Remark 2.1. }  It is easy to see that $tu\in N_{\lambda}$ if and only if $\gamma_{u}'(t)=0$,
which shows that the elements in $N_{\lambda}$ correspond to stationary points of the fibering map $\gamma_{u}$.
More specifically, $u\in N_{\lambda}$ if and only if $\gamma_{u}'(1)=0$.
Based on this connection, we expect to split $N_{\lambda}$ by the higher derivative value of $\gamma_{u}$ at $t=1$.
\vskip2mm
\par
As shown in Appendix A.2, if $\phi\in C^{1}([0,+\infty), \mathbb{R})$, $(\phi_{1})$, $(\phi_{2})$ and $(H)$ hold.
Then $\gamma_{u}\in C^{2}((0,+\infty), \mathbb{R})$ with the second order derivative given by
\begin{eqnarray}\label{2.0.26}
\gamma_{u}''(t)
& := &
\int_{\Omega}\left[ \phi\left( \frac{u^{2}+|\nabla u|^{2}}{2}t^{2}\right)(u^{2}+|\nabla u|^{2})
+t^{2}\phi'\left( \frac{u^{2}+|\nabla u|^{2}}{2}t^{2}\right)(u^{2}+|\nabla u|^{2})^{2}\right]dx\nonumber\\
& &
- \lambda q t^{q-1}\int_{\Omega} a(x)|u|^{q+1}dx-pt^{p-1}\int_{\Omega}b(x)|u|^{p+1}dx.
\end{eqnarray}
Clearly, for any $u\in N_{\lambda}$,
\begin{eqnarray}\label{2.6}
\gamma_{u}''(1)
 &= &
 \int_{\Omega}\left[\phi'\left( \frac{u^{2}+|\nabla u|^{2}}{2}\right)(u^{2}+|\nabla u|^{2})^{2}
 +(1-q)\phi\left( \frac{u^{2}+|\nabla u|^{2}}{2}\right)(u^{2}+|\nabla u|^{2})\right]dx\nonumber \\
 & &
 -(p-q)\int_{\Omega}b(x)|u|^{p+1}dx
 \end{eqnarray}
is equivalent to
\begin{eqnarray}\label{2.7}
 \gamma_{u}''(1)
 &= &
 \int_{\Omega}\left[\phi'\left( \frac{u^{2}+|\nabla u|^{2}}{2}\right)(u^{2}+|\nabla u|^{2})^{2}
 -(p-1)\phi\left( \frac{u^{2}+|\nabla u|^{2}}{2}\right)(u^{2}+|\nabla u|^{2})\right]dx\nonumber \\
 & &
+\lambda(p-q)\int_{\Omega}a(x)|u|^{q+1}dx.
\end{eqnarray}
As was pointed by Brown-Wu in \cite{Brown2007}, it is natural to divide $N_{\lambda}$ into three submanifolds:
\begin{eqnarray*}
  N_{\lambda}^{+}:=\{u\in N_{\lambda}: \gamma_{u}''(1)>0 \};\\
  N_{\lambda}^{-}:=\{u\in N_{\lambda}: \gamma_{u}''(1)<0 \};\\
  N_{\lambda}^{0}:=\{u\in N_{\lambda}: \gamma_{u}''(1)=0 \}.
\end{eqnarray*}

 \par
 We shall use the minimization method to find the critical points on these submanifolds, respectively.
If we can obtain a minimizer $u$ on the submanifold, it remains to prove that the submanifold is a natural constraint on $u$, and
the Lagrange multiplier method is a very useful tool to prove the fact for $u\not\in N_{\lambda}^{0}$.
However, this method is invalid for $u\in N_{\lambda}^{0}$, and
it is unknown whether the functional $J_{\lambda}$ has critical points on $N_{\lambda}^{0}$.
As in \cite{Brown2007}, we will control the parameter $\lambda$ to be small enough such that $N_{\lambda}^{0}=\emptyset$ so that we shall seek the minimizers $u$ only on the submanifolds $N_{\lambda}^{+}$ and $N_{\lambda}^{-}$.
Then, we have the following results.
\vskip2mm
 \noindent
{\bf Lemma 2.1. } {\it Suppose that $(\phi_1)$, $(\phi_3)$, $(\phi_4)$ and $(H)$ hold.
Then $N_{\lambda}^{0}=\emptyset$ for any $\lambda\in (0,\lambda_{1})$ where $\lambda_{1}$ is given by (\ref{aa1}).
}
\vskip0mm
\noindent
{\bf Proof.} Let $\lambda\in (0,\lambda_{1})$ be fixed.
Arguing by contradiction, we assume that $N_{\lambda}^{0}\neq\emptyset$.
Taking $u\in N_{\lambda}^{0}$ as a fixed function, it follows that $\gamma_{u}'(1)=\gamma_{u}''(1)=0$.
Firstly, by (\ref{2.6}), we obtain that
\begin{eqnarray*}
 \int_{\Omega}\left[\phi'\left( \frac{u^{2}+|\nabla u|^{2}}{2}\right)(u^{2}+|\nabla u|^{2})^{2}
 +(1-q)\phi\left( \frac{u^{2}+|\nabla u|^{2}}{2}\right)(u^{2}+|\nabla u|^{2})\right]dx
 =
 (p-q)\int_{\Omega}b(x)|u|^{p+1}dx.
\end{eqnarray*}
In view of $(\phi_3)$, we obtain that
\begin{eqnarray}\label{2.0.7}
\int_{\Omega}\left[\phi'\left( \frac{u^{2}+|\nabla u|^{2}}{2}\right)(u^{2}+|\nabla u|^{2})^{2}
 +(1-q)\phi\left( \frac{u^{2}+|\nabla u|^{2}}{2}\right)(u^{2}+|\nabla u|^{2})\right]dx
\geq
\rho_{3}\|\nabla u\|_{2}^{2}.
\end{eqnarray}
Then it follows from $q+1<p+1$, $(H)$ and Proposition 2.1 that
\begin{eqnarray*}
(p-q)\int_{\Omega}b(x)|u|^{p+1}dx
\leq
(p-q)\|b\|_{\infty}S^{p+1}_{p+1}\|\nabla u\|_{2}^{p+1},
\end{eqnarray*}
which together with (\ref{2.0.7}) implies that
\begin{eqnarray}\label{2.0.8}
\|\nabla u\|_{2}
\geq
\left[  \frac{\rho_{3}}{(p-q)\|b\|_{\infty}S_{p+1}^{p+1}}\right]^{\frac{1}{p-1}}.
 \end{eqnarray}
Furthermore, by (\ref{2.7}), we get
\begin{eqnarray*}
 \int_{\Omega}\left[ (p-1)\phi\left( \frac{u^{2}+|\nabla u|^{2}}{2}\right)(u^{2}+|\nabla u|^{2})
 -\phi'\left( \frac{u^{2}+|\nabla u|^{2}}{2}\right)(u^{2}+|\nabla u|^{2})^{2}\right]dx
=
\lambda(p-q)\int_{\Omega}a(x)|u|^{q+1}dx.
\end{eqnarray*}
In view of $(\phi_4)$, we have
\begin{eqnarray}\label{2.0.9}
\int_{\Omega}\left[ (p-1)\phi\left( \frac{u^{2}+|\nabla u|^{2}}{2}\right)(u^{2}+|\nabla u|^{2})
 -\phi'\left( \frac{u^{2}+|\nabla u|^{2}}{2}\right)(u^{2}+|\nabla u|^{2})^{2}\right]dx
\geq
\rho_{5}\|\nabla u\|_{2}^{2}.
\end{eqnarray}
Then taking into account that $q+1<p+1$, $(H)$ and Proposition 2.1, we observe that
\begin{eqnarray*}
\lambda(p-q)\int_{\Omega}a(x)|u|^{q+1}dx
 \leq
\lambda(p-q)\|a\|_{\infty}S^{q+1}_{q+1}\|\nabla u\|_{2}^{q+1},
\end{eqnarray*}
which together with (\ref{2.0.9}) implies that
\begin{eqnarray}\label{2.0.10}
\lambda
\geq
\frac{ \rho_{5}\|\nabla u\|_{2}^{1-q}}{(p-q)\|a\|_{\infty}S_{q+1}^{q+1}}.
 \end{eqnarray}
Now, thanks to (\ref{2.0.8}) and (\ref{2.0.10}), we conclude that
\begin{eqnarray*}
\lambda
\geq
\frac{ \rho_{5}}{(p-q)\|a\|_{\infty}S_{q+1}^{q+1}}\left(  \frac{\rho_{3}}{(p-q)\|b\|_{\infty}S_{p+1}^{p+1}}\right)^{\frac{1-q}{p-1}}
=\lambda_{1}.
 \end{eqnarray*}
This is a contradiction with $\lambda\in (0,\lambda_{1})$.
Hence $N_{\lambda}^{0}=\emptyset$.
\qed
\vskip2mm
\par
Next, we shall prove that $N_{\lambda}$ is a $C^{1}$-manifold.
\vskip2mm
 \noindent
{\bf Lemma 2.2.} {\it Suppose that $(\phi_1)$, $(\phi_2)$, $(\phi_3)$, $(\phi_4)$ and $(H)$ hold.
Then $N_{\lambda}$ is a $C^{1}$-manifold for any $\lambda\in (0,\lambda_{1})$.}
\vskip0mm
 \noindent
{\bf Proof.} Fix $\lambda$ such that $\lambda\in (0,\lambda_{1})$.
By Lemma 2.1, we have $N_{\lambda}=N_{\lambda}^{+}\cup N_{\lambda}^{-}$.
We define $G:N_{\lambda}\rightarrow\mathbb{R}$ given by
\begin{eqnarray*}
 G(u)
& = & \langle J'_{\lambda}(u),u\rangle \\
& = & \int_{\Omega}\phi\left(\frac{u^{2}+|\nabla u|^{2}}{2}\right)(u^{2}+|\nabla u|^{2} )dx
  -\lambda\int_{\Omega}a(x)|u|^{q+1}dx-\int_{\Omega}b(x)|u|^{p+1}dx,
 \;\;u\in N_{\lambda}.
 \end{eqnarray*}
Firstly, we prove that $G\in C^{1}(N_\lambda,\mathbb R)$ and
\begin{eqnarray}\label{2.0.11}
\langle G'(u),v\rangle
& = &
\int_{\Omega}\left[\phi'\left(\frac{u^{2}+|\nabla u|^{2}}{2}\right)(u^{2}+|\nabla u|^{2})
+2\phi\left(\frac{u^{2}+|\nabla u|^{2}}{2}\right)\right]
\cdot(uv+\nabla u\nabla v)dx\nonumber\\
& &
-(q+1)\lambda\int_{\Omega} a(x)|u|^{q-1}uvdx
-(p+1)\int_{\Omega}b(x)|u|^{p-1}uvdx,\;\;u,v\in  N_{\lambda}.
 \end{eqnarray}
Defining $F:\Omega\times\mathbb{R}\times\mathbb{R}^{N}\rightarrow\mathbb{R}$ by
\begin{eqnarray*}
F(x,u,\xi)
=
\phi\left(\frac{u^{2}+|\xi|^{2}}{2}\right)(u^{2}+|\xi|^{2} )
-\lambda a(x)|u|^{q+1}
-b(x)|u|^{p+1}.
\end{eqnarray*}
It follows that
\begin{eqnarray*}
& &
F_{u}(x,u,\xi)
=
\phi'\left(\frac{u^{2}+|\xi|^{2}}{2}\right)(u^{2}+|\xi|^{2})u
+\phi\left(\frac{u^{2}+|\xi|^{2}}{2}\right)2u
-(q+1)\lambda a(x)|u|^{q-1}u
-(p+1)b(x)|u|^{p-1}u,\\
& &
F_{\xi}(x,u,\xi)
=
\phi'\left(\frac{u^{2}+|\xi|^{2}}{2}\right)(u^{2}+|\xi|^{2})\xi
+\phi\left(\frac{u^{2}+|\xi|^{2}}{2}\right)2\xi.
\end{eqnarray*}
By $(\phi_1)$ and $(H)$, we have
\begin{eqnarray*}
|F(x,u,\xi)|
\leq
\rho_{1}(u^{2}+|\xi|^{2})
+\lambda\|a\|_{\infty}|u|^{q+1}
+\|b\|_{\infty}|u|^{p+1}.
\end{eqnarray*}
Combing with $1<q+1<2<p+1<2^{\ast}$, we also mention that
\begin{eqnarray*}
& &
|F(x,u,\xi)|
\leq
C(1+|\xi|^{2}),\;\mbox{if}\; |u|\leq1,\\
& &
|F(x,u,\xi)|
\leq
C(|u|^{p+1}+|\xi|^{2}),\;\mbox{if}\; |u|>1,
\end{eqnarray*}
which imply that
\begin{eqnarray*}
|F(x,u,\xi)|
\leq
C(1+|u|^{p+1}+|\xi|^{2}).
\end{eqnarray*}
Using $(\phi_1)$, $(\phi_2)$ and $(H)$, we have
\begin{eqnarray*}
|F_{u}(x,u,\xi)|
\leq
2(\rho_{2}+\rho_{1})|u|+(q+1)\lambda\|a\|_{\infty}|u|^{q}+(p+1)\|b\|_{\infty}|u|^{p},
\end{eqnarray*}
which together with $1<q+1<2<p+1<2^{\ast}$ implies that
\begin{eqnarray*}
|F_{u}(x,u,\xi)|
\leq
C(1+|u|^{p}).
\end{eqnarray*}
Moreover, in view of $(\phi_1)$ and $(\phi_2)$, we obtain that
\begin{eqnarray*}
|F_{\xi}(x,u,\xi)|
\leq
2(\rho_{2}+\rho_{1})|\xi|.
\end{eqnarray*}
Hence, using \cite[Theorem C.1]{Struwe1996}, we have $G\in C^{1}(N_\lambda,\mathbb R)$ and (\ref{2.0.11}) holds.
\par
Next, we present the proof that $\langle G'(u),u\rangle\neq 0$ for any $u\in N_{\lambda}$.
Clearly, $u\in N_{\lambda}$ implies that $\gamma'_{u}(1)=0$.
In view of (\ref{2.0.25}), (\ref{2.0.26}) and (\ref{2.0.11}), we deduce that
\begin{eqnarray*}
\langle G'(u),u\rangle
& = &
\int_{\Omega}\left[\phi'\left(\frac{u^{2}+|\nabla u|^{2}}{2}\right)(u^{2}+|\nabla u|^{2})^{2}
+\phi\left(\frac{u^{2}+|\nabla u|^{2}}{2}\right)(u^{2}+|\nabla u|^{2})\right]dx\\
& &
-q\lambda \int_{\Omega}a(x)|u|^{q+1}dx-p\int_{\Omega}b(x)|u|^{p+1}dx\\
& = &
\gamma''_{u}(1).
 \end{eqnarray*}
Clearly, $\langle G'(u),u\rangle\neq 0$ for any $u\in N_{\lambda}$.
Moreover, using (\ref{2.0.1}) we also have that $G(u)=\langle J'_{\lambda}(u),u\rangle=0$.
Hence, $0\in \mathbb{R}$ is a regular value for $G$ and $N_{\lambda}=G^{-1}(0)$ is a $C^{1}$-manifold.
\qed
\vskip2mm
\par
Now we are able to prove that any critical points for $J_{\lambda}$ on $N_{\lambda}$ is a free critical points, i.e.,  a critical point in the whole space $H^{1}_{0}(\Omega)$.
\vskip2mm
 \noindent
{\bf Lemma 2.3.} {\it Assume that $(\phi_1)$, $(\phi_2)$, $(\phi_3)$, $(\phi_4)$ and $(H)$ hold.
For any $\lambda\in (0,\lambda_{1})$, if $ u_{0}\in N_{\lambda} $ is a local minimum or maximum of $J_{\lambda}$.
Then $u_{0}$ is a free critical point for $J_{\lambda}$. }
\vskip0mm
\noindent
{\bf Proof.} Fixed $\lambda\in (0,\lambda_{1})$.
By Lemma 2.1 and Lemma 2.2, we know that $N_{\lambda}=N_{\lambda}^{+}\cup N_{\lambda}^{-}$ and $N_{\lambda}$ is a $C^{1}$-manifold.
Without any loss of generality, we assume that $u_{0}\in N_{\lambda} $ is a local minimum of $J_{\lambda}$.
It is easy to verify that $u_{0}$ is a solution for the minimization problem
 \begin{equation}\label{2.0.27}
 \min\{J_{\lambda}(u),G(u)=0\},
 \end{equation}
where the definition of $G(u)$ is given in Lemma 2.2.
Since $N_{\lambda}$ is a $C^{1}$-manifold, we get
\begin{eqnarray}\label{2.0.28}
J_{\lambda}'(u_{0})=\mu G'(u_{0}),
 \end{eqnarray}
where $\mu\in \mathbb{R}$ is given by Lagrange multiplies theorem.
We are able to prove that $\mu=0$.
Indeed, firstly, by (\ref{2.0.28}) and $u_{0}\in N_{\lambda}$, we obtain that
\begin{eqnarray}\label{2.0.29}
0=\langle J_{\lambda}'(u_{0}),u_{0}\rangle=\mu\langle G'(u_{0}),u_{0}\rangle.
 \end{eqnarray}
Then, taking $u=v=u_{0}$ in (\ref{2.0.11}) we see that
\begin{eqnarray*}
\langle G'(u_{0}),u_{0}\rangle
&=&
\int_{\Omega}\left[\phi'\left(\frac{u_{0}^{2}+|\nabla u_{0}|^{2}}{2}\right)(u_{0}^{2}+|\nabla u_{0}|^{2})
+2\phi\left(\frac{u_{0}^{2}+|\nabla u_{0}|^{2}}{2}\right)\right]\cdot(u_{0}^{2}+|\nabla u_{0}|^{2})dx\\
& &
-(q+1)\lambda\int_{\Omega} a(x)|u_{0}|^{q+1}dx-(p+1)\int_{\Omega}b(x)|u_{0}|^{p+1}dx,
  \end{eqnarray*}
which together with (\ref{2.0.25}) and (\ref{2.0.26}) implies that
\begin{eqnarray*}
\langle G'(u_{0}),u_{0}\rangle
& = &
\int_{\Omega}\left[\phi'\left(\frac{u_{0}^{2}+|\nabla u_{0}|^{2}}{2}\right)(u_{0}^{2}+|\nabla u_{0}|^{2})^{2}
+\phi\left(\frac{u_{0}^{2}+|\nabla u_{0}|^{2}}{2}\right)(u_{0}^{2}+|\nabla u_{0}|^{2})\right]dx\\
& &
-q\lambda \int_{\Omega}a(x)|u_{0}|^{q+1}dx-p\int_{\Omega}b(x)|u_{0}|^{p+1}dx\\
& = &
\gamma''_{u_{0}}(1).
 \end{eqnarray*}
Since $N_{\lambda}=N_{\lambda}^{+}\cup N_{\lambda}^{-}$, $\langle G'(u_{0}),u_{0}\rangle\neq0$.
Combing with (\ref{2.0.29}), we mention that $\mu=0$.
Taking $\mu=0$ in (\ref{2.0.28}) we get that $J_{\lambda}'(u_{0})=0$.
Thus $u_{0}$ is a critical point for $J_{\lambda}$ on $H^{1}_{0}(\Omega)$.
\qed

\section{Analysis of the fibering map}
\par
In this section, we shall present a complete description of the fibering map $\gamma_{u}$.
\par
Firstly, by the definition of $\gamma_{u}'(t)$ we have that
\begin{eqnarray}\label{4.0.1}
\gamma_{u}'(t)
= t^{q}\left(m_{u}(t)-\lambda \int_{\Omega}a(x)|u|^{q+1}dx\right)
 \end{eqnarray}
where
\begin{eqnarray*}
m_{u}(t)
&= &\int_{\Omega}t^{1-q}\phi\left( \frac{u^{2}+|\nabla u|^{2}}{2}t^{2}\right)(u^{2}+|\nabla u|^{2})dx
-t^{p-q}\int_{\Omega}b(x)|u|^{p+1}dx.
 \end{eqnarray*}
Clearly, for any $t>0$, $tu\in N_{\lambda}$ if and only if $t$ is a solution of equation $m_{u}(t)=\lambda \int_{\Omega}a(x)|u|^{q+1}dx$.
Thus, in order to obtain the behavior of the fibering map $\gamma_{u}$, it is necessary to analyze the properties of the auxiliary function $m_{u}$.
To this end, we introduce the auxiliary function $\eta_{u}$ and obtain the following results.
\vskip2mm
 \noindent
{\bf Lemma 3.1. } {\it
For any  $u\in H_{0}^{1}(\Omega)\backslash\{0\}$, define $\eta_{u}:(0,+\infty)\rightarrow\mathbb{R}$  by
\begin{eqnarray*}
\eta_{u}(t)
=
\int_{\Omega}t^{1-p}(u^{2}+|\nabla u|^{2})\left[(1-q)\phi\left( \frac{u^{2}+|\nabla u|^{2}}{2}t^{2}\right)
+2\phi'\left( \frac{u^{2}+|\nabla u|^{2}}{2}t^{2}\right) \frac{u^{2}+|\nabla u|^{2}}{2}t^{2}\right]dx.
\end{eqnarray*}
Suppose that $(\phi_1)$, $(\phi_2)$, $(\phi_3)$ and $(\phi_5)$ hold.
Then $\lim_{t\rightarrow0^{+}}\eta_{u}(t)=+\infty$ and $\lim_{t\rightarrow+\infty}\eta_{u}(t)=0$.
Moreover,
\begin{eqnarray}\label{551}
\eta'_{u}(t)
& =&\int_{\Omega}(1-q)(1-p)t^{-p}\phi\left( \frac{u^{2}+|\nabla u|^{2}}{2}t^{2}\right)(u^{2}+|\nabla u|^{2})dx\nonumber\\
&  & +\int_{\Omega}(4-p-q)t^{2-p}\phi'\left( \frac{u^{2}+|\nabla u|^{2}}{2}t^{2}\right)(u^{2}+|\nabla u|^{2})^{2}dx\nonumber\\
& &
+\int_{\Omega}t^{4-p}\phi''\left( \frac{u^{2}+|\nabla u|^{2}}{2}t^{2}\right)(u^{2}+|\nabla u|^{2})^{3}dx
\end{eqnarray}
and $\forall t>0$, $\eta'_{u}(t)<0$.}
\vskip2mm
\noindent
{\bf Proof.} Initially, in view of $(\phi_3)$, we have
\begin{eqnarray*}
\rho_{3}t^{1-p}\int_{\Omega}(u^{2}+|\nabla u|^{2})dx
\leq
\eta_{u}(t)
\leq
\rho_{4}t^{1-p}\int_{\Omega}(u^{2}+|\nabla u|^{2})dx,
\end{eqnarray*}
which together with $p>1$ implies that $\lim_{t\rightarrow0^{+}}\eta_{u}(t)=+\infty$ and $\lim_{t\rightarrow+\infty}\eta_{u}(t)=0$.
Next, we present the proof of (\ref{551}).
Actually, it follows that
\begin{eqnarray}\label{552}
& &\eta'_{u}(t)\nonumber\\
&= &
\lim_{h\rightarrow0}\frac{\eta_{u}(t+h)-\eta_{u}(t)}{h}\nonumber\\
&= &\lim_{h\rightarrow0}\int_{\Omega}\frac{1}{h}\Bigg[(1-q)(t+h)^{1-p}\phi\left( \frac{u^{2}+|\nabla u|^{2}}{2}(t+h)^{2}\right)
+(t+h)^{3-p}\phi'\left( \frac{u^{2}+|\nabla u|^{2}}{2}(t+h)^{2}\right)(u^{2}+|\nabla u|^{2})\nonumber\\
& &-(1-q)t^{1-p}\phi\left( \frac{u^{2}+|\nabla u|^{2}}{2}t^{2}\right)
-t^{3-p}\phi'\left( \frac{u^{2}+|\nabla u|^{2}}{2}t^{2}\right)(u^{2}+|\nabla u|^{2})\Bigg](u^{2}+|\nabla u|^{2})dx.
\end{eqnarray}
For any given $t>0$ and  $u\in H_{0}^{1}(\Omega)\backslash\{0\}$, let
\begin{eqnarray*}
f(h)=(1-q)(t+h)^{1-p}\phi\left( \frac{u^{2}+|\nabla u|^{2}}{2}(t+h)^{2}\right)
+(t+h)^{3-p}\phi'\left( \frac{u^{2}+|\nabla u|^{2}}{2}(t+h)^{2}\right)(u^{2}+|\nabla u|^{2}),
\;\;
0<|h|\leq1.
\end{eqnarray*}
Since $\phi\in C^{2}([0,+\infty),\mathbb{R})$, we obtain that $f$ is of class $C^{1}$ and
\begin{eqnarray*}
& & f'(h)\\
&= &
(1-q)(1-p)(t+h)^{-p}\phi\left( \frac{u^{2}+|\nabla u|^{2}}{2}(t+h)^{2}\right)
+(1-q)(t+h)^{2-p}\phi'\left( \frac{u^{2}+|\nabla u|^{2}}{2}(t+h)^{2}\right)(u^{2}+|\nabla u|^{2})\\
& &
+(3-p)(t+h)^{2-p}\phi'\left( \frac{u^{2}+|\nabla u|^{2}}{2}(t+h)^{2}\right)(u^{2}+|\nabla u|^{2})
+(t+h)^{4-p}\phi''\left( \frac{u^{2}+|\nabla u|^{2}}{2}(t+h)^{2}\right)(u^{2}+|\nabla u|^{2})^{2}.
\end{eqnarray*}
Applying the mean value theorem, there exists $\theta\in\mathbb{R}$ with $0<\theta<h\leq1$ such that
\begin{eqnarray}\label{548}
& &
\frac{1}{h}\Bigg[(1-q)(t+h)^{1-p}\phi\left( \frac{u^{2}+|\nabla u|^{2}}{2}(t+h)^{2}\right)
+(t+h)^{3-p}\phi'\left( \frac{u^{2}+|\nabla u|^{2}}{2}(t+h)^{2}\right)(u^{2}+|\nabla u|^{2})\nonumber\\
& &
-(1-q)t^{1-p}\phi\left( \frac{u^{2}+|\nabla u|^{2}}{2}t^{2}\right)
-t^{3-p}\phi'\left( \frac{u^{2}+|\nabla u|^{2}}{2}t^{2}\right)(u^{2}+|\nabla u|^{2})\Bigg](u^{2}+|\nabla u|^{2})\nonumber\\
&= &
 \frac{f(h)-f(0)}{h-0}(u^{2}+|\nabla u|^{2}) \nonumber\\
&= &
f'(\theta)(u^{2}+|\nabla u|^{2}) \nonumber\\
&= &
(1-q)(1-p)(t+\theta)^{-p}\phi\left( \frac{u^{2}+|\nabla u|^{2}}{2}(t+\theta)^{2}\right)(u^{2}+|\nabla u|^{2})
+(4-p-q)(t+\theta)^{2-p}\phi'\left( \frac{u^{2}+|\nabla u|^{2}}{2}(t+\theta)^{2}\right)\nonumber\\
& &
 \cdot(u^{2}+|\nabla u|^{2})^{2}
+(t+\theta)^{4-p}\phi''\left( \frac{u^{2}+|\nabla u|^{2}}{2}(t+\theta)^{2}\right)(u^{2}+|\nabla u|^{2})^{3}.
\end{eqnarray}
Next we show that there is $l\in L^{1}(\Omega)$ such that $|f'(\theta)|(u^{2}+|\nabla u|^{2})\leq l$.
In view of $0<\theta<h\leq1$, $0<q<1<p$, $(\phi_1)$, $(\phi_2)$ and Proposition 2.1, we have
\begin{eqnarray}\label{549}
& &|f'(\theta)|(u^{2}+|\nabla u|^{2})\nonumber\\
&\leq&
t^{-p}(u^{2}+|\nabla u|^{2})\cdot\Bigg[(1-q)(p-1)\phi\left( \frac{u^{2}+|\nabla u|^{2}}{2}(t+\theta)^{2}\right)\nonumber\\
& &+2|4-p-q|\left|\phi'\left( \frac{u^{2}+|\nabla u|^{2}}{2}(t+\theta)^{2}\right)\right|\frac{u^{2}+|\nabla u|^{2}}{2}(t+\theta)^{2}\nonumber\\
& &
+4\left|\phi''\left( \frac{u^{2}+|\nabla u|^{2}}{2}(t+\theta)^{2}\right)\right|\left(\frac{u^{2}+|\nabla u|^{2}}{2}(t+\theta)^{2}\right)^{2}\Bigg]\nonumber\\
&\leq&
t^{-p}(u^{2}+|\nabla u|^{2})\cdot\Bigg[(1-q)(p-1)\rho_{1}+\max\{2|4-p-q|, 4\}
\left|\phi'\left( \frac{u^{2}+|\nabla u|^{2}}{2}(t+\theta)^{2}\right)\right|\frac{u^{2}+|\nabla u|^{2}}{2}(t+\theta)^{2}\nonumber\\
& &
+\max\{2|4-p-q|, 4\}\left|\phi''\left( \frac{u^{2}+|\nabla u|^{2}}{2}(t+\theta)^{2}\right)\right|\left(\frac{u^{2}+|\nabla u|^{2}}{2}(t+\theta)^{2}\right)^{2}\Bigg]\nonumber\\
&\leq&
\left((1-q)(p-1)\rho_{1}+\rho_{2}\max\{2|4-p-q|, 4\}\right)t^{-p}(u^{2}+|\nabla u|^{2}):=l
\in L^{1}(\Omega).
 \end{eqnarray}
In addition, using the continuity of $\phi''$, we get
\begin{eqnarray}\label{550}
& &
\lim_{\theta\rightarrow0}f'(\theta)(u^{2}+|\nabla u|^{2})\nonumber\\
&= &
(1-q)(1-p)t^{-p}\phi\left( \frac{u^{2}+|\nabla u|^{2}}{2}t^{2}\right)(u^{2}+|\nabla u|^{2})
+(4-p-q)t^{2-p}\phi'\left( \frac{u^{2}+|\nabla u|^{2}}{2}t^{2}\right)\nonumber\\
& &
\cdot(u^{2}+|\nabla u|^{2})^{2}
+t^{4-p}\phi''\left( \frac{u^{2}+|\nabla u|^{2}}{2}t^{2}\right)(u^{2}+|\nabla u|^{2})^{3},
\;\mbox{a.e.}\; x \in \Omega.
 \end{eqnarray}
Note that $\theta\rightarrow0$ as $h\rightarrow0$.
Hence, it follows from (\ref{552}), (\ref{548}), (\ref{549}), (\ref{550}) and the Lebesgue dominated convergence theorem that (\ref{551}) holds.
Then using (\ref{551}) and $(\phi_5)$ we deduce that
\begin{eqnarray*}
\eta'_{u}(t)
\leq
-\rho_{6}t^{-p}\int_{\Omega}(u^{2}+|\nabla u|^{2})dx
<0,
\end{eqnarray*}
which implies that $\forall t>0$, $\eta'_{u}(t)<0$.
\qed

\vskip2mm
By Lemma 3.1, we obtain that the following conclusions of $m_{u}$ holds.
\vskip2mm
 \noindent
{\bf Lemma 3.2. } {\it Assume that $(\phi_1)$, $(\phi_2)$, $(\phi_3)$ and $(\phi_5)$ hold.
Then for any $\lambda>0$, $t\in(0,+\infty)$,  $u\in H_{0}^{1}(\Omega)\backslash\{0\}$, we have
\begin{eqnarray}\label{543}
m'_{u}(t)
&= & \int_{\Omega}\left[(1-q)t^{-q}\phi\left( \frac{u^{2}+|\nabla u|^{2}}{2}t^{2}\right)(u^{2}+|\nabla u|^{2})
+t^{2-q}\phi'\left( \frac{u^{2}+|\nabla u|^{2}}{2}t^{2}\right)(u^{2}+|\nabla u|^{2})^{2} \right]dx\nonumber\\
& &-(p-q)t^{p-q-1}\int_{\Omega}b(x)|u|^{p+1}dx,
\end{eqnarray}
and $m_{u}$ satisfies the following properties:
\begin{itemize}
 \item[$(i)$]
Assume that  $\int_{\Omega}b(x)|u|^{p+1}dx\leq0$ holds.
Then $m_{u}'(t)>0$ for any $t>0$. Moreover, $m_{u}(0):=\lim_{t\rightarrow0^{+}}m_{u}(t)=0$, $m_{u}(\infty):=\lim_{t\rightarrow+\infty}m_{u}(t)=+\infty$.
 \item[$(ii)$]
Suppose that $\int_{\Omega}b(x)|u|^{p+1}dx>0$ holds.
Then $m_{u}(0)=0$, $m_{u}(\infty)=-\infty$.
$m_{u}'(t)>0$ for any $t>0$ small enough. $m_{u}'(t)<0$ for any $t>0$ large enough.
Furthermore, there is an unique maximum point $\tilde{t}>0$ for $m_{u}$ such that $m_{u}'(\tilde{t})=0$.
\end{itemize}}
 \noindent
{\bf Proof.} Initially, we show that (\ref{543}) holds.
Let
\begin{eqnarray*}
\tilde{f}(h)=(t+h)^{1-q}\phi\left( \frac{u^{2}+|\nabla u|^{2}}{2}(t+h)^{2}\right),
\;\;
 0<|h|\leq1,
\;\;
 t\in(0,+\infty),
\;\;
u\in H_{0}^{1}(\Omega)/\{0\}.
 \end{eqnarray*}
By the continuity of $\phi'$, we get $\tilde{f}$ is of class $C^{1}$ and
\begin{eqnarray*}
\tilde{f}'(h)
= (1-q)(t+h)^{-q}\phi\left( \frac{u^{2}+|\nabla u|^{2}}{2}(t+h)^{2}\right)
+(t+h)^{2-q}\phi'\left( \frac{u^{2}+|\nabla u|^{2}}{2}(t+h)^{2}\right)(u^{2}+|\nabla u|^{2}).
 \end{eqnarray*}
Applying the mean value theorem, there is $\theta\in \mathbb{R}$ with $0<\theta<h\leq1$ such that
\begin{eqnarray}\label{540}
& &\frac{(t+h)^{1-q}\phi\left( \frac{u^{2}+|\nabla u|^{2}}{2}(t+h)^{2}\right)
-t^{1-q}\phi\left( \frac{u^{2}+|\nabla u|^{2}}{2}t^{2}\right)}{h}(u^{2}+|\nabla u|^{2})\nonumber\\
&= &\frac{\tilde{f}(h)-\tilde{f}(0)}{h-0}(u^{2}+|\nabla u|^{2})\nonumber\\
&= &\tilde{f}'(\theta)(u^{2}+|\nabla u|^{2}).
 \end{eqnarray}
Next we show that there is $h\in L^{1}(\Omega)$ such that $|\tilde{f}'(\theta)|(u^{2}+|\nabla u|^{2})\leq h$.
Taking into account that $0<\theta<h\leq1$, $0<q<1$, $(\phi_1)$, $(\phi_2)$ and Proposition 2.1, we have
\begin{eqnarray}\label{541}
|\tilde{f}'(\theta)|(u^{2}+|\nabla u|^{2})
\leq
[(1-q)\rho_{1}+2\rho_{2}]t^{-q}(u^{2}+|\nabla u|^{2})
\in L^{1}(\Omega).
 \end{eqnarray}
Additionally, by the continuity of $\phi$, we conclude that
\begin{eqnarray}\label{542}
\lim_{\theta\rightarrow 0}\tilde{f}'(\theta)(u^{2}+|\nabla u|^{2})
=  (1-q)t^{-q}\phi\left( \frac{u^{2}+|\nabla u|^{2}}{2}t^{2}\right)(u^{2}+|\nabla u|^{2})
+t^{2-q}\phi'\left( \frac{u^{2}+|\nabla u|^{2}}{2}t^{2}\right)(u^{2}+|\nabla u|^{2})^{2},
 \end{eqnarray}
for a.e. $x\in \Omega$. Note that $\theta\rightarrow0$ as $h\rightarrow0$.
Thus, in view of (\ref{540}), (\ref{541}), (\ref{542}) and the Lebesgue dominated convergence theorem, we obtain that
\begin{eqnarray*}
& &\lim_{h\rightarrow 0}\int_{\Omega}\frac{(t+h)^{1-q}\phi\left( \frac{u^{2}+|\nabla u|^{2}}{2}(t+h)^{2}\right)
-t^{1-q}\phi\left( \frac{u^{2}+|\nabla u|^{2}}{2}t^{2}\right)}{h}(u^{2}+|\nabla u|^{2})dx\\
&= & \int_{\Omega}\left[(1-q)t^{-q}\phi\left( \frac{u^{2}+|\nabla u|^{2}}{2}t^{2}\right)(u^{2}+|\nabla u|^{2})
+t^{2-q}\phi'\left( \frac{u^{2}+|\nabla u|^{2}}{2}t^{2}\right)(u^{2}+|\nabla u|^{2})^{2}\right]dx,
 \end{eqnarray*}
which implies that (\ref{543}) holds.
Next we shall prove the item (i).
For any $t>0$,
in view of $q<p$ and $(\phi_3)$, we mention that
\begin{eqnarray*}
m'_{u}(t)
\geq t^{-q}\int_{\Omega}\rho_{3}(u^{2}+|\nabla u|^{2})dx
-(p-q)t^{p-q-1}\int_{\Omega}b(x)|u|^{p+1}dx
>0.
\end{eqnarray*}
Moreover, taking into account $(\phi_1)$, it is easy to verify that
\begin{eqnarray}\label{544}
t^{1-q}\rho_{0}\int_{\Omega}(u^{2}+|\nabla u|^{2})dx
-t^{p-q}\int_{\Omega}b(x)|u|^{p+1}dx
\leq
 m_{u}(t)
\end{eqnarray}
and
\begin{eqnarray}\label{545}
m_{u}(t)
\leq
t^{1-q}\rho_{1}\int_{\Omega}(u^{2}+|\nabla u|^{2})dx
-t^{p-q}\int_{\Omega}b(x)|u|^{p+1}dx.
\end{eqnarray}
It follows from (\ref{544}), (\ref{545}) and $0<q<1<p$ that $\lim_{t\rightarrow0^{+}}m_{u}(t)=0$ and $\lim_{t\rightarrow+\infty}m_{u}(t)=+\infty$.
Next, we  prove the item (ii).
Firstly, in view of (\ref{544}), (\ref{545}) and $0<q<1<p$, we mention that $\lim_{t\rightarrow0^{+}}m_{u}(t)=0$ and $\lim_{t\rightarrow+\infty}m_{u}(t)=-\infty$.
Then by $0<q<1<p$ and $(\phi_3)$, we have
\begin{eqnarray}\label{546}
m'_{u}(t)
\geq
\rho_{3} t^{-q}\int_{\Omega}(u^{2}+|\nabla u|^{2})dx
-(p-q)t^{p-q-1}\int_{\Omega}b(x)|u|^{p+1}dx
\end{eqnarray}
and
\begin{eqnarray}\label{547}
m'_{u}(t)
\leq
\rho_{4} t^{-q}\int_{\Omega}(u^{2}+|\nabla u|^{2})dx
-(p-q)t^{p-q-1}\int_{\Omega}b(x)|u|^{p+1}dx.
\end{eqnarray}
Note that $q>q+1-p$. Then it follows from (\ref{546}) that $m_{u}'(t)>0$ for any $t>0$ small enough, and
(\ref{547}) implies that $m_{u}'(t)<0$ for any $t>0$ large enough.
Thus, there exists $\tilde{t}>0$ such that $m'_{u}(\tilde{t})=0$ and $\tilde{t}>0$ is a global maximum point for $m_{u}$.
Since $m'_{u}(t)=0$ if and only if
\begin{eqnarray*}
 \int_{\Omega}t^{1-p}(u^{2}+|\nabla u|^{2})\left[(1-q)\phi\left( \frac{u^{2}+|\nabla u|^{2}}{2}t^{2}\right)
+2\phi'\left( \frac{u^{2}+|\nabla u|^{2}}{2}t^{2}\right) \frac{u^{2}+|\nabla u|^{2}}{2}t^{2}\right]dx
=
(p-q)\int_{\Omega}b(x)|u|^{p+1}dx,
\end{eqnarray*}
we get the unique $\tilde{t}>0$ by Lemma 3.1.
Hence, there is an unique maximum point $\tilde{t}>0$ for $m_{u}(t)$ such that $m_{u}'(\tilde{t})=0$.
\qed

\vskip2mm
\par
It follows from the above lemma that $m_{u}$ is a strictly increasing function for $t>0$ whenever $\int_{\Omega}b(x)|u|^{p+1}dx\leq0$
and $m_{u}$ is initially increasing and eventually decreasing  when $\int_{\Omega}b(x)|u|^{p+1}dx>0$.
Furthermore, it is easy to see that the essential nature of $m_{u}$ is determined by the sign of $\int_{\Omega}b(x)|u|^{p+1}dx$.
More precisely, it follows from (\ref{4.0.1}) that the essential nature of the fibering map $\gamma_{u}$ is determined by the signs of $\int_{\Omega}b(x)|u|^{p+1}dx$ and $\int_{\Omega}\lambda a(x)|u|^{q+1}dx$.
Thus, we shall give a fairly complete description of the fibering map $\gamma_{u}$ for every possible choice of signs of $\int_{\Omega}b(x)|u|^{p+1}dx$ and $\int_{\Omega}\lambda a(x)|u|^{q+1}dx$.
\par
 $\clubsuit^{1}$ Assume that $\int_{\Omega}b(x)|u|^{p+1}dx\leq0$ and $\int_{\Omega}\lambda a(x)|u|^{q+1}dx\leq0$ hold.
By Lemma 3.2.(i), it is easy to verify that
$m_{u}(0):=\lim_{t\rightarrow0^{+}}m_{u}(t)=0$,
$m_{u}(\infty):=\lim_{t\rightarrow+\infty}m_{u}(t)=+\infty$
and
$m_{u}'(t)>0$, for all $t>0$.
Under these results we mention that
\begin{eqnarray*}
m_{u}(t)>0\ge \int_{\Omega}\lambda a(x)|u|^{q+1}dx\;\mbox{for any }\;t>0, \lambda>0.
\end{eqnarray*}
Thus $tu\not\in N_{\lambda}(\Omega)$ for any $t>0$.
In particular, we also see that $\gamma_{u}'(t)>0$ for any $t>0$.
\par
 $\clubsuit^{2}$ Assume that $\int_{\Omega}b(x)|u|^{p+1}dx\leq0$ and $\int_{\Omega}\lambda a(x)|u|^{q+1}dx>0$ hold.
Using one more time Lemma 3.2.(i), we obtain that $m_{u}(0)=0$, $m_{u}(\infty)=+\infty$ and $m_{u}'(t)>0$ for any $t>0$.
Thus, the equation $m_{u}(t)=\int_{\Omega}\lambda a(x)|u|^{q+1}dx$ admits exactly one solution $t_{1}=t_{1}(u,\lambda)>0$
such that $\gamma_{u}'(t_{1})=0$ and $t_{1}u\in N_{\lambda}(\Omega)$ for any $\lambda>0$.
Firstly, $\gamma_{u}'(t)<0$ for any $t\in (0,t_{1})$ and $\gamma_{u}'(t)>0$ for any $t\in (t_{1}, +\infty)$.
Moreover, by (\ref{2.0.4}), $(\phi_{1})$ and $1<q+1<2<p+1$, we have
$\lim_{t\rightarrow0^{+}}\gamma_{u}(t)=0$ and $\lim_{t\rightarrow+\infty}\gamma_{u}(t)=+\infty$.
Hence, $t_{1}$ is a global minimum point of $\gamma_{u}(t)$ with $\gamma_{u}(t_{1})<0$.
Additionally, it follows from (\ref{2.6}) and (\ref{543}) that
\begin{eqnarray}\label{560}
\gamma''_{tu}(1)=t^{q+2}m'_{u}(t),\;\;t>0,
\end{eqnarray}
which together with $m_{u}'(t_{1})>0$ imply that $\gamma''_{t_{1}u}(1)>0$ and then $t_{1}u\in N^{+}_{\lambda}$.
\par
 $\clubsuit^{3}$  Assume that $\int_{\Omega}b(x)|u|^{p+1}dx>0$ and $\int_{\Omega}\lambda a(x)|u|^{q+1}dx\leq0$ hold.
By Lemma 3.2.(ii), we know that $m_{u}(0)=0$, $m_{u}(\infty)=-\infty$,
$m_{u}'(t)>0$ for any $t>0$ small enough and $m_{u}'(t)<0$ for any $t>0$ large enough.
Furthermore, $m_{u}(t)$ has an unique critical point $\tilde{t}>0$ which is a local maximum.
As a consequence, the equation $m_{u}(t)=\int_{\Omega}\lambda a(x)|u|^{q+1}dx$ admits one solution $0<\tilde{t}<t_{2}$
such that $\gamma_{u}'(t_{2})=0$ and $t_{2}u\in N_{\lambda}(\Omega)$ for any $\lambda>0$.
First of all, $\gamma_{u}'(t)>0$ for any $t\in (0,t_{2})$ and $\gamma_{u}'(t)<0$ for any $t\in (t_{2}, +\infty)$.
Due to the fact that (\ref{2.0.4}), $(\phi_{1})$ and $1<q+1<2<p+1$, we deduce that
$\lim_{t\rightarrow+\infty}\gamma_{u}(t)=-\infty$.
Hence, $t_{2}$ is a global maximum point of $\gamma_{u}(t)$ with $\gamma_{u}(t_{2})>0$.
Additionally, it follows from (\ref{560}) and $m_{u}'(t_{2})<0$ that $\gamma''_{t_{2}u}(1)<0$ and $t_{2}u\in N^{-}_{\lambda}$.
\par
 $\clubsuit^{4}$  Assume that $\int_{\Omega}b(x)|u|^{p+1}dx>0$ and $\int_{\Omega}\lambda a(x)|u|^{q+1}dx>0$ hold.
It is easily seen that the behavior of $\gamma_{u}$ depends on $\lambda$.
\par
If $\lambda>\frac{m_{u}(\tilde{t})}{\int_{\Omega}\lambda a(x)|u|^{q+1}dx}$,
it follows that
\begin{eqnarray*}
m_{u}(\tilde{t})<\int_{\Omega}\lambda a(x)|u|^{q+1}dx,
\end{eqnarray*}
which implies that $\gamma_{u}$ has no critical points.
Moreover, by (\ref{2.0.4}), $(\phi_{1})$ and $1<q+1<2<p+1$, we have
$\lim_{t\rightarrow0^{+}}\gamma_{u}(t)=0$ and $\lim_{t\rightarrow+\infty}\gamma_{u}(t)=-\infty$.
Hence, $\gamma_{u}(t)<0$ for any $t\in (0, +\infty)$ i.e., $\gamma_{u}$ is a decreasing function.
\par
If $\lambda=\frac{m_{u}(\tilde{t})}{\int_{\Omega}\lambda a(x)|u|^{q+1}dx}$,
we have
\begin{eqnarray*}
m_{u}(\tilde{t})=\lambda \int_{\Omega}a(x)|u|^{q+1}dx,
\end{eqnarray*}
which implies that the fibering map $\gamma_{u}$ admits an unique critical point $\tilde{t}>0$.
$\gamma'_{u}(t)<0$ for any $t\in (0,\tilde{t})\cup (\tilde{t}, +\infty)$.
It follows from (\ref{2.0.4}), $(\phi_{1})$ and $1<q+1<2<p+1$ that
$\lim_{t\rightarrow0^{+}}\gamma_{u}(t)=0$ and $\lim_{t\rightarrow+\infty}\gamma_{u}(t)=-\infty$.
\par
If $\lambda<\frac{m_{u}(\tilde{t})}{\int_{\Omega}\lambda a(x)|u|^{q+1}dx}$,
it follows that
\begin{eqnarray*}
m_{u}(\tilde{t})>\int_{\Omega}\lambda a(x)|u|^{q+1}dx,
\end{eqnarray*}
which implies that $m_{u}(t)=\int_{\Omega}\lambda a(x)|u|^{q+1}dx$ has two solutions $t_{3}(u,\lambda)$ and $t_{4}(u,\lambda)$ with $0<t_{3}(u,\lambda)<\tilde{t}<t_{4}(u,\lambda)$
such that $\gamma_{u}'(t_{3})=\gamma_{u}'(t_{4})=0$ and $t_{3}u, t_{4}u\in N_{\lambda}(\Omega)$.
Firstly, $\gamma'_{u}(t)<0$ for any $t\in (0,t_{3})$, $\gamma'_{u}(t)>0$ for any $t\in (t_{3},t_{4})$ and $\gamma'_{u}(t)<0$ for any $(t_{4},+\infty)$.
Then by (\ref{2.0.4}), $(\phi_{1})$ and $1<q+1<2<p+1$, we mention that
$\lim_{t\rightarrow+\infty}\gamma_{u}(t)=-\infty$.
Therefore, $t_{3}$ is a local minimum point of $\gamma_{u}$ and $t_{4}$ is a global maximum point of $\gamma_{u}$.
Furthermore, it follows from $m'_{u}(t_{3})>0$, $m'_{u}(t_{4})<0$ and (\ref{560}) that $t_{3}u\in N^{+}_{\lambda}$ and $t_{4}u\in N^{-}_{\lambda}$.
The possible distributions for $\gamma_{u}(t_{3})$ and $\gamma_{u}(t_{4})$ are
$\gamma_{u}(t_{3})<0<\gamma_{u}(t_{4})$,
$\gamma_{u}(t_{3})<\gamma_{u}(t_{4})<0$,
or
$\gamma_{u}(t_{3})<\gamma_{u}(t_{4})=0$.
\vskip2mm
 \noindent
{\bf Remark 3.1.} We expect find the critical points of the extremal type of the functional $J_{\lambda}$ by standard minimization procedure in each part $N_{\lambda}^{\pm}$.
To this end, it is necessary to estimate the energy of the functional $J_{\lambda}$ in each part $N_{\lambda}^{\pm}$.
To estimate the energy of the functional $J_{\lambda}$ in $N_{\lambda}^{-}$,
it is necessary to show that $\gamma_{u}(t_{3})<0<\gamma_{u}(t_{4})$.
To this end, we shall give a upper bounded $\lambda_{2}$ of $\lambda$ such that
for any $u\in H_{0}^{1}(\Omega)\backslash \{0\}$ satisfies $\int_{\Omega}b(x)|u|^{p+1}dx>0$,
there exists $t_{\max}>0$ and $\delta_{\lambda}>0$,
for any $\lambda\in (0,\lambda_{2})$,
we have $\gamma_{u}(t_{\max})\geq \delta_{\lambda}$.
The definitions of $\lambda_{2}$ and $\delta_{\lambda}$ will be given in Lemma 3.3.
\vskip2mm
 \noindent
{\bf Lemma 3.3. } {\it Suppose that $(\phi_{1})$ and $(H)$ hold.
Then for any $\lambda \in(0,\lambda_{2})$ and $u\in H_{0}^{1}(\Omega)\backslash\{0\}$ which satisfies $\int_{\Omega}b(x)|u|^{p+1}dx>0$,
there exists $t_{\max}>0$ such that $\gamma_{u}(t_{\max})\geq\delta_{\lambda}>0$, where
$
 \delta_{\lambda} =\delta^{\frac{q+1}{2}}\left(\delta^{1-\frac{q+1}{2}}-\lambda c_{1}\right)
$
and $\lambda_2$, $c_1$ and $\delta$ are defined by (\ref{aa1}) and (\ref{aa2}). }
\vskip0mm
\noindent
{\bf Proof.} Let $\lambda\in (0,\lambda_{2})$, $u\in H_{0}^{1}(\Omega)\backslash \{0\}$ which satisfies $\int_{\Omega}b(x)|u|^{p+1}dx>0$.
Define the auxiliary function $h_{u}:\mathbb{R}\rightarrow\mathbb{R}$ given by
\begin{eqnarray}\label{4.4.1}
 h_{u}(t)
 =
 \int_{\Omega}\Phi\left(\frac{u^{2}+|\nabla u|^{2}}{2}t^{2}\right)dx
 -\frac{1}{p+1}t^{p+1}\int_{\Omega}b(x)|u|^{p+1}dx.
\end{eqnarray}
Firstly, we shall prove that there exists a maximum point $t_{\max}>0$ for $h_{u}(t)$ such that $ h'_{u}(t_{\max})=0$.
Due to $(\phi_1)$, we obtain that
\begin{eqnarray*}
 h_{u}(t)
 \leq
 \frac{\rho_{1}}{2}t^{2}\int_{\Omega}(u^{2}+|\nabla u|^{2})dx
 -\frac{1}{p+1}t^{p+1}\int_{\Omega}b(x)|u|^{p+1}dx,
\end{eqnarray*}
which together with $2<p+1$ implies that
\begin{eqnarray}\label{561}
   \lim_{t\rightarrow+\infty}h_{u}(t)=-\infty.
\end{eqnarray}
Then it follows that
\begin{eqnarray*}
 h'_{u}(t)
 =
 t\int_{\Omega}\phi\left( \frac{u^{2}+|\nabla u|^{2}}{2}t^{2}\right)(u^{2}+|\nabla u|^{2})dx
 -t^{p}\int_{\Omega}b(x)|u|^{p+1}dx,
 \;\;t>0.
\end{eqnarray*}
In view of $(\phi_1)$, we have
\begin{eqnarray}\label{562}
 h'_{u}(t)
 \geq
 t\rho_{0}\int_{\Omega}(u^{2}+|\nabla u|^{2})dx
 -t^{p}\int_{\Omega}b(x)|u|^{p+1}dx
\end{eqnarray}
and
\begin{eqnarray}\label{563}
  h'_{u}(t)
  \leq
  t\rho_{1}\int_{\Omega}(u^{2}+|\nabla u|^{2})
  -t^{p}\int_{\Omega}b(x)|u|^{p+1}dx.
\end{eqnarray}
It follows from (\ref{562}), (\ref{563}) and $1<p$ that $h'_{u}(t)>0$ for any $t>0$ small enough and $h'_{u}(t)<0$ for any $t>0$ large enough.
At this moment, taking into account the estimates above and (\ref{561}) we easily see that there exists $t_{\max}>0$ such that $ h'_{u}(t_{\max})=0$ and
$t_{\max}$ is a global maximum point for $h_{u}$.
Next, we prove that
\begin{eqnarray}\label{564}
 \frac{\lambda}{q+1}t_{\max}^{q+1}\int_{\Omega}a(x)|u|^{q+1}dx
\leq
 \frac{\lambda\|a\|_{\infty}S_{q+1}^{q+1}}{(q+1)(\frac{\rho_{0}}{2}-\frac{\rho_{1}}{p+1})^{\frac{q+1}{2}}}h_{u}^{\frac{q+1}{2}}(t_{\max}).
\end{eqnarray}
In fact, by $q+1<2^{\ast}$, $(H)$ and Proposition 2.1, we have
\begin{eqnarray}\label{565}
 \frac{\lambda}{q+1}t_{\max}^{q+1}\int_{\Omega}a(x)|u|^{q+1}dx
 &\leq&\frac{\lambda}{q+1}t_{\max}^{q+1}\|a\|_{\infty}\|u\|_{q+1}^{q+1}\nonumber\\
 &\leq&\frac{\lambda}{q+1}\|a\|_{\infty}S_{q+1}^{q+1}\|t_{\max}\nabla u\|_{2}^{q+1}.
\end{eqnarray}
Note that for any $t>0$, $h'_{u}(t)=0$ if and only if
$t\int_{\Omega}\phi\left( \frac{u^{2}+|\nabla u|^{2}}{2}t^{2}\right)(u^{2}+|\nabla u|^{2})dx=t^{p}\int_{\Omega}b(x)|u|^{p+1}dx$,
and so
\begin{eqnarray}\label{566}
\frac{1}{p+1}t_{\max}^{p+1}\int_{\Omega}b(x)|u|^{p+1}dx
=
\frac{2}{p+1}\int_{\Omega}\phi\left( \frac{u^{2}+|\nabla u|^{2}}{2}t_{\max}^{2}\right)\frac{u^{2}+|\nabla u|^{2}}{2}t_{\max}^{2}dx,
\end{eqnarray}
which together with $(\phi_1)$ shows that
\begin{eqnarray}\label{567}
 h_{u}(t_{\max})
 & \geq &
 \left(\frac{\rho_{0}}{2}-\frac{\rho_{1}}{p+1}\right)t_{\max}^{2}(\|u\|_{2}^{2}+\|\nabla u\|_{2}^{2})\nonumber\\
 & \geq &
\left(\frac{\rho_{0}}{2}-\frac{\rho_{1}}{p+1}\right)\|t_{\max}\nabla u\|_{2}^{2}.
\end{eqnarray}
It follows from (\ref{565}) and (\ref{567}) that (\ref{564}) holds.
Next, we prove that
\begin{eqnarray}\label{568}
 h_{u}(t_{\max})
\geq
\delta>0,
\;\;\forall u\in H^{1}_{0}(\Omega)\backslash\{0\},
\end{eqnarray}
where
\begin{eqnarray*}
 \delta=\left(\frac{\rho_{0}}{2}-\frac{\rho_{1}}{p+1}\right)\left(  \frac{\rho_{0}}{\|b\|_{\infty}S_{p+1}^{p+1}} \right)^{\frac{2}{p-1}}.
\end{eqnarray*}
Firstly, in view of $p+1<2^{\ast}$, $(H)$ and Proposition 2.1, we conclude that
\begin{eqnarray}\label{569}
 \frac{1}{p+1}t_{\max}^{p+1}\int_{\Omega}b(x)|u|^{p+1}dx
 &\leq&
 \frac{1}{p+1}t_{\max}^{p+1}\|b\|_{\infty}\|u\|_{p+1}^{p+1}\nonumber \\
 &\leq&
 \frac{1}{p+1}\|b\|_{\infty}S_{p+1}^{p+1}\|t_{\max}\nabla u\|_{2}^{p+1}.
\end{eqnarray}
Then by (\ref{566}) and $(\phi_1)$, we have
\begin{eqnarray*}
\frac{1}{p+1}t_{\max}^{p+1}\int_{\Omega}b(x)|u|^{p+1}dx
 \geq
 \frac{\rho_{0}}{p+1}\|t_{\max}\nabla u\|_{2}^{2},
\end{eqnarray*}
which together with (\ref{569}) implies that
\begin{eqnarray}\label{570}
\| t_{\max}\nabla u\|_{2}
 \geq
 \left(  \frac{\rho_{0}}{\|b\|_{\infty}S_{p+1}^{p+1}} \right)^{\frac{1}{p-1}}.
\end{eqnarray}
Hence, combining with (\ref{567}), we have (\ref{568}) holds.
Now, in view of (\ref{2.0.4}), (\ref{4.4.1}), (\ref{564}) and (\ref{568}), we observe that
\begin{eqnarray*}
 \gamma_{u}(t_{\max})
 & = &
 h_{u}(t_{\max})-\frac{\lambda}{q+1}t_{\max}^{q+1}\int_{\Omega}a(x)|u|^{q+1}dx\\
 & \geq &
 h_{u}(t_{\max})
- \frac{\lambda\|a\|_{\infty}S_{q+1}^{q+1}}{(q+1)(\frac{\rho_{0}}{2}-\frac{\rho_{1}}{p+1})^{\frac{q+1}{2}}}h_{u}^{\frac{q+1}{2}}(t_{\max})\\
 & =&
 h_{u}^{\frac{q+1}{2}}(t_{\max})\left[h_{u}^{1-\frac{q+1}{2}}(t_{\max})-\lambda c_{1}\right]\\
 & \geq &
 \delta^{\frac{q+1}{2}}\left[\delta^{1-\frac{q+1}{2}}-\lambda c_{1}\right]
 = \delta_{\lambda}
 > 0
 \end{eqnarray*}
when $0<\lambda<\lambda_{2}$.
\qed
\vskip2mm
\par
Now, we fix $\lambda \in (0,\lambda_{0})$, where $\lambda_{0}=\min\{\lambda_{1},\lambda_{2}\}$. Then by Lemma 2.1, $N_{\lambda}=N_{\lambda}^{+}\cup N_{\lambda}^{-}$.
\vskip2mm
 \noindent
{\bf Lemma 3.4. } {\it Assume that $(\phi_{1})$, $(\phi_{3})$, $(\phi_{4})$ and $(H)$ hold.
Then for any $0<\lambda<\lambda_{0}$, the following conclusions hold:
\begin{itemize}
\item[$(i)$]
For all $ u \in N_{\lambda}^{+}$, $\int_{\Omega}a(x)|u|^{q+1}dx>0$, $\gamma_{u}'(t)<0$ for any  $t\in (0,1)$ and $\gamma_{u}(1)<0$. Furthermore, $J_{\lambda}(u)<0$.
 \item[$(ii)$]
For any $ u \in N_{\lambda}^{-}$, $\int_{\Omega}b(x)|u|^{p+1}dx>0$, and
$t=1$ is a global maximum point for $\gamma_{u}$ such that $\gamma_{u}(1)\geq  \delta_{\lambda}>0$.
Furthermore, $J_{\lambda}(u)\geq  \delta_{\lambda}>0$.
\end{itemize}}
 \noindent{\bf Proof.} Initially, we prove the item (i). Fixed $ u \in N_{\lambda}^{+}$.
 By (\ref{2.7}) and $(\phi_4)$, we have
\begin{eqnarray*}
0
<
\gamma_{u}''(1)
 \leq
-\rho_{5}\int_{\Omega}(u^{2}+|\nabla u|^{2})dx
+\lambda(p -q)\int_{\Omega}a(x)|u|^{q+1}dx,
 \end{eqnarray*}
which together with $q<p$ implies that $\int_{\Omega}a(x)|u|^{q+1}dx>0$.
Additionally, $\gamma_{u}'(1)=0$ and $\gamma_{u}''(1)>0$ show that $t=1$ is a local minimum point of $\gamma_{u}$.
Hence, the following conclusions can be obtained from $\int_{\Omega}a(x)|u|^{q+1}dx>0$ and the previous analysis of fibering mapping $\gamma_{u}$:
\par
$\bullet$ Assume that $\int_{\Omega}b(x)|u|^{p+1}dx\leq0$ holds.
Then $t=1$ is an unique critical point of $\gamma_{u}$.
First of all, $\gamma_{u}'(t)<0$ for any $t\in (0,1)$ and $\gamma_{u}'(t)>0$ for any $t\in (1, +\infty)$.
Then, $\lim_{t\rightarrow0^{+}}\gamma_{u}(t)=0$, $\lim_{t\rightarrow+\infty}\gamma_{u}(t)=+\infty$.
Hence, $t=1$ is a global minimum point of $\gamma_{u}(t)$ with $\gamma_{u}(1)<0$.
\par
$\bullet$  Assume that $\int_{\Omega}b(x)|u|^{p+1}dx>0$ holds.
Then there exists $t_{4}(u,\lambda)>0$ with $0<1<\tilde{t}<t_{4}(u,\lambda)$ such that $\gamma_{u}'(t_{4})=0$ and $t_{4}u\in N^{-}_{\lambda}$.
Firstly, $\gamma_{u}(t)$ is decreasing in $(0, 1)$, increasing in $(1,t_{4})$ and decreasing in $t\in (t_{4},+\infty)$.
Then, it follows that $\lim_{t\rightarrow+\infty}\gamma_{u}(t)=-\infty$.
Thus $t=1$ is a local minimum point of $\gamma_{u}$, $t_{4}$ is a global maximum point of $\gamma_{u}$ and $\gamma_{u}(1)<0<\gamma_{u}(t_{4})$.
\par
As a consequence, $J_{\lambda}(u)=\gamma_{u}(1)<0$ for any $ u \in N_{\lambda}^{+}$.
Next, we prove the item (ii).
Let $ u \in N_{\lambda}^{-}$.
In view of (\ref{2.6}) and $(\phi_3)$, we have
\begin{eqnarray*}
0
>
\gamma_{u}''(1)
 \geq
\rho_{3}\int_{\Omega}(u^{2}+|\nabla u|^{2})dx
-(p-q)\int_{\Omega}b(x)|u|^{p+1}dx,
 \end{eqnarray*}
which together with $q<p$ implies that $\int_{\Omega}b(x)|u|^{p+1}dx>0$.
Furthermore, the facts that $\gamma_{u}'(1)=0$ and $\gamma_{u}''(1)<0$ imply that $t=1$ is a local maximum point of $\gamma_{u}$.
Then the following conclusions can be obtained from $\int_{\Omega}b(x)|u|^{p+1}dx>0$ and the previous analysis of fibering mapping $\gamma_{u}$:
\par
$\bullet$ Assume that $\int_{\Omega}a(x)|u|^{q+1}dx\leq0$ holds.
Then $t=1$ is an unique critical point of $\gamma_{u}$.
Initially, $\gamma_{u}'(t)>0$ for any $t\in (0,1)$ and $\gamma_{u}'(t)<0$ for any $t\in (1, +\infty)$.
Then it follows that $\lim_{t\rightarrow+\infty}\gamma_{u}(t)=-\infty$.
Thus, $t=1$ is a global maximum point of $\gamma_{u}(t)$ with $\gamma_{u}(1)>0$.
\par
$\bullet$ Assume that $\int_{\Omega}a(x)|u|^{q+1}dx>0$ holds.
Then there exists $t_{3}(u,\lambda)>0$ with $0<t_{3}(u,\lambda)<\tilde{t}<1$ such that $t_{3}u\in N^{+}_{\lambda}$.
Firstly, $\gamma_{u}(t)$ is decreasing in $(0, t_{3})$, increasing in $(t_{3},1)$ and decreasing in $(1,+\infty)$.
Then it follows that $\lim_{t\rightarrow+\infty}\gamma_{u}(t)=-\infty$.
Hence, $t_{3}$ is a local minimum point of $\gamma_{u}$, $t=1$ is a global maximum point of $\gamma_{u}$ and then by Remark 3.1 and Lemma 3.3, we have $\gamma_{u}(t_{3})<0<\gamma_{u}(1)$.
\par
As a consequence, we obtain that $t=1$ is a golbal maximum point for $\gamma_{u}$.
Then by Lemma 3.3, we have that $\gamma_{u}(1)\geq\delta_{\lambda}>0$.
\qed
\section{The Palais-Smale sequence}
\par
In this section, in order to prove that any minimizer sequences on the Nehari manifold $N_{\lambda}^{+}$ or $N_{\lambda}^{-}$ provides us a Palais-Smale sequences, we follow the ideas discussed in \cite{Tarantello1992} and obtain the following result.
\vskip2mm
 \noindent
{\bf Proposition 4.1. } {\it Suppose that $(\phi_{1})$-$(\phi_{5})$ and $(H)$ hold.
Then for any $\lambda \in (0,\lambda_{0})$, we have the following assertions:
\begin{itemize}
\item[$(i)$] there exists a minimizing sequence $\{u_{n}\}\subset N_{\lambda}^{+}$ such that
\begin{eqnarray*}
\lim_{n\rightarrow\infty}J_{\lambda}(u_{n})= \inf_{u\in N_{\lambda}^{+}}J_{\lambda}(u)<0,
\;\;
\lim_{n\rightarrow\infty}J'_{\lambda}(u_{n})=0;
\end{eqnarray*}
\item[$(ii)$] there exists a minimizing sequence $\{u_{n}\}\subset N_{\lambda}^{-}$ such that
\begin{eqnarray*}
\lim_{n\rightarrow\infty}J_{\lambda}(u_{n})= \inf_{u\in N_{\lambda}^{-}}J_{\lambda}(u)>0,
\;\;
\lim_{n\rightarrow\infty}J'_{\lambda}(u_{n})=0.
\end{eqnarray*}
\end{itemize}}
\par
Initially, we shall prove some auxiliary results in order to prove Proposition 4.1.
\vskip2mm
 \noindent
{\bf Lemma 4.2. } {\it Suppose that $(\phi_{1})$, $(\phi_{2})$ and $(H)$ hold.
Then for each $\lambda \in (0,\lambda_{0})$ and $u\in N_{\lambda}^{+}$, there exist $\epsilon>0$ and a differentiable function $\xi^{+}:B(0,\epsilon)\subset H^{1}_{0}(\Omega) \rightarrow \mathbb{R}^{+}$ such that $\xi^{+}(0)=1$, the function
$\xi^{+}(v)(u-v)\in N_{\lambda}^{+}$ and
\begin{eqnarray}\label{4.1.1}
 \langle (\xi^{+})'(0),v\rangle
& =&\frac{1}{\gamma_{u}''(1)}\int_{\Omega}
\Bigg[2\phi\left(\frac{u^{2}+|\nabla u|^{2}}{2}\right)(uv+\nabla u\nabla v)
 +\phi'\left(\frac{u^{2}+|\nabla u|^{2}}{2}\right)(uv+\nabla u\nabla v)(u^{2}+|\nabla u|^{2})
dx\nonumber\\
& &
-(q+1)\lambda a(x)|u|^{q-1}uv
-(p+1)b(x)|u|^{p-1}uv\Bigg]dx
\end{eqnarray}
for all $v\in B(0,\epsilon)$.\\}
{\bf Proof.} Let $\lambda \in (0,\lambda_{0})$ be fixed. For $u\in N_{\lambda}^{+}$, define a function $F_{u}:\mathbb{R}^{+}\times H^{1}_{0}(\Omega)\setminus\{0\}\rightarrow\mathbb{R}$ by
\begin{eqnarray}\label{4.1.2}
 F_{u}(\xi,w)
& =&
\langle J'_{\lambda}(\xi(u-w)),\xi(u-w)\rangle\nonumber\\
& =&
\int_{\Omega}\xi^{2}\phi\left(\frac{(u-w)^{2}+|\nabla (u-w)|^{2}}{2}\xi^{2}\right)\left((u-w)^{2}+|\nabla (u-w)|^{2}\right)dx\nonumber\\
& &
-\xi^{q+1}\lambda \int_{\Omega}a(x)|u-w|^{q+1}dx
-\xi^{p+1}\int_{\Omega}b(x)|u-w|^{p+1}dx
\end{eqnarray}
for all $\xi\in \mathbb{R}^{+}$ and $w\in H^{1}_{0}(\Omega)\setminus\{0\}$.
Then $F_{u}(1,0)=\langle J'_{\lambda}(u),u\rangle=0$ and
\begin{eqnarray*}
& &
\partial_{1}F_{u}(\xi,w)\\
& =&
\int_{\Omega}\Bigg[2\xi\phi\left(\frac{(u-w)^{2}+|\nabla (u-w)|^{2}}{2}\xi^{2}\right)\left((u-w)^{2}+|\nabla (u-w)|^{2}\right)+\xi^{2}\phi'\left(\frac{(u-w)^{2}+|\nabla (u-w)|^{2}}{2}\xi^{2}\right)\\
& &
\cdot\left((u-w)^{2}+|\nabla (u-w)|^{2}\right)^{2}\Bigg]dx
-(q+1)\xi^{q}\lambda \int_{\Omega}a(x)|u-w|^{q+1}dx
-(p+1)\xi^{p}\int_{\Omega}b(x)|u-w|^{p+1}dx.
\end{eqnarray*}
By (\ref{2.0.25}) and (\ref{2.0.26}), we mention that
\begin{eqnarray}\label{4.1.3}
& &
\partial_{1}F_{u}(1,0)\nonumber\\
& =&
\int_{\Omega}\Bigg[2\phi\left(\frac{u^{2}+|\nabla u|^{2}}{2}\right)\left(u^{2}+|\nabla u|^{2}\right)+\phi'\left(\frac{u^{2}+|\nabla u|^{2}}{2}\right)\left(u^{2}+|\nabla u|^{2}\right)^{2}\Bigg]dx
-(q+1)\lambda \int_{\Omega}a(x)|u|^{q+1}dx\nonumber\\
& &
-(p+1)\int_{\Omega}b(x)|u|^{p+1}dx\nonumber\\
& =&
\int_{\Omega}\Bigg[\phi\left(\frac{u^{2}+|\nabla u|^{2}}{2}\right)\left(u^{2}+|\nabla u|^{2}\right)+\phi'\left(\frac{u^{2}+|\nabla u|^{2}}{2}\right)\left(u^{2}+|\nabla u|^{2}\right)^{2}\Bigg]dx
-q\lambda \int_{\Omega}a(x)|u|^{q+1}dx\nonumber\\
& &
-p\int_{\Omega}b(x)|u|^{p+1}dx\nonumber\\
& =&\gamma_{u}''(1)>0.
\end{eqnarray}
According to the implicit function theorem, there exist $\epsilon>0$ and a differentiable function
$\xi^{+}:B(0,\epsilon)\subset H^{1}_{0}(\Omega) \rightarrow \mathbb{R}^{+}$ such that $\xi^{+}(0)=1$, $F_{u}(\xi^{+}(w),w)=0$
and
\begin{eqnarray}\label{4.1.4}
 \langle (\xi^{+})'(w),v\rangle
=-\frac{\langle \partial_{2}F_{u}(\xi^{+}(w),w),v\rangle}{\partial_{1}F_{u}(\xi^{+}(w),w)},
\;\; \forall w,v\in B(0,\epsilon),
\end{eqnarray}
where $\partial_{2}F_{u}$ and $\partial_{1}F_{u}$ denote the partial derivatives on the first and second variable, respectively.
Next, we calculate $\langle \partial_{2}F_{u}(\xi^{+}(w),w),v\rangle$.
In fact, for any $w,v\in B(0,\epsilon)$, it follows that
\begin{eqnarray}\label{4.1.5}
& &
\langle \partial_{2}F_{u}(\xi^{+}(w),w),v\rangle\nonumber\\
&= &
\lim_{t\rightarrow0}\frac{F_{u}(\xi^{+},w+tv)-F_{u}(\xi^{+},w)}{t}\nonumber\\
&= &
 \lim_{t\rightarrow0}(\xi^{+})^{2}\int_{\Omega}\frac{1}{t}
\Bigg[\phi\left(\frac{(u-(w+tv))^{2}+|\nabla(u-(w+tv))|^{2}}{2}(\xi^{+})^{2}\right)\left((u-(w+tv))^{2}+|\nabla (u-(w+tv))|^{2}\right)\nonumber\\
& &
-\phi\left(\frac{(u-w)^{2}+|\nabla (u-w)|^{2}}{2}(\xi^{+})^{2}\right)((u-w)^{2}+|\nabla (u-w)|^{2} )\Bigg]dx
\nonumber\\
& &
 -\lim_{t\rightarrow0}(\xi^{+})^{q+1}\int_{\Omega}\lambda a(x)\frac{|u-(w+tv)|^{q+1}-|u-w|^{q+1}}{t}dx\nonumber\\
& &
 -\lim_{t\rightarrow0}(\xi^{+})^{p+1}\int_{\Omega}b(x)\frac{|u-(w+tv)|^{p+1}-|u-w|^{p+1}}{t}dx.
\end{eqnarray}
For any $u\in  N^{+}_{\lambda}$ and $w,v\in B(0,\epsilon)$, let
\begin{eqnarray*}
 \tilde{f}(t)=\phi\left(\frac{(u-(w+tv))^{2}+|\nabla(u-(w+tv))|^{2}}{2}(\xi^{+})^{2}\right)\left((u-(w+tv))^{2}+|\nabla (u-(w+tv))|^{2}\right),
\;\;0<|t|\leq1.
\end{eqnarray*}
It follows from $\phi\in C^{1}([0,+\infty),\mathbb{R})$ that $\tilde{f}$ is of class $C^{1}$ and
\begin{eqnarray*}
& &\tilde{f}'(t)\\
& = &
-\Bigg[2\phi\left(\frac{(u-(w+tv))^{2}+|\nabla(u-(w+tv))|^{2}}{2}(\xi^{+})^{2}\right)
+(\xi^{+})^{2}\phi'\left(\frac{(u-(w+tv))^{2}+|\nabla(u-(w+tv))|^{2}}{2}(\xi^{+})^{2}\right)\\
& &
\cdot((u-(w+tv))^{2}+|\nabla (u-(w+tv))|^{2})\Bigg]
\left[(u-(w+tv))v+\nabla (u-(w+tv))\cdot\nabla v\right].
 \end{eqnarray*}
By the mean value theorem, there exists $\theta\in\mathbb{R}$ with $0<\theta<t\leq1$ such that
\begin{eqnarray}\label{4.1.9}
& &
\frac{1}{t}
\Bigg[\phi\left(\frac{(u-(w+tv))^{2}+|\nabla(u-(w+tv))|^{2}}{2}(\xi^{+})^{2}\right)\left((u-(w+tv))^{2}+|\nabla (u-(w+tv))|^{2}\right)\nonumber\\
& &
-\phi\left(\frac{(u-w)^{2}+|\nabla (u-w)|^{2}}{2}(\xi^{+})^{2}\right)((u-w)^{2}+|\nabla (u-w)|^{2} )\Bigg]\nonumber\\
& = & \frac{\tilde{f}(t)-\tilde{f}(0)}{t-0}\nonumber\\
& = & \tilde{f}'(\theta).
 \end{eqnarray}
Using $0<\theta<t\leq1$, $(\phi_1)$, $(\phi_2)$, H\"older inequality and Proposition 2.1, we have
\begin{eqnarray}\label{4.1.10}
|\tilde{f}'(\theta)|
\leq 2(\rho_{2}+\rho_{1})[(|u|+|v|+|w|)|v|+(|\nabla u|+|\nabla v|+|\nabla w|)|\nabla v|]
\in L^{1}(\Omega).
 \end{eqnarray}
In addition, it follows from the continuity of $\phi'$ that
\begin{eqnarray}\label{4.1.11}
& & \lim_{\theta\rightarrow0}\tilde{f}'(\theta)\nonumber\\
& = &
-(\xi^{+})^{2}\phi'\left(\frac{(u-w)^{2}+|\nabla(u-w)|^{2}}{2}(\xi^{+})^{2}\right)
((u-w)^{2}+|\nabla (u-w)|^{2})\cdot\left[(u-w)v+\nabla (u-w)\cdot\nabla v\right]\nonumber\\
& &
-2\phi\left(\frac{(u-w)^{2}+|\nabla(u-w)|^{2}}{2}(\xi^{+})^{2}\right)\left[(u-w)v+\nabla (u-w)\cdot\nabla v\right]
\;\mbox{a.e.}\; x \in \Omega.
\end{eqnarray}
Note that $\theta\rightarrow0$ as $t\rightarrow0$.
At this moment, by (\ref{4.1.9}), (\ref{4.1.10}), (\ref{4.1.11}) and the Lebesgue dominated convergence theorem, we have
\begin{eqnarray}\label{4.1.20}
& & \lim_{t\rightarrow0}(\xi^{+})^{2}\int_{\Omega}\frac{1}{t}
\Bigg[\phi\left(\frac{(u-(w+tv))^{2}+|\nabla(u-(w+tv))|^{2}}{2}(\xi^{+})^{2}\right)\left((u-(w+tv))^{2}+|\nabla (u-(w+tv))|^{2}\right)\nonumber\\
& &
-\phi\left(\frac{(u-w)^{2}+|\nabla (u-w)|^{2}}{2}(\xi^{+})^{2}\right)((u-w)^{2}+|\nabla (u-w)|^{2} )\Bigg]dx
\nonumber\\
& =&
-\int_{\Omega}(\xi^{+})^{4}\phi'\left(\frac{(u-w)^{2}+|\nabla(u-w)|^{2}}{2}(\xi^{+})^{2}\right)
((u-w)^{2}+|\nabla (u-w)|^{2})
\cdot\left[(u-w)v+\nabla (u-w)\cdot\nabla v\right]dx\nonumber\\
& &
-\int_{\Omega}2(\xi^{+})^{2}\phi\left(\frac{(u-w)^{2}+|\nabla(u-w)|^{2}}{2}(\xi^{+})^{2}\right)
\left[(u-w)v+\nabla (u-w)\cdot\nabla v\right]dx.
\end{eqnarray}
Next, we consider
\begin{eqnarray*}
\lim_{t\rightarrow0}\lambda (\xi^{+})^{q+1}\int_{\Omega}a(x)\frac{|u-(w+tv)|^{q+1}-|u-w|^{q+1}}{t}dx.
  \end{eqnarray*}
Set $f_{2}(t)=|u-(w+tv)|^{q+1}$, $0<|t|\leq1$, $u\in  N^{+}_{\lambda}$ and $w,v\in B(0,\epsilon)$.
It follows that
\begin{eqnarray*}
 f_{2}'(t)=-(q+1)|u-(w+tv)|^{q-1}(u-(w+tv))v.
  \end{eqnarray*}
By the mean value theorem, there exists $\theta\in \mathbb{R}$ with $0<\theta<t\leq1$ such that
\begin{eqnarray}\label{507}
a(x)\frac{|u-(w+tv)|^{q+1}-|u-w|^{q+1}}{t}
& = &  a(x)\frac{f_{2}(t)-f_{2}(0)}{t-0}\nonumber\\
& = &  a(x)f_{2}'(\theta).
  \end{eqnarray}
In view of $1<q+1<2<p+1<2^{\ast}$, $0<\theta<t\leq1$, $(H)$,  H\"older inequality and Proposition 2.1, we also have
\begin{eqnarray}\label{508}
| a(x)f_{2}'(\theta)|
  &\leq &  (q+1)\|a\|_{\infty}(|u|+|w|+|v|)^{q}|v|
 \in L^{1}(\Omega).
  \end{eqnarray}
Additionally, it follows that
\begin{eqnarray}\label{A.0.1}
\lim_{\theta\rightarrow0} a(x)f_{2}'(\theta)
=-(q+1)a(x)|u-w|^{q-1}(u-w)v,
\;\mbox{a.e.} \;x \in \Omega.
  \end{eqnarray}
Note that $\theta\rightarrow0$ as $t\rightarrow0$.
Thus, in view of (\ref{507}), (\ref{508}), (\ref{A.0.1}) and the Lebesgue dominated convergence theorem, we have
\begin{eqnarray}\label{4.1.7}
\lim_{t\rightarrow0}\lambda (\xi^{+})^{q+1}\int_{\Omega}a(x)\frac{|u-(w+tv)|^{q+1}-|u-w|^{q+1}}{t}dx
= -(q+1)(\xi^{+})^{q+1}\lambda \int_{\Omega}a(x)|u-w|^{q-1}(u-w)vdx.
  \end{eqnarray}
Similarly, we also have
\begin{eqnarray}\label{4.1.8}
\lim_{t\rightarrow0}(\xi^{+})^{p+1}\int_{\Omega}b(x)\frac{|u-(w+tv)|^{p+1}-|u-w|^{p+1}}{t}dx
 =-(p+1)(\xi^{+})^{p+1}\int_{\Omega}b(x)|u-w|^{p-1}(u-w)vdx.
\end{eqnarray}
Now, combing with (\ref{4.1.5}), (\ref{4.1.20}), (\ref{4.1.7}) and (\ref{4.1.8}), we have
\begin{eqnarray*}
& &
\langle \partial_{2}F_{u}(\xi^{+}(w),w),v\rangle\\
&= &
-2(\xi^{+})^{2}\int_{\Omega}\phi\left(\frac{(u-w)^{2}+|\nabla(u-w)|^{2}}{2}(\xi^{+})^{2}\right)
\left((u-w)v+\nabla (u-w)\cdot\nabla v\right)dx\\
& &
-(\xi^{+})^{4}\int_{\Omega}\phi'\left(\frac{(u-w)^{2}+|\nabla(u-w)|^{2}}{2}(\xi^{+})^{2}\right)
((u-w)^{2}+|\nabla (u-w)|^{2})\left((u-w)v+\nabla (u-w)\cdot\nabla v\right)dx\\
& &
+(q+1)(\xi^{+})^{q+1}\lambda\int_{\Omega} a(x)|u-w|^{q-1}(u-w)vdx
+(p+1)(\xi^{+})^{p+1}\int_{\Omega}b(x)|u-w|^{p-1}(u-w)vdx.
\end{eqnarray*}
Then, putting $w=0$ and $\xi^{+}=\xi^{+}(0)=1$, the last identity just above shows that
\begin{eqnarray}\label{4.1.12}
& &
\langle \partial_{2}F_{u}(1,0),v\rangle\nonumber\\
&= &
-2\int_{\Omega}\phi\left(\frac{u^{2}+|\nabla u|^{2}}{2}\right)
\left(uv+\nabla u\cdot\nabla v\right)dx
-\int_{\Omega}\phi'\left(\frac{u^{2}+|\nabla u|^{2}}{2}\right)
(u^{2}+|\nabla u|^{2})\left(uv+\nabla u\cdot\nabla v\right)dx\nonumber\\
& &
+(q+1)\lambda\int_{\Omega} a(x)|u|^{q-1}uvdx
+(p+1)\int_{\Omega}b(x)|u|^{p-1}uvdx.
\end{eqnarray}
Thus, it follows from (\ref{4.1.3}), (\ref{4.1.4}) and (\ref{4.1.12}) that (\ref{4.1.1}) holds.
\par
Next, we shall show that $\xi^{+}(w)(u-w)\in N_{\lambda}^{+}$.
It is easy to see that $\xi^{+}(w)(u-w)\in H^{1}_{0}(\Omega)\setminus\{0\}$.
Then, it follows from (\ref{4.1.2}) and $F_{u}(\xi(w),w)=0$ that
$\langle J'_{\lambda}(\xi^{+}(u-w)),\xi^{+}(u-w)\rangle=0$, that is $\xi^{+}(w)(u-w)\in N_{\lambda}$.
Then, we define $\Psi:N_{\lambda}\rightarrow \mathbb{R}$ given by
$\Psi(u)=\langle J'_{\lambda}(u),u\rangle$ for $u\in N_{\lambda}$.
By Lemma 2.2, we know that $\Psi\in C^{1}(N_{\lambda},\mathbb{R})$.
In particular, for any $u\in N_{\lambda}^{+}$, in view of (\ref{2.0.25}), (\ref{2.0.26}) and (\ref{2.0.11}), we deduce that
\begin{eqnarray*}
\langle \Psi'(u),u\rangle
& = &
\int_{\Omega}\left[\phi'\left(\frac{u^{2}+|\nabla u|^{2}}{2}\right)(u^{2}+|\nabla u|^{2})^{2}
+\phi\left(\frac{u^{2}+|\nabla u|^{2}}{2}\right)(u^{2}+|\nabla u|^{2})\right]dx\\
& &
-q\lambda \int_{\Omega}a(x)|u|^{q+1}dx-p\int_{\Omega}b(x)|u|^{p+1}dx\\
& = &
\gamma''_{u}(1)<0.
 \end{eqnarray*}
Thus, by the continuity of the functions $\Psi'$ and $\xi^{+}$, we have
\begin{eqnarray*}
&   & \gamma''_{\xi^{+}(w)(u-w)}(1)\\
& = &
\langle \Psi'(\xi^{+}(w)(u-w)),\xi^{+}(w)(u-w)\rangle\\
& = &
\int_{\Omega}\phi'\left(\frac{(\xi^{+}(w)(u-w))^{2}+|\nabla (\xi^{+}(w)(u-w))|^{2}}{2}\right)((\xi^{+}(w)(u-w))^{2}+|\nabla (\xi^{+}(w)(u-w))|^{2})^{2}dx\\
& &
+\int_{\Omega}\phi\left(\frac{(\xi^{+}(w)(u-w))^{2}+|\nabla (\xi^{+}(w)(u-w))|^{2}}{2}\right)((\xi^{+}(w)(u-w))^{2}+|\nabla (\xi^{+}(w)(u-w))|^{2})dx\\
& &
-q\lambda \int_{\Omega}a(x)|\xi^{+}(w)(u-w)|^{q+1}dx-p\int_{\Omega}b(x)|\xi^{+}(w)(u-w)|^{p+1}dx
<0
 \end{eqnarray*}
if $\epsilon$ sufficiently small, this implies that $\xi^{+}(w)(u-w)\in N_{\lambda}^{+}$.
\qed
\par
Analogously, using the same ideas discussed in the previous result, we get the following result.
\vskip2mm
 \noindent
{\bf Lemma 4.3. } {\it Suppose that $(\phi_{1})$, $(\phi_{2})$ and $(H)$ hold.
Then for each  $\lambda \in (0,\lambda_{0})$ and $u\in N_{\lambda}^{-}$, there exist $\epsilon>0$ and a differentiable function $\xi^{-}:B(0,\epsilon)\subset H^{1}_{0}(\Omega) \rightarrow \mathbb{R}^{+}$ such that $\xi^{-}(0)=1$, the function
$\xi^{-}(v)(u-v)\in N_{\lambda}^{-}$ and
\begin{eqnarray}\label{4.2.1}
 \langle (\xi^{-})'(0),v\rangle
& =&\frac{1}{\gamma_{u}''(1)}\int_{\Omega}
\Bigg[2\phi\left(\frac{u^{2}+|\nabla u|^{2}}{2}\right)(uv+\nabla u\nabla v)
 +\phi'\left(\frac{u^{2}+|\nabla u|^{2}}{2}\right)(uv+\nabla u\nabla v)(u^{2}+|\nabla u|^{2})\nonumber\\
& &
-(q+1)\lambda a(x)|u|^{q-1}uv
-(p+1)b(x)|u|^{p-1}uv\Bigg]dx
\end{eqnarray}
for all $v\in B(0,\epsilon)$.\\}
\vskip2mm
 \noindent
{\bf Proof of Proposition 4.1} Here we shall prove the item $(i)$.
The proof of item $(ii)$ follows the same lines with using Lemma 4.3 instead of Lemma 4.2 and
using Lemma 3.4.(ii) instead of Lemma 3.4.(i).
Let $\lambda \in (0,\lambda_{0})$ be fixed.
Proposition 2.2 implies that $J_{\lambda}$ is bounded from below on $N_{\lambda}$ and so on $N_{\lambda}^{+}\subset N_{\lambda}$.
Thus, we have $\inf_{N_{\lambda}^{+}}J_{\lambda}<+\infty$.
Then by Lemma 3.4.(i), we know that $J_{\lambda}(u)<0$ for any $u\in N_{\lambda}^{+}$, and so $\inf_{N_{\lambda}^{+}}J_{\lambda}<0$.
Applying the Ekeland variational principle \cite{Ekeland1974}, there exists a minimizing sequence $\{u_{n}\}\subset N_{\lambda}^{+}$ such that
\begin{eqnarray}\label{4.3.1}
J_{\lambda}(u_{n})= \inf_{u\in N_{\lambda}^{+}}J_{\lambda}(u)+o_{n}(1)<0
 \end{eqnarray}
and
\begin{eqnarray}\label{4.3.2}
J_{\lambda}(u_{n})< J_{\lambda}(v)+\frac{1}{n}\|\nabla (v-u_{n})\|_{2},\;\;\forall v\in  N_{\lambda}^{+},
 \end{eqnarray}
where $o_{n}(1)$ denotes a quantity that goes to zero as $n$ goes to $+\infty$.
Since the functional $J_{\lambda}$ is coercive on $N_{\lambda}^{+}\subset N_{\lambda}$,
$\{u_{n}\}$ is a bounded sequence in $H^{1}_{0}(\Omega)$.
\par
Next, we will show that
\begin{eqnarray}\label{4.3.3}
\lim_{n\rightarrow\infty}J'_{\lambda}(u_{n})=0.
\end{eqnarray}
Since $\{u_{n}\}\subset N_{\lambda}^{+}$, according to Lemma 4.2, there exist $\epsilon_{n}>0$ and differentiable functions $\xi_{n}^{+}:B(0,\epsilon_{n})\rightarrow \mathbb{R}^{+}$ such that $\xi_{n}^{+}(0)=1$, the functions
$\xi_{n}^{+}(w)(u_{n}-w)\in N_{\lambda}^{+}$ and
\begin{eqnarray}\label{4.3.4}
 \langle (\xi_{n}^{+})'(0),w\rangle
=\frac{\chi_{n}(w)}{\gamma_{u_{n}}''(1)}, \;\mbox{for all}\; w\in B(0,\epsilon_{n}),
\end{eqnarray}
where $\chi_{n}:H^{1}_{0}(\Omega)\rightarrow \mathbb{R}$ given by
\begin{eqnarray}\label{4.3.4.1}
\chi_{n}(v)
& =&\int_{\Omega}
\Bigg[2\phi\left(\frac{u_{n}^{2}+|\nabla u_{n}|^{2}}{2}\right)(u_{n}v+\nabla u_{n}\nabla v)
 +\phi'\left(\frac{u_{n}^{2}+|\nabla u_{n}|^{2}}{2}\right)(u_{n}v+\nabla u_{n}\nabla v)(u_{n}^{2}+|\nabla u_{n}|^{2})\nonumber\\
& &
-(q+1)\lambda a(x)|u_{n}|^{q-1}u_{n}v
-(p+1)b(x)|u_{n}|^{p-1}u_{n}v\Bigg]dx,
\;\;\forall v\in H^{1}_{0}(\Omega).
\end{eqnarray}
We put $\rho\in (0,\epsilon_{n})$.
Now we claim that there exists a constant $C>0$ which is independent on $\rho>0$, such that
\begin{eqnarray}\label{4.3.6}
\Bigg\langle J'_{\lambda}(u_{n}),\frac{u}{\|\nabla u\|_{2}}\Bigg\rangle
\leq
\frac{C}{n}(\|(\xi_{n}^{+})'(0)\|_{(H^{1}_{0}(\Omega))^{'}}+1),
\;\;
\forall u\in H^{1}_{0}(\Omega)\setminus\{0\}.
\end{eqnarray}
Furthermore, we claim also that $\|(\xi_{n}^{+})'(0)\|_{(H^{1}_{0}(\Omega))^{'}}$ is uniformly bounded in $n$, that is, there exists a constant $\tilde{C}>0$ which is independent on $n$, such that
\begin{eqnarray}\label{4.3.14}
\|(\xi_{n}^{+})'(0)\|_{(H^{1}_{0}(\Omega))^{'}}\leq \tilde{C} \;\mbox{for each}\; n\in \mathbb{N}.
\end{eqnarray}
\par
In what follows we shall prove the claims given just above.
Define the auxiliary function
\begin{eqnarray}\label{4.3.5}
w_{\rho}=\frac{\rho u}{\|\nabla u\|_{2}}\in B(0,\epsilon_{n}).
\end{eqnarray}
Using one more time Lemma 4.2, we mention that
\begin{eqnarray}\label{4.3.7}
\mu_{\rho}=\xi_{n}^{+}(w_{\rho})(u_{n}-w_{\rho})\in N_{\lambda}^{+}.
\end{eqnarray}
Then by (\ref{4.3.2}), we get
\begin{eqnarray}\label{4.3.8}
-\frac{1}{n}\|\nabla (\mu_{\rho}-u_{n})\|_{2}
\leq
J_{\lambda}(\mu_{\rho})-J_{\lambda}(u_{n}).
 \end{eqnarray}
Notice also that
\begin{eqnarray}
\label{4.3.9} & & w_{\rho}\rightarrow0\;\mbox{and}\;\mu_{\rho}\rightarrow u_{n}\;\mbox{in}\; H^{1}_{0}(\Omega),\\
\label{4.3.60}& &\xi_{n}^{+}(w_{\rho})\rightarrow1 \;\mbox{in}\; \mathbb{R},\\
\label{4.3.61}& & J'_{\lambda}(\mu_{\rho})\rightarrow J'_{\lambda}(u_{n})\;\mbox{in}\; (H^{1}_{0}(\Omega))'
 \end{eqnarray}
as $\rho\rightarrow0$ hold true for any $n\in\mathbb{N}$.
Then applying mean value theorem, there exists $t\in (0,1)$ such that
\begin{eqnarray*}
& &J_{\lambda}(\mu_{\rho})-J_{\lambda}(u_{n})\nonumber\\
& =&
\langle J'_{\lambda}((1-t)\mu_{\rho}+tu_{n}),\mu_{\rho}-u_{n}\rangle\nonumber\\
& =&
\langle J'_{\lambda}(\mu_{\rho}+t(u_{n}-\mu_{\rho}))-J'_{\lambda}(u_{n}),\mu_{\rho}-u_{n}\rangle
+\langle J'_{\lambda}(u_{n}),\mu_{\rho}-u_{n}\rangle\nonumber\\
& \leq&
|\langle J'_{\lambda}(\mu_{\rho}+t(u_{n}-\mu_{\rho}))-J'_{\lambda}(u_{n}),\mu_{\rho}-u_{n}\rangle|
+\langle J'_{\lambda}(u_{n}),\mu_{\rho}-u_{n}\rangle\nonumber\\
& \leq&
\|\langle J'_{\lambda}(\mu_{\rho}+t(u_{n}-\mu_{\rho}))-J'_{\lambda}(u_{n})\rangle\|_{(H^{1}_{0}(\Omega))'}
\|\nabla (\mu_{\rho}-u_{n})\|_{2}
+\langle J'_{\lambda}(u_{n}),\mu_{\rho}-u_{n}\rangle\nonumber\\
& =& o_{\rho}(1)\|\nabla (\mu_{\rho}-u_{n})\|_{2}+\langle J'_{\lambda}(u_{n}),\mu_{\rho}-u_{n}\rangle,
\end{eqnarray*}
where $o_{\rho}(1)$ denotes a quantity that goes to zero as $\rho$ goes to zero.
Therefore, combing with the estimate just above and (\ref{4.3.8}), we conclude that
\begin{eqnarray*}
-\frac{1}{n}\|\nabla (\mu_{\rho}-u_{n})\|_{2}
\leq
o_{\rho}(1)\|\nabla (\mu_{\rho}-u_{n})\|_{2}+\langle J'_{\lambda}(u_{n}),\mu_{\rho}-u_{n}\rangle.
 \end{eqnarray*}
Then combing with (\ref{4.3.7}), we have
\begin{eqnarray*}
& &-\frac{1}{n}\|\nabla (\mu_{\rho}-u_{n})\|_{2}\\
& \leq&
\langle J'_{\lambda}(u_{n}), \xi_{n}^{+}(w_{\rho})(u_{n}-w_{\rho})-u_{n}\rangle +o_{\rho}(1)\|\nabla (\mu_{\rho}-u_{n})\|_{2}\\
& =&
\langle J'_{\lambda}(u_{n}), \xi_{n}^{+}(w_{\rho})u_{n}-\xi_{n}^{+}(w_{\rho})w_{\rho}-u_{n}\rangle +o_{\rho}(1)\|\nabla (\mu_{\rho}-u_{n})\|_{2}\\
& =&
\langle J'_{\lambda}(u_{n}), (\xi_{n}^{+}(w_{\rho})-1)u_{n}-\xi_{n}^{+}(w_{\rho})w_{\rho}\rangle
+ \langle J'_{\lambda}(u_{n}),w_{\rho}\rangle
-\langle J'_{\lambda}(u_{n}),w_{\rho}\rangle
+o_{\rho}(1)\|\nabla (\mu_{\rho}-u_{n})\|_{2}\\
& =&
\langle J'_{\lambda}(u_{n}), (\xi_{n}^{+}(w_{\rho})-1)u_{n}-\xi_{n}^{+}(w_{\rho})w_{\rho}+w_{\rho}\rangle
-\langle J'_{\lambda}(u_{n}),w_{\rho}\rangle
+o_{\rho}(1)\|\nabla (\mu_{\rho}-u_{n})\|_{2}\\
& =&
\langle J'_{\lambda}(u_{n}), (\xi_{n}^{+}(w_{\rho})-1)u_{n}-(\xi_{n}^{+}(w_{\rho})-1)w_{\rho}\rangle
-\langle J'_{\lambda}(u_{n}),w_{\rho}\rangle
+o_{\rho}(1)\|\nabla (\mu_{\rho}-u_{n})\|_{2}\\
& =&
(\xi_{n}^{+}(w_{\rho})-1)\langle J'_{\lambda}(u_{n}),u_{n}-w_{\rho}\rangle
-\langle J'_{\lambda}(u_{n}),w_{\rho}\rangle
+o_{\rho}(1)\|\nabla (\mu_{\rho}-u_{n})\|_{2}\\
& =&
(\xi_{n}^{+}(w_{\rho})-1)\langle J'_{\lambda}(u_{n})-J'_{\lambda}(\mu_{\rho}),u_{n}-w_{\rho}\rangle
+(\xi_{n}^{+}(w_{\rho})-1)\langle J'_{\lambda}(\mu_{\rho}),u_{n}-w_{\rho}\rangle
-\langle J'_{\lambda}(u_{n}),w_{\rho}\rangle\\
& &+o_{\rho}(1)\|\nabla (\mu_{\rho}-u_{n})\|_{2}\\
& =&
(\xi_{n}^{+}(w_{\rho})-1)\langle J'_{\lambda}(u_{n})-J'_{\lambda}(\mu_{\rho}),u_{n}-w_{\rho}\rangle
+\frac{1}{\xi_{n}^{+}}(\xi_{n}^{+}(w_{\rho})-1)\langle J'_{\lambda}(\mu_{\rho}),\xi_{n}^{+}(u_{n}-w_{\rho})\rangle
-\langle J'_{\lambda}(u_{n}),w_{\rho}\rangle\\
& &+o_{\rho}(1)\|\nabla (\mu_{\rho}-u_{n})\|_{2}\\
& =&
(\xi_{n}^{+}(w_{\rho})-1)\langle J'_{\lambda}(u_{n})-J'_{\lambda}(\mu_{\rho}),u_{n}-w_{\rho}\rangle
+\frac{1}{\xi_{n}^{+}}(\xi_{n}^{+}(w_{\rho})-1)\langle J'_{\lambda}(\mu_{\rho}),\mu_{\rho}\rangle
-\langle J'_{\lambda}(u_{n}),w_{\rho}\rangle\\
& &+o_{\rho}(1)\|\nabla (\mu_{\rho}-u_{n})\|_{2}\\
& =&
(\xi_{n}^{+}(w_{\rho})-1)\langle J'_{\lambda}(u_{n})-J'_{\lambda}(\mu_{\rho}),u_{n}-w_{\rho}\rangle
-\langle J'_{\lambda}(u_{n}),w_{\rho}\rangle
+o_{\rho}(1)\|\nabla (\mu_{\rho}-u_{n})\|_{2}.
 \end{eqnarray*}
Then, by (\ref{4.3.5}), we have
\begin{eqnarray*}
-\frac{1}{n}\|\nabla (\mu_{\rho}-u_{n})\|_{2}
\leq
(\xi_{n}^{+}(w_{\rho})-1)\left\langle J'_{\lambda}(u_{n})-J'_{\lambda}(\mu_{\rho}),u_{n}-\frac{\rho u}{\|\nabla u\|_{2}}\right\rangle
-\left\langle J'_{\lambda}(u_{n}),\frac{\rho u}{\|\nabla u\|_{2}}\right\rangle
+o_{\rho}(1)\|\nabla (\mu_{\rho}-u_{n})\|_{2},
 \end{eqnarray*}
which implies that
\begin{eqnarray}\label{4.3.51}
& &\left\langle J'_{\lambda}(u_{n}),\frac{ u}{\|\nabla u\|_{2}}\right\rangle\nonumber\\
& \leq&
\frac{\|\nabla (\mu_{\rho}-u_{n})\|_{2}}{n\rho}
+\frac{\xi_{n}^{+}(w_{\rho})-1}{\rho}\left\langle J'_{\lambda}(u_{n})-J'_{\lambda}(\mu_{\rho}),u_{n}-\frac{\rho u}{\|\nabla u\|_{2}}\right\rangle
+o_{\rho}(1)\frac{\|\nabla (\mu_{\rho}-u_{n})\|_{2}}{\rho}\nonumber\\
& \leq&
\frac{\|\nabla (\mu_{\rho}-u_{n})\|_{2}}{n\rho}
+\frac{|\xi_{n}^{+}(w_{\rho})-1|}{\rho}\left|\left\langle J'_{\lambda}(u_{n})-J'_{\lambda}(\mu_{\rho}),u_{n}-\frac{\rho u}{\|\nabla u\|_{2}}\right\rangle\right|
+o_{\rho}(1)\frac{\|\nabla (\mu_{\rho}-u_{n})\|_{2}}{\rho}\nonumber\\
& \leq&
\frac{\|\nabla (\mu_{\rho}-u_{n})\|_{2}}{n\rho}
+\frac{|\xi_{n}^{+}(w_{\rho})-1|}{\rho}\|J'_{\lambda}(u_{n})-J'_{\lambda}(\mu_{\rho})\|_{(H^{1}_{0}(\Omega))'}
\left\|\nabla (u_{n}-\frac{\rho u}{\|\nabla u\|_{2}})\right\|_{2}
+o_{\rho}(1)\frac{\|\nabla (\mu_{\rho}-u_{n})\|_{2}}{\rho}\nonumber\\
& \leq&
\frac{\|\nabla (\mu_{\rho}-u_{n})\|_{2}}{n\rho}
+\frac{|\xi_{n}^{+}(w_{\rho})-1|}{\rho}\|J'_{\lambda}(u_{n})-J'_{\lambda}(\mu_{\rho})\|_{(H^{1}_{0}(\Omega))'}
(\|\nabla u_{n}\|_{2}+\rho)
+o_{\rho}(1)\frac{\|\nabla (\mu_{\rho}-u_{n})\|_{2}}{\rho}.
 \end{eqnarray}
It is not hard to see that
\begin{eqnarray}\label{4.3.50}
\frac{\|\nabla (\mu_{\rho}-u_{n})\|_{2}}{\rho}
&= &
\frac{1}{\rho}\left\|\nabla \left(\xi_{n}^{+}(w_{\rho})(u_{n}-w_{\rho})-u_{n} \right)\right\|_{2}\nonumber\\
&= &
\frac{1}{\rho}\left\|\nabla \left((\xi_{n}^{+}(w_{\rho})-1)u_{n}-\xi_{n}^{+}(w_{\rho})w_{\rho}\right)\right\|_{2}\nonumber\\
&= &
\frac{1}{\rho}\left\|\nabla \left((\xi_{n}^{+}(w_{\rho})-1)u_{n}-\xi_{n}^{+}(w_{\rho})\frac{\rho u}{\|\nabla u\|_{2}}\right)\right\|_{2}\nonumber\\
&\leq &
\frac{|\xi_{n}^{+}(w_{\rho})-1|}{\rho}\|\nabla u_{n}\|_{2}+|\xi_{n}^{+}(w_{\rho})|,
\end{eqnarray}
and
\begin{eqnarray}\label{4.3.12}
\lim_{\rho\rightarrow0}\frac{|\xi_{n}^{+}(w_{\rho})-1|}{\rho}
=
\Bigg\langle (\xi_{n}^{+})'(0),\frac{u_{n}}{\|\nabla u_{n}\|_{2}}\Bigg\rangle
\leq
\|(\xi_{n}^{+})'(0)\|_{(H^{1}_{0}(\Omega))^{'}}.
\end{eqnarray}
Let $\rho\rightarrow0$ in (\ref{4.3.50}), using the fact that (\ref{4.3.60}), (\ref{4.3.12}) and $\{u_{n}\}$ is bounded in $H^{1}_{0}(\Omega)$, we infer that
\begin{eqnarray}\label{4.3.13}
\lim_{\rho\rightarrow0}\frac{\|\nabla (\mu_{\rho}-u_{n})\|_{2}}{\rho}
&\leq &
\lim_{\rho\rightarrow0} \Bigg[ |\xi_{n}^{+}(w_{\rho})|+\frac{|\xi_{n}^{+}(w_{\rho})-1|}{\rho}\|\nabla u_{n}\|_{2}\Bigg]\nonumber\\
&\leq &
1+\|\nabla u_{n}\|_{2}\|(\xi_{n}^{+})'(0)\|_{(H^{1}_{0}(\Omega))^{'}}\nonumber\\
&\leq &
C\left(1+\|(\xi_{n}^{+})'(0)\|_{(H^{1}_{0}(\Omega))^{'}}\right),
\end{eqnarray}
where $C$ is a positive constant which is independent on $\rho>0$.
Furthermore, by (\ref{4.3.12}) and (\ref{4.3.13}), we have
\begin{eqnarray}\label{4.3.53}
\frac{|\xi_{n}^{+}(w_{\rho})-1|}{\rho}
\leq
\|(\xi_{n}^{+})'(0)\|_{(H^{1}_{0}(\Omega))^{'}}
\;\;\mbox{and}\;\;
\frac{\|\nabla (\mu_{\rho}-u_{n})\|_{2}}{\rho}
\leq
C\left(1+\|(\xi_{n}^{+})'(0)\|_{(H^{1}_{0}(\Omega))^{'}}\right)
\end{eqnarray}
for $\rho>0$ small enough.
Now, we let $\rho\rightarrow0$ in (\ref{4.3.51}), combing with (\ref{4.3.13}), (\ref{4.3.53}), (\ref{4.3.61}) and $\|\nabla u_{n}\|_{2}$ is bounded for $n\in\mathbb{N}$, we obtain the inequality (\ref{4.3.6}).
\par
It remains to show that there exists a constant $\tilde{C}>0$ which is independent on $n$, such that
\begin{eqnarray*}
\|(\xi_{n}^{+})'(0)\|_{(H^{1}_{0}(\Omega))^{'}}\leq \tilde{C} \;\mbox{for each}\; n\in \mathbb{N}.
\end{eqnarray*}
In view of (\ref{4.3.4}), we have
\begin{eqnarray*}
\|(\xi_{n}^{+})'(0)\|_{(H^{1}_{0}(\Omega))^{'}}
=
\sup_{v\in H^{1}_{0}(\Omega)\backslash \{0\} }
\frac{|\chi_{n}(v)|}{\|\nabla v\|_{2}}\cdot\frac{1}{\gamma_{u_{n}}''(1)}.
\end{eqnarray*}
By (\ref{4.3.4.1}), $(\phi_{1})$, $(\phi_{2})$, $(H)$, H\"older inequality, Proposition 2.1 and $\{u_{n}\}$ is bounded in $H^{1}_{0}(\Omega)$, we observe that
\begin{eqnarray*}
& &|\chi_{n}(v)|\\
& \leq&
\int_{\Omega}
\Bigg[2\phi\left(\frac{u_{n}^{2}+|\nabla u_{n}|^{2}}{2}\right)
 +2\left|\phi'\left(\frac{u_{n}^{2}+|\nabla u_{n}|^{2}}{2}\right)\right|\frac{u_{n}^{2}+|\nabla u_{n}|^{2}}{2}\Bigg]|u_{n}v+\nabla u_{n}\nabla v|dx\\
& &
+(q+1)\lambda\int_{\Omega} |a(x)||u_{n}|^{q}|v|dx
+(p+1)\int_{\Omega}|b(x)||u_{n}|^{p}|v|dx\\
& \leq&
2(\rho_{1}+\rho_{2})
\Bigg[\int_{\Omega}|u_{n}||v|dx+\int_{\Omega}|\nabla u_{n}||\nabla v|dx\Bigg]
+(q+1)\lambda\|a\|_{\infty}\int_{\Omega} |u_{n}|^{q}|v|dx\\
& &
+(p+1)\|b\|_{\infty}\int_{\Omega}|u_{n}|^{p}|v|dx\\
& \leq&
\Bigg(2(\rho_{1}+\rho_{2})(1+S_{2}^{2})\|\nabla u_{n}\|_{2}
+(q+1)\lambda\|a\|_{\infty}S_{q+1}^{q+1}\|\nabla u_{n}\|_{2}^{q}
+(p+1)\|b\|_{\infty}S_{p+1}^{p+1}\|\nabla u_{n}\|_{2}^{p}\Bigg)\|\nabla v\|_{2}\\
& \leq&
c\|\nabla v\|_{2},
\end{eqnarray*}
where $c$ is a positive constant.
Therefore, we only need to show that
\begin{eqnarray}\label{4.3.17}
\gamma_{u_{n}}''(1)>d
\end{eqnarray}
for some $d>0$ and $n$ large enough.
We argue by contradiction.
Assume that there exists a subsequence $u_{n}$ such that
\begin{eqnarray}\label{4.3.18}
\gamma_{u_{n}}''(1)=o_{n}(1).
\end{eqnarray}
Firstly, we shall give the following important estimates
\begin{eqnarray}\label{4.3.19}
\liminf_{n\rightarrow\infty}\|\nabla u_{n}\|_{2}
\geq
\Bigg[-\frac{(p+1)(q+1)}{\lambda(p-q)\|a\|_{\infty}S^{q+1}_{q+1}}\inf_{u\in N_{\lambda}^{+}}J_{\lambda}(u)            \Bigg]^{\frac{1}{q+1}}
>0.
\end{eqnarray}
By (\ref{2.0.3}), $(\phi_1)$, $(H)$, $1<q+1<2<p+1<2^{\ast}$ and Proposition 2.1, we mention that
\begin{eqnarray*}
J_{\lambda}(u_{n})
& \geq &
\left(\frac{\rho_{0}}{2}-\frac{\rho_{1}}{p+1}\right)(\|u_{n}\|_{2}^{2}+\|\nabla u_{n}\|_{2}^{2})
-\lambda\left(\frac{1}{q+1}-\frac{1}{p+1}\right)\|a\|_{\infty}\|u_{n}\|_{q+1}^{q+1}\\
& \geq &
-\frac{\lambda(p-q)}{(q+1)(p+1)}\|a\|_{\infty}S^{q+1}_{q+1}\|\nabla u_{n}\|_{2}^{q+1}.
\end{eqnarray*}
Then combing with (\ref{4.3.1}), we get
\begin{eqnarray*}
-\frac{\lambda(p-q)}{(q+1)(p+1)}\|a\|_{\infty}S^{q+1}_{q+1}\|\nabla u_{n}\|_{2}^{q+1}
\leq
\inf_{u\in N_{\lambda}^{+}}J_{\lambda}(u)   +\frac{1}{n}
\end{eqnarray*}
holds for any $n\in \mathbb{N}$ large enough, that is
\begin{eqnarray*}
\|\nabla u_{n}\|_{2}
>
\Bigg[-\left(\inf_{u\in N_{\lambda}^{+}}J_{\lambda}(u)
+\frac{1}{n}\right)\frac{(p+1)(q+1)}{\lambda(p-q)\|a\|_{\infty}S^{q+1}_{q+1}}\Bigg]^{\frac{1}{q+1}}
\end{eqnarray*}
for any $n\in \mathbb{N}$ large enough.
Taking the limits in estimate just above, we have (\ref{4.3.19}) holds since $\inf_{u\in N_{\lambda}^{+}}J_{\lambda}(u)<0$.
Then, it follows from (\ref{2.6}) and (\ref{4.3.18})  that
\begin{eqnarray*}
 & &\int_{\Omega}\left[\phi'\left( \frac{u_{n}^{2}+|\nabla u_{n}|^{2}}{2}\right)(u_{n}^{2}+|\nabla u_{n}|^{2})^{2}
 +(1-q)\phi\left( \frac{u_{n}^{2}+|\nabla u_{n}|^{2}}{2}\right)(u_{n}^{2}+|\nabla u_{n}|^{2})\right]dx\\
 & =&
 (p-q)\int_{\Omega}b(x)|u_{n}|^{p+1}dx
 +o_{n}(1).
\end{eqnarray*}
By $(\phi_3)$, we have
\begin{eqnarray}\label{4.3.20}
\int_{\Omega}\left[\phi'\left( \frac{u_{n}^{2}+|\nabla u_{n}|^{2}}{2}\right)(u_{n}^{2}+|\nabla u_{n}|^{2})^{2}
 +(1-q)\phi\left( \frac{u_{n}^{2}+|\nabla u_{n}|^{2}}{2}\right)(u_{n}^{2}+|\nabla u_{n}|^{2})\right]dx
\geq
\rho_{3}\|\nabla u_{n}\|_{2}^{2}.
\end{eqnarray}
It follows from $q+1<p+1$, $(H)$ and Proposition 2.1 that
\begin{eqnarray*}
(p-q)\int_{\Omega}b(x)|u_{n}|^{p+1}dx
\leq
(p-q)\|b\|_{\infty}S^{p+1}_{p+1}\|\nabla u_{n}\|_{2}^{p+1}.
\end{eqnarray*}
Hence, using the estimates above, we get
\begin{eqnarray*}
\rho_{3}\|\nabla u_{n}\|_{2}^{2}
\leq
(p-q)\|b\|_{\infty}S^{p+1}_{p+1}\|\nabla u_{n}\|_{2}^{p+1}
+o_{n}(1),
\end{eqnarray*}
which implies that
\begin{eqnarray*}
\rho_{3}
\leq
(p-q)\|b\|_{\infty}S^{p+1}_{p+1}\|\nabla u_{n}\|_{2}^{p-1}
+\frac{o_{n}(1)}{\|\nabla u_{n}\|_{2}^{2}},
\end{eqnarray*}
Then using (\ref{4.3.19}), we have
\begin{eqnarray}\label{4.3.21}
\|\nabla u_{n}\|_{2}
\geq
\left[  \frac{\rho_{3}}{(p-q)\|b\|_{\infty}S_{p+1}^{p+1}}\right]^{\frac{1}{p-1}}+o_{n}(1).
 \end{eqnarray}
Moreover, by (\ref{2.7}) and (\ref{4.3.18}), we deduce that
\begin{eqnarray*}
& & \int_{\Omega}\left[ (p-1)\phi\left( \frac{u_{n}^{2}+|\nabla u_{n}|^{2}}{2}\right)(u_{n}^{2}+|\nabla u_{n}|^{2})
 -\phi'\left( \frac{u_{n}^{2}+|\nabla u_{n}|^{2}}{2}\right)(u_{n}^{2}+|\nabla u_{n}|^{2})^{2}\right]dx\\
& =&
\lambda(p-q)\int_{\Omega}a(x)|u_{n}|^{q+1}dx+o_{n}(1).
\end{eqnarray*}
In view of $(\phi_4)$, we have
\begin{eqnarray}\label{4.3.22}
\int_{\Omega}\left[ (p-1)\phi\left( \frac{u_{n}^{2}+|\nabla u_{n}|^{2}}{2}\right)(u_{n}^{2}+|\nabla u_{n}|^{2})
 -\phi'\left( \frac{u_{n}^{2}+|\nabla u_{n}|^{2}}{2}\right)(u_{n}^{2}+|\nabla u_{n}|^{2})^{2}\right]dx
\geq
\rho_{5}\|\nabla u_{n}\|_{2}^{2}.
\end{eqnarray}
Taking into account that $q+1<p+1$, $(H)$ and Proposition 2.1, we observe that
\begin{eqnarray*}
\lambda(p-q)\int_{\Omega}a(x)|u_{n}|^{q+1}dx
 \leq
\lambda(p-q)\|a\|_{\infty}S^{q+1}_{q+1}\|\nabla u_{n}\|_{2}^{q+1}.
\end{eqnarray*}
To sum up, using the estimate (\ref{4.3.19}), we have
\begin{eqnarray}\label{4.3.23}
\|\nabla u_{n}\|_{2}
& \leq &
\left[  \lambda\frac{(p-q)\|a\|_{\infty}S_{q+1}^{q+1}}{\rho_{5}}\right]^{\frac{1}{1-q}}
+\frac{o_{n}(1)}{\rho_{5}^{\frac{1}{1-q}}\|\nabla u_{n}\|_{2}^{\frac{q+1}{1-q}}}\nonumber\\
& = &
\left[  \lambda\frac{(p-q)\|a\|_{\infty}S_{q+1}^{q+1}}{\rho_{5}}\right]^{\frac{1}{1-q}}
+o_{n}(1).
 \end{eqnarray}
Combining with (\ref{4.3.21}) and (\ref{4.3.23}), we have a contradiction for each
$\lambda <\lambda_{0}=\min\{\lambda_{1},\lambda_{2}\}$ where $\lambda_{1}$ and $\lambda_{2}$ were given in Lemma 2.1 and Lemma 3.3, respectively.
This finishes the proof.
\qed

\section{The proof of main results}
\vskip2mm
 \noindent
{\bf Proof of Theorem 1.1.}
We shall split the proof into two cases.
In the first one, we shall prove that the problem (\ref{eq1}) admits at least one non-trivial weak solution $u_{0}\in N_{\lambda}^{-}$ satisfying $J_{\lambda}(u_{0})>0$ for any $0<\lambda<\lambda_{0}$.
 Let $\lambda \in(0,\lambda_{0})$ be fixed.
Taking into account Proposition 4.1, there exists a minimizing sequence $\{u_{n}\}\subset N_{\lambda}^{-}$ such that
\begin{eqnarray}
\label{1.1.1}
& & \lim_{n\rightarrow\infty}J_{\lambda}(u_{n})= \inf_{u\in N_{\lambda}^{-}}J_{\lambda}(u)>0,\\
\label{1.1.2}
& &  \lim_{n\rightarrow\infty}J'_{\lambda}(u_{n})=0.
\end{eqnarray}
Note that Proposition 2.2 implies that the functional $J_{\lambda}$ is coercive on $N_{\lambda}^{-}\subset N_{\lambda}$.
Thus, $\{u_{n}\}$ is a bounded sequence in $H^{1}_{0}(\Omega)$, that is, there exists $M>0$ such that
\begin{eqnarray}\label{1.1.4}
\|\nabla u_{n}\|_{2}\leq M \;\mbox{for all} \;n\in \mathbb{N}.
 \end{eqnarray}
Without any loss of generality, we assume that $u_{n}\rightharpoonup u_{0} $ in $H^{1}_{0}(\Omega)$.
It is easy to see that
\begin{eqnarray}\label{1.1.40}
\liminf_{n\rightarrow\infty}\|\nabla u_{n}\|_{2}^{2}\geq \|\nabla u_{0}\|_{2}^{2}.
 \end{eqnarray}
By Proposition 2.1 and \cite[Theorem 4.9, p.94]{Brezis2010}, we have
\begin{eqnarray}
\label{582} & & u_{n}\rightarrow u_{0} \;\mbox{in} \;L^{\alpha}(\Omega)  \;\mbox{for}\;1<\alpha<2^{\ast},\\
\label{1.1.14}& &  u_{n}\rightarrow u_{0} \;\mbox{a.e.} \;\mbox{in}\; \Omega.
\end{eqnarray}
Let
\begin{eqnarray*}
\Gamma_{\lambda}(u)=\int_{\Omega}\Phi\left(\frac{u^{2}+|\nabla u|^{2}}{2}\right)dx,\;\;\forall u\in H_{0}^{1}(\Omega).
\end{eqnarray*}
It follows from (\ref{2.1.1}) that
\begin{eqnarray*}
\langle J'_{\lambda}(u_{n}),u_{n}-u_{0}\rangle
=
\langle \Gamma'_{\lambda}(u_{n}),u_{n}-u_{0}\rangle
-\lambda\int_{\Omega}a(x)|u_{n}|^{q-1}u_{n}(u_{n}-u_{0})dx
-\int_{\Omega}b(x)|u_{n}|^{p-1}u_{n}(u_{n}-u_{0})dx.
\end{eqnarray*}
Clearly, (\ref{1.1.2}) implies that $\langle J'_{\lambda}(u_{n}),u_{n}-u_{0}\rangle\rightarrow 0$.
By $1<q+1<2<p+1<2^{\ast}$, $(H)$, H\"older inequality and (\ref{582}), we have
\begin{eqnarray*}
\left|\int_{\Omega}a(x)|u_{n}|^{q-1}u_{n}(u_{n}-u_{0})dx\right|
& \leq &
\|a\|_{\infty}\left(\int_{\Omega}|u_{n}|^{q\cdot\frac{q+1}{q}}dx\right)^{\frac{q}{q+1}}
\left(\int_{\Omega}|u_{n}-u_{0}|^{q+1}dx\right)^{\frac{1}{q+1}}\\
& = &
\|a\|_{\infty}\|u_{n}\|_{q+1}^{q}\|u_{n}-u_{0}\|_{q+1}\rightarrow0
\end{eqnarray*}
and
\begin{eqnarray*}
\left|\int_{\Omega}b(x)|u_{n}|^{p-1}u_{n}(u_{n}-u_{0})dx\right|
& \leq &
\|b\|_{\infty}\left(\int_{\Omega}|u_{n}|^{p\cdot\frac{p+1}{p}}dx\right)^{\frac{p}{p+1}}
\left(\int_{\Omega}|u_{n}-u_{0}|^{p+1}dx\right)^{\frac{1}{p+1}}\\
& = &
\|b\|_{\infty}\|u_{n}\|_{p+1}^{p}\|u_{n}-u_{0}\|_{p+1}\rightarrow0.
\end{eqnarray*}
Hence, we have $\langle \Gamma'_{\lambda}(u_{n}),u_{n}-u_{0}\rangle\rightarrow 0$.
Then combing with $(\phi_1)$, $(\phi_6)$ and $(\phi_7)$, and applying \cite[Lemma 3.5, p.11]{Jeanjean2022}, we have
\begin{eqnarray}\label{1.1.13}
\nabla u_{n}(x)\rightarrow \nabla u_{0}(x) \;\mbox{a.e.}\; \;\mbox{in}\;\Omega.
 \end{eqnarray}
In addition, by (\ref{582}), $1<q+1<2<p+1<2^{\ast}$ and $(H)$, we also have
\begin{eqnarray}
\label{584} & &\lim_{n\rightarrow\infty}\int_{\Omega}a(x)|u_{n}|^{q+1}dx = \int_{\Omega}a(x)|u_{0}|^{q+1}dx,\\
\label{1.1.6} & &\lim_{n\rightarrow\infty}\int_{\Omega}b(x)|u_{n}|^{p+1}dx =  \int_{\Omega}b(x)|u_{0}|^{p+1}dx.
\end{eqnarray}
\par
Next, we shall prove that $\int_{\Omega}b(x)|u_{0}|^{p+1}dx>0$.
Using (\ref{2.0.2}), $1<q+1<2<p+1$ and $(\phi_1)$, we easily see that
\begin{eqnarray*}
& &
\int_{\Omega}b(x)|u_{n}|^{p+1}dx\\
& =&
\frac{1}{\frac{1}{q+1}-\frac{1}{p+1}}
\left\{J_{\lambda}(u_{n})+\int_\Omega \left[\frac{1}{q+1}\phi\left(\frac{u_{n}^{2}+|\nabla u_{n}|^{2}}{2}\right)(u_{n}^{2}+|\nabla u_{n}|^{2})
-\Phi\left(\frac{u_{n}^{2}+|\nabla u_{n}|^{2}}{2}\right)\right]dx\right\}\\
& \geq&
\frac{1}{\frac{1}{q+1}- \frac{1}{p+1}}\left\{J_{\lambda}(u_{n})
+\left(\frac{\rho_{0}}{q+1}-\frac{\rho_{1}}{2}\right)\int_\Omega(u_{n}^{2}+|\nabla u_{n}|^{2})dx\right\}\\
& \geq&
\frac{1}{\frac{1}{q+1}- \frac{1}{p+1}}\left\{J_{\lambda}(u_{n})
+\left(\frac{\rho_{0}}{q+1}-\frac{\rho_{1}}{2}\right)\|\nabla u_{n}\|_{2}^{2}\right\}.
\end{eqnarray*}
Taking the limits in estimate above and combing with (\ref{1.1.6}), (\ref{1.1.1}) and (\ref{1.1.40}), we have $\int_{\Omega}b(x)|u_{0}|^{p+1}dx>0$.
Hence, $u_{0}\not\equiv 0$.
Moreover, by $\clubsuit^{3}$ and $\clubsuit^{4}$ in Section 3, we know that there exists a critical point $t_{0}>0$ of the fibering map $\gamma_{u_{0}}$ such that $\gamma'_{u_{0}}(t_{0})=0$, $t_{0}u_{0}\in N_{\lambda}^{-}$ and $\gamma_{u_{0}}(t_{0})>0$.
\par
At this stage, we shall prove that
\begin{eqnarray}\label{1.1.53}
u_{n}\rightarrow u_{0} \;\mbox{in}\; H^{1}_{0}(\Omega).
 \end{eqnarray}
The proof follows arguing by contradiction.
Assume that $u_{n}\not\rightarrow u_{0}$ in $H^{1}_{0}(\Omega)$.
Then there exist a subsequence of $\{u_{n}\}$, still denoted by $\{u_{n}\}$, and a positive constant $\delta>0$, such that
\begin{eqnarray}\label{1.1.3}
\lim_{n\rightarrow\infty}\int_{\Omega}|\nabla u_{n}-\nabla u_{0}|^{2}dx\geq \delta>0.
 \end{eqnarray}
Define a function $j_{1}:\mathbb{R}\rightarrow\mathbb{R}$ by
\begin{eqnarray*}
j_{1}(s)=\Phi(|s|), \;\;\forall s\in \mathbb{R}.
\end{eqnarray*}
It follows that $j_{1}(0)=0$.
Let $f_{n}=\frac{|u_{n}|^{2}+|\nabla u_{n}|^{2}}{2}t_{0}^{2}$ and $f=\frac{u_{0}^{2}+|\nabla u_{0}|^{2}}{2}t_{0}^{2}$.
By $(\phi_1)$, (\ref{1.1.4}) and Proposition 2.1, we have
\begin{eqnarray*}
\int_{\Omega}|j_{1}(f_{n}(x))|dx
& =&
\int_{\Omega}\left|\Phi\left(\frac{u_{n}^{2}+|\nabla u_{n}|^{2}}{2}t_{0}^{2}\right)\right|dx\nonumber\\
& \leq&
\frac{\rho_{1}}{2}t_{0}^{2}(1+S_{2}^{2})\|\nabla u_{n}\|_{2}^{2}\nonumber\\
& \leq&
\frac{\rho_{1}}{2}t_{0}^{2}(1+S_{2}^{2})M^{2}<+\infty.
\end{eqnarray*}
It follows from (\ref{1.1.14}) and (\ref{1.1.13}) that $f_{n}(x)\rightarrow f(x)$ a.e. in $\Omega$.
Therefore, applying the Brezis-Lieb Lemma \cite{Brezis1983}, we have
\begin{eqnarray}\label{1.1.7}
\lim_{n\rightarrow\infty}\int_{\Omega}
\Bigg[\Phi\left(\frac{u_{n}^{2}+|\nabla u_{n}|^{2}}{2}t_{0}^{2}\right)
-\Phi\left(\left|\frac{u_{n}^{2}+|\nabla u_{n}|^{2}}{2}
-\frac{u_{0}^{2}+|\nabla u_{0}|^{2}}{2}\right|t_{0}^{2}\right)\Bigg]dx
=
\int_{\Omega}\Phi\left(\frac{u_{0}^{2}+|\nabla u_{0}|^{2}}{2}t_{0}^{2}\right)dx.
\end{eqnarray}
Since $\Phi\in C^{1}([0,+\infty),\mathbb{R})$, by the mean value theorem and $(\phi_1)$, we mention that
there is
\begin{eqnarray*}
\varsigma\in \left[0,\left|\frac{u_{n}^{2}+|\nabla u_{n}|^{2}}{2}-\frac{u_{0}^{2}+|\nabla u_{0}|^{2}}{2}\right|t_{0}^{2}\right]
\end{eqnarray*}
such that
\begin{eqnarray}\label{1.1.8}
\Phi\left(\left|\frac{u_{n}^{2}+|\nabla u_{n}|^{2}}{2}-\frac{u_{0}^{2}+|\nabla u_{0}|^{2}}{2}\right|t_{0}^{2}\right)
&= &
t_{0}^{2}\phi(\varsigma)\left|\frac{u_{n}^{2}+|\nabla u_{n}|^{2}}{2}-\frac{u_{0}^{2}+|\nabla u_{0}|^{2}}{2}\right|\nonumber\\
&\geq &
\frac{\rho_{0}}{2}t_{0}^{2}(u_{n}^{2}-u_{0}^{2}+|\nabla u_{n}|^{2}-|\nabla u_{0}|^{2}).
\end{eqnarray}
Applying (\ref{1.1.14}) and the Brezis-Lieb Lemma \cite{Brezis1983}, we deduce that
\begin{eqnarray*}
\lim_{n\rightarrow\infty}\int_{\Omega}(u_{n}^{2}-|u_{n}^{2}-u_{0}^{2}|)dx
=\lim_{n\rightarrow\infty}\int_{\Omega}u_{0}^{2}dx,
 \end{eqnarray*}
which together with (\ref{582}) implies that
\begin{eqnarray}\label{1.1.10}
\lim_{n\rightarrow\infty}\int_{\Omega}(u_{n}^{2}-u_{0}^{2})dx
=\lim_{n\rightarrow\infty}\int_{\Omega}|u_{n}^{2}-u_{0}^{2}|dx
=0.
\end{eqnarray}
Moreover, applying (\ref{1.1.13}) and the Brezis-Lieb Lemma \cite{Brezis1983} one more time, we have
\begin{eqnarray}\label{1.1.11}
\lim_{n\rightarrow\infty}\int_{\Omega}\left(|\nabla u_{n}|^{2}-|\nabla u_{n}-\nabla u_{0}|^{2}\right)dx=\int_{\Omega}|\nabla u_{0}|^{2}dx.
\end{eqnarray}
Hence, combing with (\ref{1.1.7}), (\ref{1.1.8}), (\ref{1.1.10}), (\ref{1.1.11}) and (\ref{1.1.3}), we obtain that
\begin{eqnarray*}
\int_{\Omega}\Phi\left(\frac{u_{0}^{2}+|\nabla u_{0}|^{2}}{2}t_{0}^{2}\right)dx
\leq
\lim_{n\rightarrow\infty}\int_{\Omega}\Phi\left(\frac{u_{n}^{2}+|\nabla u_{n}|^{2}}{2}t_{0}^{2}\right)dx
-\frac{\rho_{0}}{2}t_{0}^{2}\delta
<
\lim_{n\rightarrow\infty}\int_{\Omega}\Phi\left(\frac{u_{n}^{2}+|\nabla u_{n}|^{2}}{2}t_{0}^{2}\right)dx,
\end{eqnarray*}
which together with (\ref{584}) and (\ref{2.0.4}) implies that
\begin{eqnarray}\label{585}
& &\gamma_{u_{0}}(t_{0})\nonumber\\
&< &
\lim_{n\rightarrow\infty}  \int_{\Omega}\Phi\left(\frac{u_{n}^{2}+|\nabla u_{n}|^{2}}{2}t_{0}^{2}\right)dx
-\lim_{n\rightarrow\infty}\frac{\lambda}{q+1}t_{0}^{q+1}\int_{\Omega}a(x)|u_{n}|^{q+1}dx
-\lim_{n\rightarrow\infty}\frac{1}{p+1}t_{0}^{p+1}\int_{\Omega}b(x)|u_{n}|^{p+1}dx\nonumber\\
&= &
 \lim_{n\rightarrow\infty}  \left\{\int_{\Omega}\Phi\left(\frac{u_{n}^{2}+|\nabla u_{n}|^{2}}{2}t_{0}^{2}\right)dx
-\frac{\lambda}{q+1}t_{0}^{q+1}\int_{\Omega}a(x)|u_{n}|^{q+1}dx
-\frac{1}{p+1}t_{0}^{p+1}\int_{\Omega}b(x)|u_{n}|^{p+1}dx\right\}\nonumber\\
&= &
 \lim_{n\rightarrow\infty}\gamma_{u_{n}}(t_{0}).
 \end{eqnarray}
In addition, since $\{u_{n}\}\subset N_{\lambda}^{-}$,
it follows from Lemma 3.4.(ii) that $t=1$ is a global maximum point of $\gamma_{u_{n}}$ such that $\gamma_{u_{n}}(1)>0$ for any $n\in N^{+}$.
We also have that $\gamma_{u_{n}}(t_{0})\leq \gamma_{u_{n}}(1)$.
Taking the limits in estimate above we get
\begin{eqnarray}\label{586} \lim_{n\rightarrow\infty}\gamma_{u_{n}}(t_{0})\leq\lim_{n\rightarrow\infty}\gamma_{u_{n}}(1).
\end{eqnarray}
Now, it follows from (\ref{1.1.1}), (\ref{585}) and (\ref{586}) that
$J_{\lambda}(t_{0}u_{0})<\inf_{u\in N_{\lambda}^{-}}J_{\lambda}(u)$.
Clearly, this is a contradiction for the definition of $\inf_{u\in N_{\lambda}^{-}}J_{\lambda}(u)$.
Thus, $u_{n}\rightarrow u_{0} $ in $H^{1}_{0}(\Omega)$ holds.
\par
Next we shall prove that $u_{0}\in N_{\lambda}^{-}$.
Firstly, by (\ref{1.1.1}), (\ref{1.1.2}), (\ref{1.1.53}) and the continuity of $J_{\lambda}$ and $J'_{\lambda}$, we have
$$J_{\lambda}(u_{0})=\inf_{N_{\lambda}^{-}}J_{\lambda}(u)>0
\;\mbox{and}\;
J'_{\lambda}(u_{0})=0.$$
Then $\int_{\Omega}b(x)|u_{0}|^{p+1}dx>0$ implies that $u_{0}\not\equiv0$.
So $u_{0}\in N_{\lambda}$.
Note that $N_{\lambda}=N_{\lambda}^{+}\cup N_{\lambda}^{-}$ when $\lambda\in (0,\lambda_{0})$.
If $u_{0}\in N_{\lambda}^{+}$, by Lemma 3.4 (i), we have $J_{\lambda}(u_{0})<0$, which is a contradiction for
$J_{\lambda}(u_{0})=\inf_{N_{\lambda}^{-}}J_{\lambda}(u)>0$.
Thus, $u_{0}\in N_{\lambda}^{-}$ and $u_{0}$ is a local minimum point of $J_{\lambda}$.
Then by Lemma 2.3, we mention that $u_{0}$ is a non-trivial weak solution for (\ref{eq1}).
So we finish the proof for the first case.
\par
Next, we shall prove that problem (\ref{eq1}) admits at least one non-trivial weak solution $\tilde{u}_{0}\in N_{\lambda}^{+}$ which satisfies
$J_{\lambda}(\tilde{u}_{0})<0$ for any $0<\lambda<\lambda_{0}$.
By Proposition 4.1, there exists a minimizing sequence $\{u_{n}\}\subset N_{\lambda}^{+}$ such that
\begin{eqnarray}
\label{1.1.20}
& & \lim_{n\rightarrow\infty}J_{\lambda}(u_{n})= \inf_{u\in N_{\lambda}^{+}}J_{\lambda}(u)<0,\\
\label{1.1.21}
& &  \lim_{n\rightarrow\infty}J'_{\lambda}(u_{n})=0.
\end{eqnarray}
Note that Proposition 2.2 implies that the functional $J_{\lambda}$ is coercive on $N_{\lambda}^{+}\subset N_{\lambda}$.
Thus, $\{u_{n}\}$ is a bounded sequence in $H^{1}_{0}(\Omega)$, that is, there exists $\tilde{M}>0$ such that
\begin{eqnarray}\label{1.1.22}
\|\nabla u_{n}\|_{2}\leq \tilde{M} \;\mbox{for all} \;n\in \mathbb{N}.
 \end{eqnarray}
Without any loss of generality, we assume that $u_{n}\rightharpoonup \tilde{u}_{0} $ in $H^{1}_{0}(\Omega)$.
It is easy to see that
\begin{eqnarray}\label{1.1.41}
\liminf_{n\rightarrow\infty}\|\nabla u_{n}\|_{2}^{2}\geq \|\nabla \tilde{u}_{0}\|_{2}^{2}.
 \end{eqnarray}
Using Proposition 2.1 and \cite[Theorem 4.9, p.94]{Brezis2010}, we deduce that
\begin{eqnarray}
\label{1.1.23} & & u_{n}\rightarrow \tilde{u}_{0} \;\mbox{in} \;L^{\alpha}(\Omega)  \;\mbox{for}\;1<\alpha<2^{\ast},\\
\label{1.1.24}& &  u_{n}\rightarrow \tilde{u}_{0} \;\mbox{a.e.} \;\mbox{in}\; \Omega.
\end{eqnarray}
Analogously, by (\ref{1.1.21}), (\ref{1.1.23}), $(\phi_1)$, $(\phi_6)$, $(\phi_7)$, $(H)$ and \cite[Lemma 3.5, p.11]{Jeanjean2022}, we have
\begin{eqnarray}\label{1.1.25}
\nabla u_{n}(x)\rightarrow \nabla \tilde{u}_{0}(x) \;\mbox{a.e.}\; \;\mbox{in}\;\Omega.
 \end{eqnarray}
Additionally, by (\ref{1.1.23}), $1<q+1<2<p+1<2^{\ast}$ and $(H)$, we have
\begin{eqnarray}
\label{576}& &\lim_{n\rightarrow\infty}\int_{\Omega}a(x)|u_{n}|^{q+1}dx = \int_{\Omega}a(x)|\tilde{u}_{0}|^{q+1}dx,\\
\label{1.1.33}& &\lim_{n\rightarrow\infty}\int_{\Omega}b(x)|u_{n}|^{p+1}dx =  \int_{\Omega}b(x)|\tilde{u}_{0}|^{p+1}dx.
\end{eqnarray}
\par
At this stage, we shall prove that
\begin{eqnarray}\label{1.1.51}
u_{n}\rightarrow \tilde{u}_{0} \;\mbox{in}\; H^{1}_{0}(\Omega).
\end{eqnarray}
The proof follows arguing by contradiction.
Assume that $u_{n}\not\rightarrow \tilde{u}_{0}$ in $H^{1}_{0}(\Omega)$.
Then there exist a subsequence of $\{u_{n}\}$, still denoted by $\{u_{n}\}$, and a positive constant $\tilde{\delta}>0$ such that
\begin{eqnarray}\label{1.1.44}
\lim_{n\rightarrow\infty}\int_{\Omega}|\nabla u_{n}-\nabla \tilde{u}_{0}|^{2}dx\geq \tilde{\delta}>0.
 \end{eqnarray}
Define a function $j_{2}:\mathbb{R}\rightarrow\mathbb{R}$ by
\begin{eqnarray*}
j_{2}(s)=2\phi(|s|)|s|, \;\;\forall s\in \mathbb{R}.
\end{eqnarray*}
It follows that $j_{2}(0)=0$.
Let $t\in(0,+\infty)$ be fixed.
Taking $\tilde{f}_{n}=\frac{|u_{n}|^{2}+|\nabla u_{n}|^{2}}{2}t^{2}$ and $\tilde{f}=\frac{\tilde{u}_{0}^{2}+|\nabla \tilde{u}_{0}|^{2}}{2}t^{2}$.
By $(\phi_1)$, (\ref{1.1.22}) and Proposition 2.1, we have
\begin{eqnarray*}
\int_{\Omega}|j_{2}(\tilde{f}_{n}(x))|dx
& =&
\int_{\Omega}\left|\phi\left(\frac{u_{n}^{2}+|\nabla u_{n}|^{2}}{2}t^{2}\right)(u_{n}^{2}+|\nabla u_{n}|^{2})t^{2}\right|dx\nonumber\\
& \leq&
t^{2}\rho_{1}(1+S_{2}^{2})\tilde{M}^{2}<+\infty.
\end{eqnarray*}
It follows from (\ref{1.1.24}) and (\ref{1.1.25}) that $\tilde{f}_{n}(x)\rightarrow \tilde{f}(x)$ a.e. in $\Omega$.
Hence, applying the Brezis-Lieb Lemma \cite{Brezis1983}, we have
\begin{eqnarray}\label{1.1.26}
& & \lim_{n\rightarrow\infty}
\int_{\Omega}\phi\left(\left|\frac{u_{n}^{2}+|\nabla u_{n}|^{2}}{2}-\frac{\tilde{u}_{0}^{2}+|\nabla \tilde{u}_{0}|^{2}}{2}\right|t^{2}\right)
\left|(u_{n}^{2}+|\nabla u_{n}|^{2})-(\tilde{u}_{0}^{2}+|\nabla \tilde{u}_{0}|^{2})\right|t^{2}dx\nonumber\\
&= &
\lim_{n\rightarrow\infty}\int_{\Omega}
\phi\left(\frac{u_{n}^{2}+|\nabla u_{n}|^{2}}{2}t^{2}\right)(u_{n}^{2}+|\nabla u_{n}|^{2})t^{2}dx
-\int_{\Omega}\phi\left(\frac{\tilde{u}_{0}^{2}+|\nabla \tilde{u}_{0}|^{2}}{2}t^{2}\right)(\tilde{u}_{0}^{2}+|\nabla \tilde{u}_{0}|^{2})t^{2}dx.
\end{eqnarray}
In view of $(\phi_1)$, we mention that
\begin{eqnarray*}
\phi\left(\left|\frac{u_{n}^{2}+|\nabla u_{n}|^{2}}{2}-\frac{\tilde{u}_{0}^{2}+|\nabla \tilde{u}_{0}|^{2}}{2}\right|t^{2}\right)
\left|(u_{n}^{2}+|\nabla u_{n}|^{2})-(\tilde{u}_{0}^{2}+|\nabla \tilde{u}_{0}|^{2})\right|t^{2}
\geq
\rho_{0}t^{2}(u_{n}^{2}-\tilde{u}_{0}^{2}+|\nabla \tilde{u}_{n}|^{2}-|\nabla \tilde{u}_{0}|^{2}).
\end{eqnarray*}
Applying (\ref{1.1.23}), (\ref{1.1.24}) and the Brezis-Lieb Lemma \cite{Brezis1983}, we deduce that
\begin{eqnarray}\label{1.1.29}
\lim_{n\rightarrow\infty}\int_{\Omega}(u_{n}^{2}-\tilde{u}_{0}^{2})dx
=\lim_{n\rightarrow\infty}\int_{\Omega}|u_{n}^{2}-\tilde{u}_{0}^{2}|dx
=0.
\end{eqnarray}
Moreover, applying (\ref{1.1.25}) and the Brezis-Lieb Lemma \cite{Brezis1983} one more time, we have
\begin{eqnarray}\label{1.1.30}
\lim_{n\rightarrow\infty}\int_{\Omega}\left(|\nabla u_{n}|^{2}-|\nabla u_{n}-\nabla \tilde{u}_{0}|^{2}\right)dx=\int_{\Omega}|\nabla \tilde{u}_{0}|^{2}dx.
\end{eqnarray}
Therefore, combing with (\ref{1.1.44}), (\ref{1.1.26}), (\ref{1.1.29}) and(\ref{1.1.30}), we obtain that
\begin{eqnarray*}
\int_{\Omega}\phi\left(\frac{\tilde{u}_{0}^{2}+|\nabla \tilde{u}_{0}|^{2}}{2}t^{2}\right)(\tilde{u}_{0}^{2}+|\nabla \tilde{u}_{0}|^{2})t^{2}dx
& \leq&
\lim_{n\rightarrow\infty}\int_{\Omega}\phi\left(\frac{u_{n}^{2}+|\nabla u_{n}|^{2}}{2}t^{2}\right)(u_{n}^{2}+|\nabla u_{n}|^{2})t^{2}dx-t^{2}\rho_{0}\delta\\
&< &
\lim_{n\rightarrow\infty}\int_{\Omega}\phi\left(\frac{u_{n}^{2}+|\nabla u_{n}|^{2}}{2}t^{2}\right)(u_{n}^{2}+|\nabla u_{n}|^{2})t^{2}dx.
\end{eqnarray*}
Multiplying the above expression by $t^{-1}$, we have
\begin{eqnarray}\label{1.1.60}
\int_{\Omega}t\phi\left(\frac{\tilde{u}_{0}^{2}+|\nabla \tilde{u}_{0}|^{2}}{2}t^{2}\right)(\tilde{u}_{0}^{2}+|\nabla \tilde{u}_{0}|^{2})dx
<
\lim_{n\rightarrow\infty}\int_{\Omega}t\phi\left(\frac{u_{n}^{2}+|\nabla u_{n}|^{2}}{2}t^{2}\right)(u_{n}^{2}+|\nabla u_{n}|^{2})dx.
\end{eqnarray}
Similarly, we also conclude that
\begin{eqnarray}\label{1.1.45}
\int_{\Omega}\Phi\left(\frac{u_{0}^{2}+|\nabla u_{0}|^{2}}{2}\right)dx
<
\lim_{n\rightarrow\infty}\int_{\Omega}\Phi\left(\frac{u_{n}^{2}+|\nabla u_{n}|^{2}}{2}\right)dx.
\end{eqnarray}
Then, we shall prove that $\int_{\Omega}a(x)|u_{0}|^{q+1}dx>0$.
Using (\ref{2.0.3}), $1<q+1<2<p+1$ and $(\phi_1)$, we easily see that
\begin{eqnarray*}
\int_{\Omega}a(x)|u_{n}|^{q+1}dx
\geq
\frac{1}{\lambda(\frac{1}{q+1}-\frac{1}{p+1})}\left\{\left(\frac{\rho_{0}}{2}-\frac{\rho_{1}}{p+1}\right)
\|\nabla u_{n}\|_{2}^{2}
-J_{\lambda}(u_{n})\right\}.
\end{eqnarray*}
Taking the limits in estimate above and combing with (\ref{576}), (\ref{1.1.41}) and (\ref{1.1.20}), we have
$\int_{\Omega}a(x)|\tilde{u}_{0}|^{q+1}dx>0$.
Then by $\clubsuit^{2}$ and $\clubsuit^{4}$ in Section 3, we know that there exists a critical point $\tilde{t}_{0}>0$ of $\gamma_{\tilde{u}_{0}}$ which satisfies
$\tilde{t}_{0}\tilde{u}_{0}\in N_{\lambda}^{+}$,
$\gamma'_{\tilde{u}_{0}}(\tilde{t}_{0})=0$ and $\gamma'_{\tilde{u}_{0}}<0$ for any $t\in(0,\tilde{t}_{0})$.
Next we shall prove that
\begin{eqnarray}\label{577}
\gamma_{\tilde{u}_{0}}(\tilde{t}_{0})<\gamma_{\tilde{u}_{0}}(1).
\end{eqnarray}
Since $\gamma_{\tilde{u}_{0}}$ is strictly decreasing on $(0,\tilde{t}_{0})$, we just prove that $\tilde{t}_{0}>1$.
On one hand, let $t=\tilde{t}_{0}$ in (\ref{1.1.60}), we have
\begin{eqnarray}\label{1.1.31}
\int_{\Omega}\tilde{t}_{0}\phi\left( \frac{\tilde{u}^{2}_{0}+|\nabla \tilde{u}_{0}|^{2}}{2}\tilde{t}^{2}_{0}\right)(\tilde{u}^{2}_{0}+|\nabla \tilde{u}_{0}|^{2})dx
<
\lim_{n\rightarrow\infty}\int_{\Omega}\tilde{t}_{0}\phi\left( \frac{u_{n}^{2}+|\nabla u_{n}|^{2}}{2}\tilde{t}^{2}_{0}\right)(u_{n}^{2}+|\nabla u_{n}|^{2})dx.
\end{eqnarray}
Then in view of (\ref{1.1.31}), (\ref{576}) and (\ref{1.1.33}), we have
\begin{eqnarray*}
& &\lim_{n\rightarrow\infty}\gamma'_{u_{n}}(\tilde{t}_{0})\\
&= &
\lim_{n\rightarrow\infty} \left\{\int_{\Omega}\tilde{t}_{0}\phi\left( \frac{u_{n}^{2}+|\nabla u_{n}|^{2}}{2}\tilde{t}^{2}_{0}\right)(u_{n}^{2}+|\nabla u_{n}|^{2})dx
-\lambda \tilde{t}^{q}_{0}\int_{\Omega}a(x)|u_{n}|^{q+1}dx
-\tilde{t}^{p}_{0}\int_{\Omega}b(x)|u_{n}|^{p+1}dx\right\}\\
&= &
\lim_{n\rightarrow\infty}  \int_{\Omega}\tilde{t}_{0}\phi\left( \frac{u_{n}^{2}+|\nabla u_{n}|^{2}}{2}\tilde{t}^{2}_{0}\right)(u_{n}^{2}+|\nabla u_{n}|^{2})dx
-\lim_{n\rightarrow\infty}\lambda \tilde{t}^{q}_{0}\int_{\Omega}a(x)|u_{n}|^{q+1}dx
-\lim_{n\rightarrow\infty}\tilde{t}^{p}_{0}\int_{\Omega}b(x)|u_{n}|^{p+1}dx\\
&> &
\int_{\Omega}\tilde{t}_{0}\phi\left( \frac{\tilde{u}_{0}^{2}+|\nabla\tilde{u}_{0}|^{2}}{2}\tilde{t}_{0}^{2}\right)(\tilde{u}_{0}^{2}+|\nabla \tilde{u}_{0}|^{2})dx
-\lambda \tilde{t}_{0}^{q}\int_{\Omega}a(x)|\tilde{u}_{0}|^{q+1}dx -\tilde{t}_{0}^{p}\int_{\Omega}b(x)|\tilde{u}_{0}|^{p+1}dx\\
&= &\gamma'_{\tilde{u}_{0}}(\tilde{t}_{0})\\
&= &0.
 \end{eqnarray*}
Hence, $\gamma'_{u_{n}}(\tilde{t}_{0})>0$ for any $n\geq n_{0}$ where $n_{0}\in N^{+}$ is large enough.
On the other hand, since $\{u_{n}\}\subset N_{\lambda}^{+}$, it follows from Lemma 3.4.(i) that $\gamma'_{u_{n}}(1)=0$ for any $n\in N^{+}$ and
$\gamma'_{u_{n}}(t)<0$ for any $t\in (0,1)$.
As a consequence, we obtain that $\tilde{t}_{0}>1$, which implies that (\ref{577}) holds.
In addition, by (\ref{2.0.4}), (\ref{1.1.45}), (\ref{576}) and (\ref{1.1.33}), we mention that
\begin{eqnarray}\label{579}
\gamma_{\tilde{u}_{0}}(1)
&= & \int_{\Omega}\Phi\left(\frac{\tilde{u}_{0}^{2}+|\nabla \tilde{u}_{0}|^{2}}{2}\right)dx -\frac{\lambda}{q+1}\int_{\Omega}a(x)|\tilde{u}_{0}|^{q+1}dx -\frac{1}{p+1}\int_{\Omega}b(x)|\tilde{u}_{0}|^{p+1}dx\nonumber \\
&< &
\lim_{n\rightarrow\infty}  \int_{\Omega}\Phi\left(\frac{u_{n}^{2}+|\nabla u_{n}|^{2}}{2}\right)dx
-\lim_{n\rightarrow\infty}\frac{\lambda}{q+1}\int_{\Omega}a(x)|u_{n}|^{q+1}dx
-\lim_{n\rightarrow\infty}\frac{1}{p+1}\int_{\Omega}b(x)|u_{n}|^{p+1}dx\nonumber \\
&=&
\lim_{n\rightarrow\infty}  \left\{\int_{\Omega}\Phi\left(\frac{u_{n}^{2}+|\nabla u_{n}|^{2}}{2}\right)dx
-\frac{\lambda}{q+1}\int_{\Omega}a(x)|u_{n}|^{q+1}dx
-\frac{1}{p+1}\int_{\Omega}b(x)|u_{n}|^{p+1}dx\right\}\nonumber \\
&= &
\lim_{n\rightarrow\infty}\gamma_{u_{n}}(1)\nonumber \\
&= &
\inf_{N_{\lambda}^{+}}J_{\lambda}(u).
 \end{eqnarray}
Thus, it follows from (\ref{577}) and (\ref{579}) that
\begin{eqnarray*}
J_{\lambda}(\tilde{t}_{0}\tilde{u}_{0})
=
\gamma_{\tilde{u}_{0}}(\tilde{t}_{0})
<
\gamma_{\tilde{u}_{0}}(1)
<
\inf_{N_{\lambda}^{+}}J_{\lambda}(u).
 \end{eqnarray*}
Clearly, this is a contradiction for the definition of $\inf_{u\in N_{\lambda}^{+}}J_{\lambda}(u)$.
Hence, $u_{n}\rightarrow \tilde{u}_{0}$ in $H^{1}_{0}(\Omega)$.
\par
Next, we prove that $\tilde{u}_{0}\in N_{\lambda}^{+}$.
 Firstly, by (\ref{1.1.51}), (\ref{1.1.20}), (\ref{1.1.21}) and the continuity of $J_{\lambda}$ and $J'_{\lambda}$, we have
$$J_{\lambda}(\tilde{u}_{0})=\inf_{N_{\lambda}^{+}}J_{\lambda}(u)<0
\;\mbox{and}\;
J'_{\lambda}(\tilde{u}_{0})=0.$$
Then $\int_{\Omega}a(x)|\tilde{u}_{0}|^{q+1}dx>0$ implies that $\tilde{u}_{0}\not\equiv0$.
So $\tilde{u}_{0}\in N_{\lambda}$.
Note that $N_{\lambda}=N_{\lambda}^{+}\cup N_{\lambda}^{-}$ when $\lambda\in (0,\lambda_{0})$.
If $\tilde{u}_{0}\in N_{\lambda}^{-}$, by Lemma 3.4 (ii), we have $J_{\lambda}(\tilde{u}_{0})>0$.
Clearly, this is a contradiction for $J_{\lambda}(\tilde{u}_{0})=\inf_{N_{\lambda}^{-}}J_{\lambda}(u)<0$.
Therefore, $\tilde{u}_{0}\in N_{\lambda}^{+}$ and $\tilde{u}_{0}$ is a local minimum point of $J_{\lambda}$.
Then by Lemma 2.3, we obtain that $\tilde{u}_{0}$ is a non-trivial weak solution of (\ref{eq1}).
Note that $N_{\lambda}^{+}\cap N_{\lambda}^{-}=\emptyset$, so $u_{0}\not\equiv \tilde{u}_{0}$.
Then by
\begin{eqnarray*}
J_{\lambda}(\tilde{u}_{0})< 0<J_{\lambda}(u_{0}),
 \end{eqnarray*}
we obtain that $\tilde{u}_{0}$ is a ground state solution of (\ref{eq1}).
\qed

\section{Example}
\par
Let $\phi(s)=(1+s)^{-3}+A$, $s\in[0,+\infty)$, where
\begin{eqnarray*}
A
>
\max\left\{\frac{ q+1}{1-q}, \frac{2}{p+1}, 5, \frac{81}{128(1-q)}, \frac{1}{(1-q)(p-1)}\left[\frac{5184}{3125}+(p+q)\frac{81}{128}\right]\right\}.
 \end{eqnarray*}
 \par
Initially, we  prove that $(\phi_1)$ holds.
If $(p+1)(q+1)>4$, by $A>\frac{ q+1}{1-q}$, we have $ A+1<\frac{2 }{q+1}A$.
Taking $\rho_{1}\in\left[A+1,\frac{2 }{q+1}A\right)$, we conclude that
\begin{eqnarray*}
0<\frac{q+1}{2}\rho_{1}<A< A+1\leq \rho_1.
 \end{eqnarray*}
Then taking $\rho_{0}\in\left(\frac{q+1}{2}\rho_{1},A\right)$ we deduce that
\begin{eqnarray*}
0<\max\left\{\frac{q+1}{2}, \frac{2}{p+1}\right\}\rho_{1}=\frac{q+1}{2}\rho_{1}<\rho_0<A< A+1\leq \rho_1,
 \end{eqnarray*}
which together with
\begin{eqnarray*}
A<\phi(s)\leq A+1,\;\;s\geq0,
 \end{eqnarray*}
implies that $(\phi_1)$ holds.
If $(p+1)(q+1)<4$, in view of $A>\frac{ 2}{p+1}$, we mention that $ A+1<\frac{ 2}{p+1}A$.
Taking $\rho_{1}\in\left[A+1,\frac{ 2}{p+1}\right)$, we obtain that
\begin{eqnarray*}
0<\frac{ 2}{p+1}\rho_{1}<A< A+1\leq \rho_1.
 \end{eqnarray*}
Then let $\rho_{0}\in\left(\frac{ 2}{p+1},\rho_{1}\right)$. We deduce that
\begin{eqnarray*}
0<\max\left\{\frac{q+1}{2}, \frac{2}{p+1}\right\}\rho_{1}=\frac{ 2}{p+1}\rho_{1}<\rho_0<A< A+1\leq \rho_1,
 \end{eqnarray*}
which together with
\begin{eqnarray*}
A<\phi(s)\leq A+1,\;\;s\geq0,
 \end{eqnarray*}
implies that $(\phi_1)$ holds.
\par
Next, we  prove that $(\phi_2)$ holds.
It is easy to see that
\begin{eqnarray*}
\phi'(s)=-3(1+s)^{-4},
\;\;
\phi''(s)=12(1+s)^{-5},
\;\;
s\in[0,+\infty).
\end{eqnarray*}
Let
\begin{eqnarray*}
f_{4}(s)=s(1+s)^{-4},
\;\;
f_{5}(s)=s^{2}(1+s)^{-5},
\;\;s\in[0,+\infty).
\end{eqnarray*}
It is not hard to verify that
\begin{eqnarray}\label{402}
& &f_{4}\left(\frac{1}{3}\right)=\max_{s\geq0}f_{4}(s)=\frac{27}{256},\;\;f_{4}(0)=\min_{s\geq0}f_{4}(s)=0,\nonumber\\
& &f_{5}\left(\frac{2}{3}\right)=\max_{s\geq0}f_{5}(s)=\frac{108}{3125},\;\;f_{5}(0)=\min_{s\geq0}f_{5}(s)=0.
\end{eqnarray}
Then
\begin{eqnarray*}
|\phi'(s)|s+|\phi''(s)|s^{2}\leq\frac{584901}{800000}.
\end{eqnarray*}
We can see that $(\phi_2)$ holds.
\par
Next, we  prove that $(\phi_3)$ holds.
Using the fact that $A>\frac{81}{128(1-q)}$ and (\ref{402}), we mention that
\begin{eqnarray*}
(1-q)(A+1)
&\geq&
(1-q)(1+s)^{-3}+(1-q)A-6s(1+s)^{-4}\\
& =&
(1-q)\phi(s)+2\phi'(s)s\\
& >&
(1-q)A-6s(1+s)^{-4}\\
& \geq&
(1-q)A-\frac{81}{128}
>0.
\end{eqnarray*}
Thus, $(\phi_3)$ holds.
\par
Next, we prove that $(\phi_4)$ holds.
Taking into account that $A>0$ and $p>1$, we conclude that
\begin{eqnarray*}
2\phi'\left( s\right)s+(1-p)\phi\left( s\right)
& =&
(1-p)(1+t)^{-3}+(1-p)A-6t(1+t)^{-4}\\
&<&
-(p-1)A\\
&<&
0.
\end{eqnarray*}
Thus, $(\phi_4)$ holds.
\par
Next, we prove that $(\phi_5)$ holds.
It follows from $0<q<1<p$ and (\ref{402}) that
\begin{eqnarray*}
A
>
\frac{1}{(1-q)(p-1)}\left[\frac{5184}{3125}+(p+q)\frac{81}{128}\right].
\end{eqnarray*}
As a consequence, we have
\begin{eqnarray*}
& &(1-q)(1-p)\phi\left( s\right)+2(4-p-q)\phi'\left( s\right)s+4\phi''\left( s\right)s^{2}\\
& =&
 6(p+q)s(1+s)^{-4}+48s^{2}(1+s)^{-5}-24s(1+s)^{-4}-(1-q)(p-1)(1+s)^{-3}-(1-q)(p-1)A\\
& <&
(p+q)\frac{81}{128}+\frac{5184}{3125}-(1-q)(p-1)A\\
& <&
0,
\end{eqnarray*}
which implies that $(\phi_5)$ holds.
\par
Next, we  prove that $(\phi_6)$ holds.
Let
\begin{eqnarray*}
f_{6}(s)=\Phi\left(s^{2}\right)
=As^{2}-\frac{1}{2}\left(1+s^{2}\right)^{-2}+\frac{1}{2},\;\;s\in[0,+\infty).
 \end{eqnarray*}
It follows that
\begin{eqnarray*}
f_{6}'(s)=2As+2s\left(1+s^{2}\right)^{-3},
\;\;
f_{6}''(s)=2A+2\left(1+s^{2}\right)^{-3}-12s^{2}\left(1+s^{2}\right)^{-4}.
 \end{eqnarray*}
Due to the fact that $A>5$, we have
\begin{eqnarray*}
f_{6}''(s)
& =&
2[A+(1+s^{2})^{-3}(1-6s^{2}(1+s^{2})^{-1})]\\
& >&
2(A-5)\\
& >&
0,
 \end{eqnarray*}
which implies that $(\phi_6)$ holds.
\par
Finally, it follows that
\begin{eqnarray*}
\lim_{s\rightarrow+\infty}\phi(s)=\lim_{s\rightarrow+\infty}((1+s)^{-3}+A)=A>0
 \end{eqnarray*}
which implies that $(\phi_7)$ holds.
\section{Appendix A}
\vskip2mm
 \noindent
{\bf Appendix A.1.} {\it Define $J_{\lambda}:H_{0}^{1}(\Omega)\rightarrow\mathbb{R}$ given by
\begin{eqnarray*}
J_{\lambda}(u)
=
\int_{\Omega}\Phi\left(\frac{u^{2}+|\nabla u|^{2}}{2}\right)dx
-\frac{\lambda}{q+1}\int_{\Omega}a(x)|u|^{q+1}dx-\frac{1}{p+1}\int_{\Omega}b(x)|u|^{p+1}dx,
\;\;
\forall u\in H_{0}^{1}(\Omega).
\end{eqnarray*}
Assume that $(\phi_1)$ and $(H)$ hold.
Then $J_{\lambda}\in C^{1}(H_{0}^{1}(\Omega),\mathbb{R})$ with derivative $J'_{\lambda}:H_{0}^{1}(\Omega)\rightarrow (H_{0}^{1}(\Omega))'$ given by
\begin{eqnarray}\label{500}
\langle J'_{\lambda}(u),v\rangle
=
\int_{\Omega}\phi\left(\frac{u^{2}+|\nabla u|^{2}}{2}\right)(uv+\nabla u\cdot \nabla v )dx
-\lambda\int_{\Omega}a(x)|u|^{q-1}uvdx-\int_{\Omega}b(x)|u|^{p-1}uvdx,
\end{eqnarray}
$\forall u,v\in H_{0}^{1}(\Omega)$. \\}
{\bf Proof.} Defining $F:\Omega\times\mathbb{R}\times\mathbb{R}^{N}\rightarrow\mathbb{R}$ by
\begin{eqnarray*}
F(x,u,\xi)
=
\Phi\left(\frac{u^{2}+|\xi|^{2}}{2}\right)
-\frac{\lambda}{q+1}a(x)|u|^{q+1}
-\frac{1}{p+1}b(x)|u|^{p+1}.
\end{eqnarray*}
It follows that
\begin{eqnarray*}
F_{u}(x,u,\xi)
=
\phi\left(\frac{u^{2}+|\xi|^{2}}{2}\right)u
-\lambda a(x)|u|^{q-1}u
-b(x)|u|^{p-1}u,\;\;
F_{\xi}(x,u,\xi)
=
\phi\left(\frac{u^{2}+|\xi|^{2}}{2}\right)\xi.
\end{eqnarray*}
By $(\phi_1)$ and $(H)$, we have
\begin{eqnarray*}
|F(x,u,\xi)|
\leq
\frac{\rho_{1}}{2}(u^{2}+|\xi|^{2})
+\frac{\lambda}{q+1}\|a\|_{\infty}|u|^{q+1}
+\frac{1}{p+1}\|b\|_{\infty}|u|^{p+1},
\end{eqnarray*}
combing with $1<q+1<2<p+1<2^{\ast}$, we mention that
\begin{eqnarray*}
& &
|F(x,u,\xi)|
\leq
C(1+|\xi|^{2}),\;\mbox{if}\; |u|\leq1,\\
& &
|F(x,u,\xi)|
\leq
C(|u|^{p+1}+|\xi|^{2}),\;\mbox{if}\; |u|>1,
\end{eqnarray*}
which imply that
\begin{eqnarray*}
|F(x,u,\xi)|
\leq
C(1+|u|^{p+1}+|\xi|^{2}).
\end{eqnarray*}
Using $(\phi_1)$ and $(H)$ one more time, we have
\begin{eqnarray*}
|F_{u}(x,u,\xi)|
\leq
\rho_{1}|u|+\lambda\|a\|_{\infty}|u|^{q}+\|b\|_{\infty}|u|^{p},
\end{eqnarray*}
which together with $1<q+1<2<p+1<2^{\ast}$ implies that
\begin{eqnarray*}
|F_{u}(x,u,\xi)|
\leq
C(1+|u|^{p}).
\end{eqnarray*}
Moreover, in view of $(\phi_1)$, we deduce that
\begin{eqnarray*}
|F_{\xi}(x,u,\xi)|
\leq
\rho_{1}|\xi|.
\end{eqnarray*}
Therefore, applying \cite[Theorem C.1]{Struwe1996}, $J_{\lambda}$ is of class $C^{1}$. Moreover, $J'_{\lambda}:H_{0}^{1}(\Omega)\rightarrow (H_{0}^{1}(\Omega))'$ is given by
\begin{eqnarray*}
\langle J'_{\lambda}(u),v\rangle
=
\int_{\Omega}\phi\left(\frac{u^{2}+|\nabla u|^{2}}{2}\right)(uv+\nabla u\cdot \nabla v )dx
-\lambda\int_{\Omega}a(x)|u|^{q-1}uvdx-\int_{\Omega}b(x)|u|^{p-1}uvdx,
\end{eqnarray*}
$\forall u,v\in H_{0}^{1}(\Omega)$.

\qed
\vskip2mm
 \noindent
{\bf Appendix A.2.} {\it Suppose that $(\phi_1)$, $(\phi_2)$ and $(H)$ hold.
Then $\gamma_{u}\in C^{2}((0,+\infty),\mathbb{R})$  with the first order derivative $\gamma_{u}':(0,+\infty)\rightarrow \mathbb{R}$
and the second order derivative $\gamma_{u}'':(0,+\infty)\rightarrow \mathbb{R}$ given by
\begin{eqnarray}\label{520}
\gamma_{u}'(t)
=
\int_{\Omega} t\phi\left( \frac{u^{2}+|\nabla u|^{2}}{2}t^{2}\right)(u^{2}+|\nabla u|^{2})dx
-\lambda t^{q}\int a(x)|u|^{q+1}dx-t^{p}\int b(x)|u|^{p+1}dx.
 \end{eqnarray}
and
\begin{eqnarray}\label{521}
\gamma_{u}''(t)
& := &
\int_{\Omega}\left[ \phi\left( \frac{u^{2}+|\nabla u|^{2}}{2}t^{2}\right)(u^{2}+|\nabla u|^{2})
+t^{2}\phi'\left( \frac{u^{2}+|\nabla u|^{2}}{2}t^{2}\right)(u^{2}+|\nabla u|^{2})^{2}\right]dx\nonumber\\
& &- \lambda q t^{q-1}\int_{\Omega} a(x)|u|^{q+1}dx
-pt^{p-1}\int_{\Omega}b(x)|u|^{p+1}dx
\end{eqnarray}
respectively.\\}
{\bf Proof.} Firstly, we present the proof for (\ref{520}).
It follows that
\begin{eqnarray*}
\gamma_{u}'(t)
& = &
\lim_{h\rightarrow 0} \frac{\gamma_{u}(t+h)-\gamma_{u}(t)}{h}\\
& = &
\lim_{h\rightarrow 0} \frac{J(tu+hu)-J(tu)}{h}\\
& = &\langle J'(tu),u\rangle,
\end{eqnarray*}
which together with (\ref{500}) implies that (\ref{520}) holds.
\par
Next, we  prove that $\gamma_{u}'$ is continuous.
Let $t_{n},t\in(0,+\infty)$ such that $t_{n}\rightarrow t$ in $\mathbb{R}$.
It is easily seen that there exists $M>0$ such that
\begin{eqnarray}\label{522}
|t_{n}|\leq M,\;\;\forall n\in \mathbb{N}.
\end{eqnarray}
Furthermore, it follows that
\begin{eqnarray}\label{523}
|t^{q}-t_{n}^{q}| \rightarrow0,
\;\;
|t^{p}-t_{n}^{p}|\rightarrow0.
\end{eqnarray}
Then in view of (\ref{520}), $(H)$, $1<q+1<2<p+1<2^{\ast}$ and Proposition 2.1, we have
\begin{eqnarray}\label{524}
&&|\gamma_{u}'(t_{n})-\gamma_{u}'(t)|\nonumber\\
&= &
\Big|\int_{\Omega} \left[ t_{n}\phi\left( \frac{u^{2}+|\nabla u|^{2}}{2} t_{n}^{2}\right)(u^{2}+|\nabla u|^{2})
-t\phi\left( \frac{u^{2}+|\nabla u|^{2}}{2}t^{2}\right)(u^{2}+|\nabla u|^{2})\right]dx\nonumber\\
&  &
+(t^{q}-t_{n}^{q}) \lambda\int_{\Omega} a(x)|u|^{q+1}dx
+ (t^{p}-t_{n}^{p})\int_{\Omega} b(x)|u|^{p+1}dx\Big|\nonumber\\
&\leq &
\int_{\Omega} \left| t_{n}\phi\left( \frac{u^{2}+|\nabla u|^{2}}{2} t_{n}^{2}\right)(u^{2}+|\nabla u|^{2})
-t\phi\left( \frac{u^{2}+|\nabla u|^{2}}{2}t^{2}\right)(u^{2}+|\nabla u|^{2})\right|dx\nonumber\\
& &
+ |t^{q}-t_{n}^{q}| \lambda \|a\|_{\infty}S^{q+1}_{q+1}\|\nabla u\|_{2}^{q+1}
+ |t^{p}-t_{n}^{p}|\|b\|_{\infty}S^{p+1}_{p+1}\|\nabla u\|_{2}^{p+1}.
 \end{eqnarray}
By $(\phi_{1})$, (\ref{522}) and Proposition 2.1, we have
\begin{eqnarray}\label{525}
&  &
\left| t_{n}\phi\left( \frac{u^{2}+|\nabla u|^{2}}{2} t_{n}^{2}\right)(u^{2}+|\nabla u|^{2})
-t\phi\left( \frac{u^{2}+|\nabla u|^{2}}{2}t^{2}\right)(u^{2}+|\nabla u|^{2})\right|\nonumber\\
&= &
\left| t_{n}\phi\left( \frac{u^{2}+|\nabla u|^{2}}{2} t_{n}^{2}\right)
-t\phi\left( \frac{u^{2}+|\nabla u|^{2}}{2}t^{2}\right)\right|(u^{2}+|\nabla u|^{2})\nonumber\\
&\leq &
\left( t_{n}\left|\phi\left( \frac{u^{2}+|\nabla u|^{2}}{2} t_{n}^{2}\right)\right|
+t\left|\phi\left( \frac{u^{2}+|\nabla u|^{2}}{2}t^{2}\right)\right|\right)(u^{2}+|\nabla u|^{2})\nonumber\\
&\leq &
(M+t)\rho_{1}(u^{2}+|\nabla u|^{2})
\in L^{1}(\Omega).
 \end{eqnarray}
Then by the continuity of $\phi$, we also have
\begin{eqnarray}\label{526}
\left| t_{n}\phi\left( \frac{u^{2}+|\nabla u|^{2}}{2} t_{n}^{2}\right)(u^{2}+|\nabla u|^{2})
-t\phi\left( \frac{u^{2}+|\nabla u|^{2}}{2}t^{2}\right)(u^{2}+|\nabla u|^{2})\right|\rightarrow0, \;\mbox{a.e.}\; x \in \Omega.
 \end{eqnarray}
Thus, it follows from (\ref{525}), (\ref{526}) and the Lebesgue dominated convergence theorem that
\begin{eqnarray}\label{527}
\lim_{n\rightarrow\infty}\int_{\Omega} \left| t_{n}\phi\left( \frac{u^{2}+|\nabla u|^{2}}{2} t_{n}^{2}\right)(u^{2}+|\nabla u|^{2})
-t\phi\left( \frac{u^{2}+|\nabla u|^{2}}{2}t^{2}\right)(u^{2}+|\nabla u|^{2})\right|dx=0,
 \end{eqnarray}
which together with (\ref{523}) and (\ref{524}) implies that $\gamma_{u}'$ is continuous.
\par
Next, we  consider the proof for (\ref{521}).
Let $f_{3}(h)=(t+h)\phi\left( \frac{u^{2}+|\nabla u|^{2}}{2}(t+h)^{2}\right)$, $0<|h|\leq1$, $t\in(0,+\infty)$, $u\in H_{0}^{1}(\Omega)$.
Since $\Phi\in C^{1}([0,+\infty),\mathbb{R})$, $f_{3}$ is of class $C^{1}$ and
\begin{eqnarray*}
f_{3}'(h)=\phi\left( \frac{u^{2}+|\nabla u|^{2}}{2}(t+h)^{2}\right)
+\phi'\left( \frac{u^{2}+|\nabla u|^{2}}{2}(t+h)^{2}\right)(t+h)^{2}(u^{2}+|\nabla u|^{2}).
 \end{eqnarray*}
Applying the mean value theorem, there is $\theta\in\mathbb{R}$ with $0<\theta<t\leq1$ such that
\begin{eqnarray}\label{528}
& &\frac{  (t+h)\phi\left( \frac{u^{2}+|\nabla u|^{2}}{2}(t+h)^{2}\right)
  - t\phi\left( \frac{u^{2}+|\nabla u|^{2}}{2}t^{2}\right)}{h}(u^{2}+|\nabla u|^{2})\nonumber\\
&= & \frac{f_{3}(h)-f_{3}(0)}{h-0}(u^{2}+|\nabla u|^{2})\nonumber\\
&= & f_{3}'(\theta)(u^{2}+|\nabla u|^{2})\nonumber\\
&= & \left[\phi\left( \frac{u^{2}+|\nabla u|^{2}}{2}(t+\theta)^{2}\right)
+\phi'\left( \frac{u^{2}+|\nabla u|^{2}}{2}(t+\theta)^{2}\right)(t+\theta)^{2}(u^{2}+|\nabla u|^{2})\right](u^{2}+|\nabla u|^{2}).
  \end{eqnarray}
Taking into account that $(\phi_1)$, $(\phi_2)$ and Proposition 2.1, we have
\begin{eqnarray}\label{529}
& &
|f_{3}'(\theta)|(u^{2}+|\nabla u|^{2})\nonumber\\
&\leq& \left(\left|\phi\left( \frac{u^{2}+|\nabla u|^{2}}{2}(t+\theta)^{2}\right)\right|
+2\left|\phi'\left( \frac{u^{2}+|\nabla u|^{2}}{2}(t+\theta)^{2}\right)\right|\frac{u^{2}+|\nabla u|^{2}}{2}(t+\theta)^{2}\right)
(u^{2}+|\nabla u|^{2})\nonumber\\
&\leq& (\rho_{1}+2\rho_{2})(u^{2}+|\nabla u|^{2})
\in L^{1}(\Omega).
  \end{eqnarray}
Moreover, by the continuity of $\phi'$, we mention that
\begin{eqnarray}\label{530}
&  &
\lim_{\theta\rightarrow 0}\left[\phi\left( \frac{u^{2}+|\nabla u|^{2}}{2}(t+\theta)^{2}\right)
+\phi'\left( \frac{u^{2}+|\nabla u|^{2}}{2}(t+\theta)^{2}\right)(t+\theta)^{2}(u^{2}+|\nabla u|^{2})\right](u^{2}+|\nabla u|^{2})\nonumber\\
& = & \phi\left( \frac{u^{2}+|\nabla u|^{2}}{2}t^{2}\right)(u^{2}+|\nabla u|^{2})+t^{2}\phi'\left( \frac{u^{2}
+|\nabla u|^{2}}{2}t^{2}\right)(u^{2}+|\nabla u|^{2})^{2},
\;\mbox{a.e.}\; x \in \Omega.
 \end{eqnarray}
Note that $\theta\rightarrow0$ as $t\rightarrow0$.
Hence, by (\ref{528}), (\ref{529}), (\ref{530}) and the Lebesgue dominated convergence theorem, we have
\begin{eqnarray*}
&  &
\lim_{h\rightarrow 0} \int_{\Omega} \frac{ (t+h)\phi\left( \frac{u^{2}+|\nabla u|^{2}}{2}(t+h)^{2}\right)
  - t\phi\left( \frac{u^{2}+|\nabla u|^{2}}{2}t^{2}\right)}{h} (u^{2}+|\nabla u|^{2})dx\\
& = &\int_{\Omega} \left( \phi\left( \frac{u^{2}+|\nabla u|^{2}}{2}t^{2}\right)(u^{2}+|\nabla u|^{2})+t^{2}\phi'\left( \frac{u^{2}
+|\nabla u|^{2}}{2}t^{2}\right)(u^{2}+|\nabla u|^{2})^{2}\right)dx,
 \end{eqnarray*}
which implies that (\ref{521}) holds.
\par
Finally, we  prove that $\gamma_{u}''$ is continuous.
Let $t_{n},t\in(0,+\infty)$ such that $t_{n}\rightarrow t$ in $\mathbb{R}$.
It is not hard to verify that
\begin{eqnarray}\label{532}
|t^{q-1}-t_{n}^{q-1}| \rightarrow0,
\;\;
|t^{p-1}-t_{n}^{p-1}|\rightarrow0.
\end{eqnarray}
In view of (\ref{521}), $(H)$, $1<q+1<2<p+1<2^{\ast}$ and Proposition 2.1, we have
\begin{eqnarray}\label{533}
&  &|\gamma_{u}''(t_{n})-\gamma_{u}''(t)|\nonumber\\
&= & |\int_{\Omega} \left[\phi\left( \frac{u^{2}+|\nabla u|^{2}}{2}t_{n}^{2}\right)-\phi\left( \frac{u^{2}+|\nabla u|^{2}}{2}t^{2}\right)\right](u^{2}+|\nabla u|^{2})dx\nonumber\\
& &
+\int_{\Omega}\left[ t_{n}^{2}\phi'\left( \frac{u^{2}+|\nabla u|^{2}}{2}t_{n}^{2}\right)- t^{2}\phi'\left( \frac{u^{2}+|\nabla u|^{2}}{2}t^{2}\right)\right](u^{2}+|\nabla u|^{2})^{2}dx\nonumber\\
& &
+(t^{q-1}-t_{n}^{q-1})\lambda q\int_{\Omega} a(x)|u|^{q+1}dx
+(t^{p-1}-t_{n}^{p-1})p\int_{\Omega} b(x)|u|^{p+1}dx|\nonumber\\
&\leq &
\int_{\Omega} \left|\phi\left( \frac{u^{2}+|\nabla u|^{2}}{2}t_{n}^{2}\right)-\phi\left( \frac{u^{2}+|\nabla u|^{2}}{2}t^{2}\right)\right|(u^{2}+|\nabla u|^{2})dx\nonumber\\
& &
+\int_{\Omega} \left| t_{n}^{2}\phi'\left( \frac{u^{2}+|\nabla u|^{2}}{2}t_{n}^{2}\right)- t^{2}\phi'\left( \frac{u^{2}+|\nabla u|^{2}}{2}t^{2}\right)\right|(u^{2}+|\nabla u|^{2})^{2}dx\nonumber\\
& &
+|t^{q-1}-t_{n}^{q-1}|\lambda q\|a\|_{\infty}S^{q+1}_{q+1}\|\nabla u\|_{2}^{q+1}
+|t^{p-1}-t_{n}^{p-1}|p\|b\|_{\infty}S^{p+1}_{p+1}\|\nabla u\|_{2}^{p+1}.
 \end{eqnarray}
It follows from $(\phi_{1})$, $(\phi_{2})$ and Proposition 2.1 that
\begin{eqnarray}\label{534}
\left|\phi\left( \frac{u^{2}+|\nabla u|^{2}}{2}t_{n}^{2}\right)-\phi\left( \frac{u^{2}+|\nabla u|^{2}}{2}t^{2}\right)\right|(u^{2}+|\nabla u|^{2})
\leq
2\rho_{1}(u^{2}+|\nabla u|^{2})
\in L^{1}(\Omega)
 \end{eqnarray}
  and
\begin{eqnarray}\label{535}
&  &
\left| t_{n}^{2}\phi'\left( \frac{u^{2}+|\nabla u|^{2}}{2}t_{n}^{2}\right)-
t^{2}\phi'\left( \frac{u^{2}+|\nabla u|^{2}}{2}t^{2}\right)\right|(u^{2}+|\nabla u|^{2})^{2}\nonumber\\
&=  &
2\left|\phi'\left( \frac{u^{2}+|\nabla u|^{2}}{2}t_{n}^{2}\right)\frac{u^{2}+|\nabla u|^{2}}{2}t_{n}^{2}
- \phi'\left( \frac{u^{2}+|\nabla u|^{2}}{2}t^{2}\right) \frac{u^{2}+|\nabla u|^{2}}{2}t^{2}\right|(u^{2}+|\nabla u|^{2})\nonumber\\
&\leq &
2\left(\left|\phi'\left( \frac{u^{2}+|\nabla u|^{2}}{2}t_{n}^{2}\right)\right|\frac{u^{2}+|\nabla u|^{2}}{2}t_{n}^{2}
+\left|\phi'\left( \frac{u^{2}+|\nabla u|^{2}}{2}t^{2}\right)\right|\frac{u^{2}+|\nabla u|^{2}}{2}t^{2}\right)(u^{2}+|\nabla u|^{2})\nonumber\\
&\leq &
4\rho_{2}(u^{2}+|\nabla u|^{2})
\in L^{1}(\Omega).
 \end{eqnarray}
Then by the continuity of $\phi'$, we also obtain that
\begin{eqnarray}\label{536}
 \left|\phi\left( \frac{u^{2}+|\nabla u|^{2}}{2}t_{n}^{2}\right)
 -\phi\left( \frac{u^{2}+|\nabla u|^{2}}{2}t^{2}\right)\right|(u^{2}+|\nabla u|^{2})
 \rightarrow0,  &\mbox{a.e.}& x \in \Omega
 \end{eqnarray}
and
\begin{eqnarray}\label{537}
 \left| t_{n}^{2}\phi'\left( \frac{u^{2}+|\nabla u|^{2}}{2}t_{n}^{2}\right)
 - t^{2}\phi'\left( \frac{u^{2}+|\nabla u|^{2}}{2}t^{2}\right)\right|(u^{2}+|\nabla u|^{2})^{2}
 \rightarrow0,  &\mbox{a.e.}& x \in \Omega,
 \end{eqnarray}
Hence, it follows from (\ref{534}), (\ref{535}), (\ref{536}), (\ref{537}) and the Lebesgue dominated convergence theorem that
\begin{eqnarray}\label{538}
\lim_{n\rightarrow\infty}\int_{\Omega}\left|\phi\left( \frac{u^{2}+|\nabla u|^{2}}{2}t_{n}^{2}\right)
-\phi\left( \frac{u^{2}+|\nabla u|^{2}}{2}t^{2}\right)\right|(u^{2}+|\nabla u|^{2})=0,
 \end{eqnarray}
and
\begin{eqnarray}\label{539}
\lim_{n\rightarrow\infty}\int_{\Omega}\left| t_{n}^{2}\phi'\left( \frac{u^{2}+|\nabla u|^{2}}{2}t_{n}^{2}\right)
- t^{2}\phi'\left( \frac{u^{2}+|\nabla u|^{2}}{2}t^{2}\right)\right|(u^{2}+|\nabla u|^{2})^{2}=0.
 \end{eqnarray}
At this moment, by (\ref{532}), (\ref{533}), (\ref{538}) and (\ref{539}), we conclude that $\gamma_{u}''$ is continuous.
\qed

\vskip2mm
\noindent{\bf Acknowledgments}\\
\noindent
This project is  supported by Yunnan Fundamental Research Projects (grant No: 202301AT070465) and  Yunnan Ten Thousand Talents Plan Young \& Elite Talents Project.

\vskip2mm
 \noindent
{\bf Statements and Declarations}\\
\noindent
The authors state no conflict of interest.

\renewcommand\refname{References}

\end{document}